\documentclass[12pt]{amsart}
\usepackage{amssymb}
\usepackage{amsmath}
\usepackage[active]{srcltx}
\usepackage{t1enc}
\usepackage[latin2]{inputenc}
\usepackage{verbatim}
\usepackage{amsmath,amsfonts,amssymb,amsthm}
\usepackage[mathcal]{eucal}
\usepackage{enumerate}
\usepackage[centertags]{amsmath}
\usepackage{graphics}

\newtheorem{theorem}{Theorem}

\newtheorem{lemma}{Lemma}
\newtheorem{remark}{Remark}

\newtheorem{corollary}{Corollary}

\newtheorem{example}{Example}

\begin{document}
\author{George Tephnadze}
\title[Partial Sums and Fej\'er Means ]{The One-dimensional Martingale Hardy Spaces  and   Partial Sums and Fej\'er Means with respect to Walsh system}
\address{G. Tephnadze, The University of Georgia, School of science and technology, 77a Merab Kostava St, Tbilisi 0128, Georgia.}
\email{g.tephnadze@ug.edu.ge}

\thanks{The research was supported by Shota Rustaveli National Science	Foundation grant no. FR-19-676.}
\date{}
\maketitle
	
\begin{abstract}
In this paper we prove and discuss some new 
$\left( H_p,L_{p}\right)$ type inequalities for partial Sums and Fej\'er means with respect to Walsh system. It is also proved that these results are the best possible in a special sense.  As applications, both some well-known and new results are pointed out. 
\end{abstract}
	
	
	
\section{Preliminaries}

\text{ \  \ }\text{ \  \ } It is well-known that (for details see e.g. \cite{gol} and \cite{sws}) for every $p>1$ there exists an absolute constant $ c_{p} $, depending only on $p$, such that
\begin{equation*}
\left\Vert S_{n}f\right\Vert _{p}\leq c_{p}\left\Vert f\right\Vert _{p},
\text{ \ when \ }p>1 \text{ \ and \ } f\in H_{1}({G}).
\end{equation*}

Moreover,  Watari \cite{Wat1} (see also Gosselin \cite{goles} and Young \cite{Yo}) proved that there exists an absolute constant $ c $ such that, for $ n =1,2,..., $
\begin{eqnarray*}
	\lambda \mu
	\left( \vert S_n f\vert>\lambda \right)&\leq& c\left\Vert f\right\Vert_{1}, \ \ \ f\in L_1(G_m), \ \ \lambda>0.
\end{eqnarray*}

On the other hand, it is also well-known that (for details see e.g. \cite{AVD} and \cite{sws}) Walsh system is not Schauder basis in $L_{1}({G})$ space. Moreover, there exists function $f\in H_{1}({G})$, such that partial sums with respect to Walsh system are not uniforml5y bounded in $L_{1}({G})$.

By applying Lebesgue constants
\begin{equation*}
L(n):=\left\Vert D_{n}\right\Vert _{1}
\end{equation*}
we easily obtain that (for details see e.g. \cite{1} and \cite{sws}) subsequences of partial sums $S_{n_{k}}f$ with respect to Walsh system converge to $f$ in $L_{1}$ norm if and only if
\begin{equation} \label{1.0.1}
\sup_{k\in \mathbb{N}}L(n_{k})\leq c<\infty.
\end{equation}

Since $ n $-th Lebesgue constant with respect to Walsh system, where  $$n=\sum_{j=0}^{\infty }n_{j}2^{j}, (n_{j}\in Z_{2})$$  can be estimated by  variation of natural number 
\begin{equation*}
V\left( n\right) =n_{0}+\overset{\infty }{\underset{k=1}{\sum }}\left|
n_{k}-n_{k-1}\right|
\end{equation*}
and it is also well known that  (for details see e.g. \cite{BPT1} and \cite{sws}) the following two-sided estimate is true
\begin{equation*}
\frac{1}{8}V\left( n\right) \leq L(n) \leq
V\left( n\right)
\end{equation*}
to obtain convergence of subsequences of partial sums $S_{n_{k}}f$ with respect to Walsh system of $f\in L_{1}$ in  $f\in L_{1}$-norm.
Condition (\ref{1.0.1}) can be replaced by
\begin{equation*}
\sup_{k\in \mathbb{N}}V(n_{k})\leq c<\infty
\end{equation*}

It follows that (for details see e.g. \ \cite{sws} and \cite{We1}) subsequence of partial sums $S_{2^{n}}$ are bounded from $H_{p}(G)$ to $H_{p}(G)$ for every $p>0$, from which we obtain that
\begin{equation} \label{1.S2n000}
\left\Vert S_{2^{n}}f-f\right\Vert _{H_{p}(G)}\rightarrow 0,\text{ \ as \ }n\rightarrow \infty,
\end{equation}
On the other hand, (see e.g.  \cite{tep7}) there exist a martingale $f\in H_{p}(G)$ $\left( 0<p<1\right),$ such that
\begin{equation*}
\underset{n\in\mathbb{N}}{\sup }\left\Vert S_{2^{n}+1}f\right\Vert
_{weak-L_{p}(G)}=\infty.
\end{equation*}%

The main reason of divergence of subsequence  $S_{{2^n}+1}f$ of partial sums it that  (for details see \cite{tep9}) Fourier coefficients of $f\in H_{p}(G)$ are not uniformly bounded when $0<p<1$.

When $0<p<1$ in \cite{tep_thesis} was investigated boundedness of  subsequences of partial sums with respect to Walsh system from $H_{p}(G)$ to $H_{p}(G)$. In particular, the following result is true:

\textbf{Theorem T1.} Let $0<p<1$ and $f\in H_{p}(G)$. Then there exists a absolute constant $c_{p},$ depending only on $p$, such that
\begin{equation*}
\left\Vert S_{m_{k}}f\right\Vert _{H_{p}(G)}\leq c_{p}\left\Vert f\right\Vert _{H_{p}(G)}
\end{equation*}
if and only if the following condition holds
\begin{equation} \label{1.cond001}
\underset{k\in \mathbb{N}}{\sup } d\left( m_{k}\right)<c<\infty,
\end{equation}
where
\begin{equation*}
d\left( m_{k}\right) :=\left\vert m_{k}\right\vert -\left\langle m_{k}\right\rangle.
\end{equation*}

In particular, Theorem T1 immediately follows:

\textbf{Theorem T2.} Let $p>0$ and $f\in H_{p}(G)$. Then there exists a absolute constant $c_{p},$ depending only on $p$, such that
\begin{equation*}
\left\Vert S_{2^{n}}f\right\Vert_{H_p(G)}\leq c_{p}\left\Vert f\right\Vert _{H_p(G)}
\end{equation*}
and
\begin{equation*}
\left\Vert S_{2^{n}+2^{n-1}}f\right\Vert_{H_p(G)}\leq c_{p}\left\Vert f\right\Vert _{H_p(G)}.
\end{equation*}

On the other hand, we have the following result:

\textbf{Theorem T3.}  Let $p>0$. Then there exists a martingale $ f \in H_{p}(G) $, such that
$$  \sup_{n\in \mathbb{N}}\left\Vert S_{2^{n}+1}f\right\Vert_{_{H_{p}(G)}}=\infty. $$

Taking into account these results it is interesting to find behaviour of rate of divergence of subsequences of partials sums with respect to Walsh system of martingale $f\in H_{p}(G)$ in the martingale Hardy spaces $H_{p}(G)$.

In the second chapter of this thesis (see also \cite{tep12}) we investigate above mentioned problem. For $0<p<1$ we have the following result:

\textbf{Theorem \ref{th4.1.1}.} Let $f\in H_{p}(G)$. Then there exists a absolute constant $c_{p},$ depending only on $p$, such that the following inequality is true
\begin{equation} \label{1.vnk00002}
\text{ }\left\Vert S_{n}f\right\Vert _{H_{p}(G)}\leq c_{p}2^{d\left( n\right)\left( 1/p-1\right) }\left\Vert f\right\Vert _{H_{p}(G)}.
\end{equation}

On the other hand, if  $0<p<1,$ $\left\{ m_{k}:\text{ }k\geq 0\right\} $ be increasing subsequence of natural numbers, such that
\begin{equation}  \label{1.dmk}
\sup_{k\in \mathbb{N}}d\left( m_{k}\right) =\infty
\end{equation}%
and
$\Phi :\mathbb{N}_{+}\rightarrow \lbrack 1,\infty )$ be non-decreasing function satisfying the condition
\begin{equation*} \label{1010}
\overline{\underset{k\rightarrow \infty }{\lim }}\frac{2^{d\left(
		m_{k}\right) \left( 1/p-1\right) }}{\Phi \left( m_{k}\right) }=\infty,
\end{equation*}%
then there exists a martingale $f\in H_{p}(G),$ such that
\begin{equation*}
\underset{k\in \mathbb{N}}{\sup }\left\Vert \frac{S_{m_{k}}f}{\Phi \left(m_{k}\right) }\right\Vert _{weak-L_{p}(G)}=\infty .
\end{equation*}

Theorem \ref{cor4.1.1} easily follows the following corollary:

\textbf{Corollary \ref{cor4.1.1}.} Let $0<p<1$ and $f\in H_{p}(G)$. Then there exists a absolute constant $c_{p},$ depending only on $p$, such that
\begin{equation*}
\text{ }\left\Vert S_{n}f\right\Vert _{H_{p}(G)}\leq c_{p}\left( n\mu \left\{\text{supp}\left( D_{n}\right) \right\} \right) ^{1/p-1}\left\Vert f\right\Vert _{H_{p}(G)}.
\end{equation*}

On the other hand, if $0<p<1$ and $\left\{ m_{k}:\text{ }k\geq 0\right\} $ be increasing sequence of natural numbers, such that
\begin{equation*}
\sup_{k\in \mathbb{N}}m_{k}\mu \left\{ \text{supp}\left( D_{m_{k}}\right)
\right\} =\infty
\end{equation*}%
and $\Phi :\mathbb{N}_{+}\rightarrow \lbrack 1,\infty )$
be non-decreasing function satisfying the condition
\begin{equation*}  \label{12e}
\overline{\underset{k\rightarrow \infty }{\lim }}\frac{\left( m_{k}\mu\left\{ \text{supp}\left( D_{m_{k}}\right) \right\} \right) ^{1/p-1}}{\Phi
	\left( m_{k}\right) }=\infty,
\end{equation*}%
then there exists a martingale $f\in H_{p}(G),$  such that
\begin{equation*}
\underset{k\in \mathbb{N}}{\sup }\left\Vert \frac{S_{m_{k}}f}{\Phi \left(
	m_{k}\right) }\right\Vert _{weak-L_{p}(G)}=\infty .
\end{equation*}

In particular, we also get the proofs of  Theorem T1 and Theorem T2.

In the second chapter of this thesis we also investigate case $p=1$. In this case the following result is true:

\textbf{Theorem \ref{th4.1.2}.} Let $n\in \mathbb{N}_{+}$ and $f\in H_{1}(G).$ Then there exists a absolute constant $c,$ such that
\begin{equation} \label{1.vnk00001}
\left\Vert S_{n}f\right\Vert _{H_{1}(G)}\leq cV\left( n\right) \left\Vert
f\right\Vert _{H_{1}(G)}.
\end{equation}

Moreover, if $\left\{ m_{k}:\text{ }k\geq 0\right\} $ be increasing sequence of natural numbers $\mathbb{N}_{+},$ such that
\begin{equation*}
\sup_{k\in \mathbb{N}}V\left( m_{k}\right) =\infty  \label{vnk}
\end{equation*}
and $\Phi :\mathbb{N}_{+}\rightarrow \lbrack 1,\infty )$
be non-decreasing function satisfying the condition
\begin{equation*}
\overline{\underset{k\rightarrow \infty }{\lim }}\frac{V\left( m_{k}\right)}{\Phi \left( m_{k}\right) }=\infty .
\end{equation*}
Then there exists a martingale $f\in H_{1}(G),$ such that
\begin{equation*}
\underset{{k\in \mathbb{N}}}{\sup }\left\Vert \frac{S_{m_{k}}f}{\Phi \left(m_{k}\right) }\right\Vert _{1}=\infty .
\end{equation*}

When $0<p<1$ in \cite{tep_thesis} was proved boundedness of maximal operators of subsequences of partial sums from $H_{p}(G)$ to $L_{p}(G)$. In particular, the following is true:

\textbf{Theorem T4.} Let $0<p<1$ and $f\in H_{p}(G)$. Then the maximal operator
\begin{equation*}
\underset{k\in \mathbb{N}}{\sup }|S_{m_{k}}f| \label{cond00}
\end{equation*}
is bounded from $H_{p}(G)$ to $L_{p}(G)$, if and only if condition (\ref{1.cond001}) is fulfilled.

In the special cases we obtain that the following is true:

\textbf{Theorem T5.} Let $p>0$ and $f\in H_{p}(G)$. Then there exists an absolute constant $c_{p},$ depending only on $p$, such that
\begin{equation}\label{1.S2n}
\left\Vert  \sup_{n\in \mathbb{N}}|S_{2^{n}}f|\right\Vert_{p}\leq c_{p}\left\Vert f\right\Vert _{H_p(G)}
\end{equation}
and
\begin{equation*}
\left\Vert  \sup_{n\in \mathbb{N}}|S_{2^{n}+2^{n-1}}f|\right\Vert_{p}\leq c_{p}\left\Vert f\right\Vert _{H_p(G)}.
\end{equation*}

On the other hand we have the following result:

\textbf{Theorem T6.}  Let $p>0$. Then there exists a martingale $ f \in H_{p}(G) $, such that
$$ \left\Vert  \sup_{n\in \mathbb{N}}|S_{2^{n}+1}f|\right\Vert_{p}=\infty. $$

Above mentioned condition (\ref{1.cond001}) is sufficient condition for the case $p=1$ also, but there exist subsequences which do not satisfy this condition, but maximal operators of these subsequences of partial sums with respect to Walsh system  are not bounded from  $H_{1}(G)$ to  $L_{1}(G)$.

Such necessary and sufficient conditions which provides boundedness of maximal operators of subsequences of partial sums with respect to Walsh system from $H_{1}(G)$ to  $L_{1}(G)$ is open problem.

In \cite{tep9} and \cite{tep_thesis} was investigated boundedness of weighted maximal operators from $H_{p}(G)$ to $L_{p}(G)$, when $ 0<p\leq 1$:

\textbf{Theorem T7.} Let $0<p\leq 1$. Then weighted maximal operator
\begin{equation*}
\overset{\sim }{S}_{p}^{\ast }f:=\underset{n\in \mathbb{N_+}}{%
	\sup }\frac{\left\vert S_{n}f\right\vert }{\left( n+1\right) ^{1/p-1}\log^{[p]}\left( n+1\right)}
\end{equation*}
is bounded from  $H_{p}(G)$ to  $L_{p}(G)$, where $ [p] $ denotes integer part of $ p $.

Moreover, for any non-decreasing function  $\varphi :\mathbb{N}_{+}\rightarrow [1,$ $\infty )$ satisfying the condition
\begin{equation*}
\overline{\lim_{n\rightarrow \infty }}\frac{\left( n+1\right) ^{1/p-1}\log^{[p]}\left( n+1\right)}{%
	\varphi \left( n+1\right)}=+\infty,
\end{equation*}%
there exists a martingale $f \in H_{p}(G)$ $ (0<p\leq 1), $ such that
\begin{equation*}
\sup_{n\in \mathbb{N}}\left\Vert \frac{S_{n}f}{\varphi \left( n\right)
}\right\Vert _{p}=\infty.
\end{equation*}

According to negative result for weighted maximal operator of partial sums of Walsh-Fourier series we immediately get the following result:

\textbf{Theorem S1.}  There exists a martingale $f\in H_{p}(G)$, $(0<p\leq1)$, such that
\begin{equation*}
\underset{n\in \mathbb{N}}{\sup }\left\Vert S_{n}f\right\Vert
_{p}=\infty .
\end{equation*}%

On the other hand, boundedness of weighted maximal operators immediately follows the following estimation:

\textbf{Theorem S2.} Let $0<p\leq1$. Then there exists a absolute constant $c_{p}$, depending only on $p$, such that 
\begin{equation*}
\left\Vert S_{n}f\right\Vert _{p}\leq c_{p} {(n+1)^{1/p-1}\log^{[p]} (n+1)}\left\Vert f\right\Vert _{H_p(G)},
\text{ \ for \ } 0<p\leq 1,
\end{equation*}
where $ [p] $ denotes integer part of $ p $.

By applying this inequality (see \cite{tep6}) we find necessary and sufficient conditions for martingale $f\in H_{p}(G)$ for which partial sums with respect to Walsh system of martingale $f\in H_{p}(G)$ converge in $H_{p}(G)$ norm.

\textbf{Theorem T8.} Let $0<p\leq 1$,  $[p]$ denotes integer part of $p$, $f\in H_{p}(G)$ and
\begin{equation*}
\omega _{H_{p}(G)}\left( \frac{1}{2^{N}},f\right)=o\left( \frac{1}{2^{N(1/p-1)}N^{[p]}}\right) ,\text{ \ as \ }N\rightarrow
\infty.
\end{equation*}%
Then
\begin{equation*}
\left\Vert S_{n}f-f\right\Vert _{p}\rightarrow 0,\text{ as }
n\rightarrow \infty.
\end{equation*}

Moreover, there exists a martingale $f\in H_{p}(G)$, where $0<p<1$, such that
\begin{equation*}
\omega _{H_{p}(G)} \left( \frac{1}{2^{N}},f\right)=O\left( \frac{1}{2^{N(1/p-1)}N^{[p]}}\right) ,\text{ \ as \ }N\rightarrow \infty
\end{equation*}%
and
\begin{equation*}
\left\Vert S_{n}f-f\right\Vert _{weak-L_{p}(G)}\nrightarrow 0,\,\,\,\text{ \ \ as \ \ }n\rightarrow \infty .
\end{equation*}%

By taking these results into account, it is interesting to find necessary and sufficient conditions for modulus of continuity, such that subsequences of partial sums with respect to Walsh system of martingale $f\in H_{p}(G)$ converge in $H_{p}(G)$ norm.

In the second chapter of this thesis (see also \cite{tep12}) we investigate this problem. By combining inequalities (\ref{1.vnk00002}) and (\ref{1.vnk00001}) we get the following theorem:

\textbf{Theorem \ref{theorem4.2.1}.}  Let $2^{k}<n\leq 2^{k+1}.$ Then there exists an absolute constant $c_{p},$ depending only on $p$, such that
\begin{equation} \label{1.sn3}
\left\Vert S_{n}f-f\right\Vert _{H_{p}(G)}\leq c_{p}2^{d\left( n\right) \left(1/p-1\right) }\omega _{H_{p}(G)}\left( \frac{1}{2^{k}},f\right) ,\text{ \ \ \ }\left( 0<p<1\right)
\end{equation}
and
\begin{equation}  \label{1.sn2}
\left\Vert S_{n}f-f\right\Vert _{H_{1}(G)}\leq c_{1}V\left( n\right) \omega_{H_{1}(G)}\left( \frac{1}{2^{k}},f\right)
\end{equation}

By applying inequality (\ref{1.sn3}) the following result is proved in the second chapter:

\textbf{Theorem \ref{th4.2.2}.} Let $0<p<1,$ $f\in H_{p}(G)$ and $\{m_{k}:k\geq 0\}$ be increasing sequence of natural number satisfying the condition
\begin{equation*}
\omega _{H_{p}(G)}\left( \frac{1}{2^{\left\vert m_{k}\right\vert }},f\right)
=o\left( \frac{1}{2^{d\left( m_{k}\right) \left( 1/p-1\right) }}\right)
\text{ as \ }k\rightarrow \infty.  \label{1.18a}
\end{equation*}%
Then
\begin{equation} \label{1.con1}
\left\Vert S_{m_{k}}f-f\right\Vert _{H_{p}(G)}\rightarrow 0\text{ \ as \ }k\rightarrow \infty .
\end{equation}

On the other hand, if $\{m_{k}:k\geq 0\}$ be increasing sequence of natural numbers satisfying the condition (\ref{1.dmk}), then there exists a martingale $f\in H_{p}(G)$ and subsequence $\{\alpha _{k}:k\geq 0\}\subset \{m_{k}:k\geq 0\},$ for which
\begin{equation*}
\omega _{H_{p}(G)}\left( \frac{1}{2^{\left\vert \alpha _{k}\right\vert }},f\right) =O\left( \frac{1}{2^{d\left( \alpha _{k}\right) \left(
		1/p-1\right) }}\right) \text{ as \ }k\rightarrow \infty \text{\ }
\end{equation*}%
and
\begin{equation}
\limsup\limits_{k\rightarrow \infty }\left\Vert S_{\alpha_{k}}f-f\right\Vert _{weak-L_{p}(G)}>c_{p}>0\,\,\,\text{as\thinspace
	\thinspace \thinspace }k\rightarrow \infty,  \label{1.con11}
\end{equation}
where $c_{p}$ is an absolute constant depending only on $p$.

According to this theorem we immediately get that the following result is true:

\textbf{Corollary \ref{cor4.2.1}.} Let $0<p<1,$ $f\in H_{p}(G)$ and $\{m_{k}:k\geq 0\}$  be increasing sequence of natural number, satisfying the condition
\begin{equation*} \label{1.cond2}
\omega _{H_{p}(G)}\left( \frac{1}{2^{\left\vert m_{k}\right\vert }},f\right)=o\left( \frac{1}{\left( m_{k}\mu \left( \text{supp}D_{m_{k}}\right) \right)^{1/p-1}}\right), \text{ \ as \ }k\rightarrow \infty.
\end{equation*}%
Then (\ref{1.con1}) holds.

On the other hand, if $\{m_{k}:k\geq 0\}$ be increasing sequence of natural number, satisfying the condition
\begin{equation*}
\overline{\underset{k\rightarrow \infty}{\lim }}\frac{\left( m_{k}\mu
	\left\{\text{supp}\left(D_{m_{k}}\right)\right\}\right)^{1/p-1}}{\Phi
	\left(m_{k}\right)}=\infty,
\end{equation*}%
then there exists a martingale $f\in H_{p}(G)$ and subsequence $\{\alpha _{k}:k\geq 0\}\subset \{m_{k}:k\geq 0\}$ such that
\begin{equation*}
\omega _{H_{p}(G)}\left( \frac{1}{2^{\left\vert \alpha _{k}\right\vert }},f\right) =O\left( \frac{1}{\left( \alpha _{k}\mu \left( \text{supp}%
	D_{\alpha _{k}}\right) \right) ^{1/p-1}}\right), \text{\ \ as \ }k\rightarrow\infty
\end{equation*}%
and (\ref{1.con11}) holds.

By applying (\ref{1.sn2}) we prove that the following is true:

\textbf{Theorem \ref{th4.2.3}.} Let $f\in H_{1}(G)$ and $\{m_{k}:k\geq 0\}$ be increasing sequence of natural number, satisfying the condition
\begin{equation*}
\omega _{H_{1}(G)}\left( \frac{1}{2^{\left\vert m_{k}\right\vert }},f\right)=o\left( \frac{1}{V\left( m_{k}\right) }\right) \text{ \ as \ }k\rightarrow
\infty.
\end{equation*}%
Then
\begin{equation*}
\left\Vert S_{m_{k}}f-f\right\Vert _{H_{1}(G)}\rightarrow 0\text{ \ as \ }
k\rightarrow \infty.
\end{equation*}

Moreover, if $\{m_{k}:k\geq 0\}$ be increasing sequence of natural number, satisfying the condition (\ref{1.dmk}), then there exists a martingale $f\in H_{1}(G)$ and subsequence $\{\alpha _{k}:k\geq 0\}\subset \{m_{k}:k\geq 0\}$ for which
\begin{equation*}
\omega _{H_{1}(G)}\left( \frac{1}{2^{\left\vert \alpha _{k}\right\vert }}%
,f\right) =O\left( \frac{1}{V\left( \alpha _{k}\right) }\right) \text{ \ as
	\ }k\rightarrow \infty
\end{equation*}%
and
\begin{equation*} \label{1.cond10}
\limsup\limits_{k\rightarrow \infty}\left\Vert S_{\alpha_k} f-f\right\Vert_1>c>0\text{ \ \ as \ \  }
k\rightarrow \infty, 
\end{equation*}%
where  $c$ is an absolute constant.

By applying Theorem \ref{th4.2.2} and Theorem \ref{th4.2.3} we immediately get proof of Theorem T8.

Weisz \cite{We3} consider convergence in norm of Fej\'er means of the one-dimensional Walsh-Fourier and proved the following:

\textbf{Theorem We1.} Let $p>1/2$ and $f\in H_{p}(G)$. Then there exists a absolute constant $c_{p}$, depending only on $p$, such that 
\begin{equation*}
\left\Vert \sigma _{k}f\right\Vert _{H_{p}(G)}\leq c_{p}\left\Vert f\right\Vert_{H_{p}(G)}.
\end{equation*}

Weisz (for details see e.g. \cite{We1}) also consider boundedness of subsequences of Fej\'er means $\sigma_{2^{n}}$ of the one-dimensional Walsh-Fourier series  from $H_{p}(G)$ to $H_{p}(G)$ when $p>0$:

\textbf{Theorem We2.} Let $p>0$ and $f\in H_{p}(G)$. Then
\begin{equation} \label{fe22222}
\left\Vert \sigma_{2^{k}}f-f\right\Vert _{H_{p}(G)}\rightarrow 0,\text{\ \ as \ \ }k\rightarrow\infty.
\end{equation}

On the other hand, in \cite{tep3} was proved the following result:

\textbf{Theorem T9.} There exists a martingale $f\in H_{p}(G)$ $\left( 0<p\leq 1/2\right)$ such that
\begin{equation*}
\underset{n\in\mathbb{N}}{\sup}\left\Vert\sigma_{2^{n}+1}f\right\Vert_{H_p(G)}=\infty.
\end{equation*}

Goginava \cite{gog16} (see also \cite{pt1}) proved that the following result is true:

\textbf{Theorem Gog1.} Let $0<p\leq 1.$ Then the sequence of operators $\left\vert \sigma _{2^{n}}f\right\vert $ are not bounded from $H_{p}(G)$ to $H_{p}(G)$.

When $ 0<p<1/2 $ then in \cite{pt2} was proved bondedness of subsequences of Fej\'er means of the one-dimensional Walsh-Fourier from $H_{p}(G)$ to $H_{p}(G)$. In particular, the following is true:

\textbf{Theorem T10.} Let $0<p<1/2$ and $f\in H_{p}(G)$. Then there exists a absolute constant $c_{p}$, depending only on $p$, such that 
\begin{equation*}
\left\Vert \sigma_{m_{k}}f\right\Vert _{H_{p}(G)}\leq c_{p}\left\Vert f\right\Vert _{H_{p}(G)}
\end{equation*}
estimation holds if and only if the condition (\ref{1.cond001}) is fulfilled.

Theorem T10 immediately follows theorem of Weisz (see Theorem We2) and and also interesting results:

\textbf{Theorem T11.} Let $p>0$ and $f\in H_{p}(G)$. Then there exists an absolute constant $c_{p},$ depending only on $p$, such that
\begin{equation*}
\left\Vert \sigma_{2^{n}}f\right\Vert_{H_p(G)}\leq c_{p}\left\Vert f\right\Vert _{H_p(G)}
\end{equation*}
and
\begin{equation*}
\left\Vert \sigma_{2^{n}+2^{n-1}}f\right\Vert_{H_p(G)}\leq c_{p}\left\Vert f\right\Vert _{H_p(G)}.
\end{equation*}

On the other hand, we have the following result:

\textbf{Theorem T12.}  Let $p>0$. Then there exists a martingale $ f \in H_{p}(G) $, such that
$$  \sup_{n\in \mathbb{N}}\left\Vert \sigma_{2^{n}+1}f\right\Vert_{_{H_{p}(G)}}=\infty. $$

According to above mentioned results it is interesting to find rate of divergence of subsequences  $\sigma_{n_{k}}f$ of Fej\'er means of the one-dimensional Walsh-Fourier series in the Hardy spaces $H_{p}(G)$.

In the third chapter of this thesis (see also \cite{tep13}) we find rate of divergence of subsequences of Fej\'er means of the one-dimensional Walsh-Fourier series on the martingale Hardy spaces $H_{p}(G)$, when
$ 0<p\leq1/2 $.

First, we consider case $p=1/2$:

\textbf{Theorem \ref{th5.1.1}.} Let $n\in \mathbb{N}_{+}$ and $f\in H_{1/2}(G).$ Then there exists an absolute constant $c,$ such that
\begin{equation} \label{vnk000f}
\left\Vert \sigma_{n}f\right\Vert _{H_{1/2}(G)}\leq c{{V}^2}\left( n\right) \left\Vert f\right\Vert _{H_{1/2}(G)}.
\end{equation}

Moreover, if $\left\{ m_{k}:\text{ }k\geq 0\right\} $ be increasing secuence of natural numbers, such that
\begin{equation*}
\sup_{k\in \mathbb{N}}V\left( m_{k}\right) =\infty  \label{vnkf}
\end{equation*}
and $\Phi :\mathbb{N}_{+}\rightarrow [1,\infty]$
be non-decreasing function satisfying the conditions
\begin{equation*}
\overline{\underset{k\rightarrow \infty }{\lim }}\frac{{V^2}\left( m_{k}\right)}{\Phi \left( m_{k}\right) }=\infty .  \label{17aaf}
\end{equation*}
then there exists a martingale $f\in H_{1/2}(G),$ such that
\begin{equation*}
\underset{{k\in \mathbb{N}}}{\sup }\left\Vert \frac{\sigma_{m_{k}}f}{\Phi \left(m_{k}\right) }\right\Vert _{1/2}=\infty.
\end{equation*}

There was also considered case $0<p<1/2$ and was proved that the following is true:

\textbf{Theorem \ref{th5.1.2}.} Let $0<p<1/2$ and $f\in H_{p}(G)$. Then there exists an absolute constant $c_{p},$ depending only on $p$ such that
\begin{equation} \label{vnk0000f}
\text{ }\left\Vert \sigma_{n}f\right\Vert _{H_{p}(G)}\leq c_{p}2^{d\left( n\right)\left( 1/p-2\right) }\left\Vert f\right\Vert _{H_{p}(G)}.
\end{equation}

On the other hand, if  $0<p<1/2,$ $\left\{ m_{k}:\text{ }k\geq 0\right\} $ be increasing sequence of natural numbers satisfying the condition (\ref{1.dmk})  and
$\Phi :\mathbb{N}_{+}\rightarrow \lbrack 1,\infty )$ be non-decreasing function such that
\begin{equation*}
\overline{\underset{k\rightarrow \infty }{\lim }}\frac{2^{d\left(m_{k}\right) \left( 1/p-2\right) }}{\Phi \left( m_{k}\right) }=\infty,
\end{equation*}%
then there exists a martingale $f\in H_{p}(G),$ such that
\begin{equation*}
\underset{k\in \mathbb{N}}{\sup }\left\Vert \frac{\sigma_{m_{k}}f}{\Phi \left(m_{k}\right) }\right\Vert _{weak-L_{p}(G)}=\infty .
\end{equation*}

From these results also follows proof of Theorem We2.

In 1975 Schipp \cite{Sc} (see also \cite{1} and \cite{zy}) proved that the maximal operator of  Fej\'er means $\sigma ^{*}$ is of type weak-(1,1):
\begin{equation*}
\mu \left( \sigma ^{*}f>\lambda \right) \leq \frac{c}{\lambda }\left\|f\right\| _{1},\text{ \qquad }\left( \lambda >0\right).
\end{equation*}
By using  Marcinkiewicz interpolation theorem it follows that  $\sigma ^{*}$ is of strong type-$ (p,p) $, when $p>1:$
\begin{equation*}
\left\| \sigma ^{*}f \right\|_{p} \leq c\left\|
f\right\| _{p}, \text{ \qquad }\left( p >1\right).
\end{equation*}
The boundedness does not hold for  $ p=1,$ but Fujii \cite{Fu} (see also \cite{Yano}) proved that maximal operator of  Fej\'er means is bounded from  $H_{1}(G)$ to $L_{1}(G)$. Weisz in \cite{We2} generalized result of Fujii and proved that maximal operator of  Fej\'er means is bounded from $H_{p}(G)$ to $L_{p}(G)$, when $ p>1/2. $ Simon \cite{Si1} construct the counterexample, which shows that boundedness does not hold when  $0<p<1/2$. Goginava \cite{GoAMH} (see also \cite{BGG} and \cite{BGG2}) generalized this result for   $0<p\leq1/2$ and proved that the following is true:

\textbf{Theorem Gog2.} There exists a martingale $f\in H_{p}(G)$ $\left( 0<p\leq1/2\right)$ such that
\begin{equation*}
\underset{n\in \mathbb{N}}{\sup }\left\Vert \sigma_{n}f\right\Vert_{p}=\infty .
\end{equation*}%

Weisz \cite{we4} (see also Goginava \cite{GoPubl}) proved that the following is true:

\textbf{Theorem We3.} Let $f\in H_{1/2}(G)$. Then there exists an absolute constant $c,$ such that
\begin{equation*}
\left\|\sigma^{*}f\right\|_{weak-L_{1/2}(G)}\leq c\left\|f\right\|_{H_{1/2}(G)}.
\end{equation*}

In \cite{pt2} was considered boundedness of maximal operators of subsequences of Fej\'er means of the one-dimensional Walsh-Fourier series from $H_{p}(G)$ to $L_{p}(G)$ for  $0<p<1/2$. In particular, the following is true:

\textbf{Theorem T13.} Let $0<p<1/2$ and $f\in H_{p}(G)$. Then the maximal operator 
\begin{equation*}
\overset{\sim }{\sigma}^{*}f:=\underset{k\in \mathbb{N}}{\sup }|\sigma_{m_{k}}f|
\end{equation*}
is bounded from $H_{p}(G)$ to $L_{p}(G)$ if and only if when condition (\ref{1.cond001}) is fulfilled.

As consequences the following results are true:

\textbf{Theorem T14.} Let $p>0$ and $f\in H_{p}(G)$. Then there exists an absolute constant $c_p$ depending only on $ p, $ such that
\begin{equation}\label{sigmamax}
\left\Vert  \sup_{n\in \mathbb{N}}|\sigma_{2^{n}}f|\right\Vert_{p}\leq c_{p}\left\Vert f\right\Vert _{H_p(G)}
\end{equation}
and
\begin{equation*}
\left\Vert  \sup_{n\in \mathbb{N}}|\sigma_{2^{n}+2^{n-1}}f|\right\Vert_{p}\leq c_{p}\left\Vert f\right\Vert _{H_p(G)}.
\end{equation*}

On the other hand, we have the following negative result:

\textbf{Theorem T15.}  Let $0<p<1/2.$ Then there exists a martingale $ f \in H_{p}(G) $, such that
$$ \left\Vert  \sup_{n\in \mathbb{N}}|\sigma_{2^{n}+1}f|\right\Vert_{p}=\infty. $$

above mentioned condition is sufficient for the case $p=1/2$ also, but there exists subsequences, which do not satisfy condition (\ref{1.cond001}), but maximal operator of subsequences of Fej\'er means of the one-dimensional Walsh-Fourier series are bounded from $H_{1/2}(G)$ to $L_{1/2}(G).$

However, it is open problem to find necessary and sufficient conditions on the indexes, which provide boundedness of maximal operator of subsequences of Fej\'er means of the one-dimensional Walsh-Fourier series from $H_{1/2}(G)$ to $L_{1/2}(G)$.

In \cite{GoSzeged} and \cite{tep3} (see also \cite{ptw4}, \cite{tepkack4}, \cite{GNCz} and \cite{tep2}) is proved that the following is true:

\textbf{Theorem GT1.} Let $0<p\leq 1/2$ and $ f \in H_{p}(G) $. Then the maximal operator
\begin{equation*}
\overset{\sim }{\sigma }_{p}^{\ast }f:=\underset{n\in \mathbb{N}}{\sup }\frac{\left\vert \sigma _{n}f\right\vert }{\left( n+1\right) ^{1/p-2}\log ^{2[1/2+p]}\left( n+1\right)}
\end{equation*}
is bounded from  $H_{p}(G)$ to  $L_{p}(G)$.

Moreover, for any nondecreasing function  $\varphi :\mathbb{N}_{+}\rightarrow [1,$ $\infty )$ satisfying the condition
\begin{equation*}
\overline{\lim_{n\rightarrow \infty }}\frac{\left( n+1\right) ^{1/p-2}\log ^{2[1/2+p]}\left( n+1\right)}{\varphi \left( n\right) }=+\infty ,
\end{equation*}%
there exists a martingale $f \in H_{p}(G)$, $ (0<p<1/2) $  such that
\begin{equation*}
\sup_{n\in \mathbb{N}}\left\Vert \frac{\sigma _{n}f}{\varphi \left( n\right)}\right\Vert _{p}=\infty.
\end{equation*}

From the divergence of weighted maximal operators we immediately get that there exists a martingale $f\in H_{p}(G)$ $(0<p\leq 1/2)$, such that
\begin{equation*}
\underset{n\in \mathbb{N}}{\sup }\left\Vert \sigma_{n}f\right\Vert
_{p}=\infty.
\end{equation*}
and from the boundedness results of weighted maximal operators we immediately get that for any $f\in H_{p}(G)$ there exists an absolute constant $c_{p}$, such that the following inequality holds true:
\begin{equation} \label{1.cond33}
\left\Vert \sigma _{n}f\right\Vert _{p}\leq c_{p} {n^{1/p-2}\log ^{2[1/2+p]}\left( n+1\right)}\left\Vert f\right\Vert _{H_{p}(G)},\text{ \ as \ } 0<p\leq1/2.
\end{equation}

By applying inequality (\ref{1.cond33}) in \cite{tep6} was found necessary and sufficient conditions for modulus of continuity of martingale $f\in H_{p}(G)$, for which Fej\'er means of the one-dimensional Walsh-Fourier series converge in $H_{p}(G)$ norm.

\textbf{Theorem T16.} Let $ 0<p\leq 1/2 $, $f\in H_{p}(G)$ and
\begin{equation*}
\omega _{H_{p}(G)}\left( \frac{1}{2^{N}},f\right) =o\left(\frac{1}{2^{N(1/p-2)}N^{2[1/2+p]}}\right),\text{\  as \ } N\rightarrow\infty.
\end{equation*}%
Then
\begin{equation*}
\left\Vert\sigma_{n}f-f\right\Vert _{p}\rightarrow 0,\text{ as }n\rightarrow \infty .
\end{equation*}

Moreover, there exists a martingale $f\in H_{p}(G)$, for which
\begin{equation*}
\omega _{H_{1/2}(G)}\left( \frac{1}{2^{N}},f\right) =O\left( \frac{1}{2^{N(1/p-2)}N^{2[1/2+p]}}\right) ,\text{ \ as \ }N\rightarrow \infty
\end{equation*}%
and
\begin{equation*}
\left\Vert\sigma_{n}f-f\right\Vert_{p}\nrightarrow 0,\,\,\,\text{as\ \ }n\rightarrow \infty.
\end{equation*}

According above mentioned results, it is interesting to find necessary and sufficient conditions for the modulus of continuity, for which subsequences $\sigma_{n_{k}}f$ of Fej\'er means of the one-dimensional Walsh-Fourier series converge in $H_{p}(G)$ norm.

In the third chapter of this thesis we find necessary and sufficient conditions for the modulus of continuity, for which subsequences $\sigma_{n_{k}}f$ of Fej\'er means of the one-dimensional Walsh-Fourier series converge in $H_{p}(G)$ norm  (see also \cite{tep13}).

By applying inequality (\ref{vnk000f}) for the case $p=1/2$ the following necessary and sufficient conditions are found:

\textbf{Theorem \ref{theorem5.2.1}.} Let $f\in H_{1/2}(G)$ and $\{m_{k}:k\geq 0\}$ be increasing sequence of natural numbers, such that
\begin{equation*}
\omega _{H_{1/2}(G)}\left( \frac{1}{2^{\left\vert m_{k}\right\vert }},f\right)=o\left( \frac{1}{{V^2}\left( m_{k}\right) }\right) \text{ \ as \ }k\rightarrow\infty. \label{cond1f}
\end{equation*}%
Then
\begin{equation*}
\left\Vert \sigma_{m_{k}}f-f\right\Vert _{H_{1/2}(G)}\rightarrow 0\text{ \ as \ }k\rightarrow \infty .
\end{equation*}

Moreover, if $\{m_{k}:k\geq 0\}$ be increasing sequence of natural numbers, such that (\ref{1.dmk}) holds true, then there exists a martingale $f\in H_{1/2}(G)$ and subsequence $\{\alpha _{k}:k\geq 0\}\subset \{m_{k}:k\geq 0\}$ such that
\begin{equation*}
\omega _{H_{1/2}(G)}\left( \frac{1}{2^{\left\vert \alpha _{k}\right\vert }},f\right) =O\left( \frac{1}{{V^2}\left( \alpha _{k}\right) }\right) \text{ \ as\ }k\rightarrow \infty
\end{equation*}%
and
\begin{equation*}
\limsup\limits_{k\rightarrow \infty }\left\Vert \sigma_{\alpha_{k}}f-f\right\Vert _{1/2}>c>0\,\,\,\text{as\thinspace \thinspace \thinspace }%
k\rightarrow \infty ,  \label{cond10f}
\end{equation*}
wher  $c$ is an absolute constant.

By applying inequality (\ref{vnk0000f}) we also investigate case $0<p<1/2$. In  the third chapter of this thesis we prove that the following is true:

\textbf{Theorem \ref{theorem5.2.2}.} Let $0<p<1/2,$ $f\in H_{p}(G)$ and $\{m_{k}:k\geq 0\}$ be increasing sequence of natural numbers, such that
\begin{equation*}
\omega _{H_{p}(G)}\left( \frac{1}{2^{\left\vert m_{k}\right\vert }},f\right)=o\left( \frac{1}{2^{d\left( m_{k}\right) \left( 1/p-2\right) }}\right),
\text{ as \ }k\rightarrow \infty.
\end{equation*}%
Then
\begin{equation*} \label{con1f}
\left\Vert \sigma_{m_{k}}f-f\right\Vert _{H_{p}(G)}\rightarrow 0,\text{ \ as \ }k\rightarrow \infty .  
\end{equation*}

On the other hand, if $\{m_{k}:k\geq 0\}$ be increasing sequence of natural numbers satisfying the condition (\ref{1.dmk}), then there exists a martingale $f\in H_{p}(G)$ and subsequence $\{\alpha _{k}:k\geq 0\}\subset \{m_{k}:k\geq 0\},$ for which
\begin{equation*}
\omega _{H_{p}(G)}\left( \frac{1}{2^{\left\vert \alpha _{k}\right\vert }},f\right) =O\left( \frac{1}{2^{d\left( \alpha _{k}\right) \left(
		1/p-2\right) }}\right), \text{ as \ }k\rightarrow \infty
\end{equation*}%
and
\begin{equation*}
\limsup\limits_{k\rightarrow \infty }\left\Vert \sigma_{\alpha_{k}}f-f\right\Vert _{weak-L_{p}(G)}>c_{p}>0,\text{\quad as\quad
}k\rightarrow \infty,  \label{con11f}
\end{equation*}%
where $c_{p}$ is constant depending only on $p$.

However, Simon in \cite{Si3} and \cite{si11} (see also \cite{S, cw, sw}) consider strong convergence theorems of the one-dimensional Walsh-Fouriere series and proved the following:

\textbf{Theorem Si1.}  Let $0<p\leq1$ and $f\in H_{1}\left( G\right)$. Then there exists an absolute constant $c_p$ depending only on $ p, $ such that the following inequality is true:
\begin{equation*}
\frac{1}{\log^{[p]} n}\overset{n}{\underset{k=1}{\sum }}\frac{\left\Vert
	S_{k}f\right\Vert _{H_{p}(G)}}{k^{2-p}}\leq c_p\left\Vert f\right\Vert _{H_{p}(G)},
\end{equation*}

Analigical result for trigonometric system was proved in \cite{sm}, for unbounded Walsh systems in \cite{Ga1}.

In \cite{tep4} was proved that the following is true:

\textbf{Theorem T17.} for any $ 0<p<1 $ and non-decreasing function $\varphi :\mathbb{N}_{+}\rightarrow [1,$ $\infty )$ satisfying the condition
\begin{equation*}
\overline{\lim_{n\rightarrow \infty }}\frac{{n^{2-p}}}{\varphi \left( n\right) }=+\infty,
\end{equation*}
there exists a martingale $f\in H_{p}(G)$, such that
\begin{equation*}
\text{ }\underset{k=1}{\overset{\infty }{\sum }}\frac{\left\Vert S_{k}f\right\Vert _{weak-L_{p}(G)}^{p}}{\varphi (k)}=\infty ,\text{\qquad} (0<p<1).
\end{equation*}

Theorem Si1 follows that if $f\in H_{1}(G)$ then the following equalities are true:
\begin{equation*}
\underset{n\rightarrow \infty }{\lim }\frac{1}{\log n}\overset{n}{\underset{%
		k=1}{\sum }}\frac{\left\Vert S_{k}f-f\right\Vert _{1}}{k}=0
\end{equation*}
and
\begin{equation*}
\underset{n\rightarrow \infty }{\lim }\frac{1}{\log n}\overset{n}{\underset{k=1}{\sum }}\frac{\left\Vert S_{k}f\right\Vert _{1}}{k}=\left\Vert
f\right\Vert _{H_{1}(G)}.
\end{equation*}

When $0<p<1$ and $f\in H_{p}(G),$ then Theorem Si1 follows that there exists an absolute constant $c_p$ depending only on $ p, $ such that 

\begin{equation*}
\frac{1}{n^{1/2-p/2}}\overset{n}{\underset{k=1}{\sum }}\frac{\left\Vert
	S_{k}f\right\Vert _{H_{p}(G)}^{p}}{k^{3/2-p/2}}\leq c_{p}\left\Vert f\right\Vert
_{H_{p}(G)}^{p}.
\end{equation*}
Moreover,
\begin{equation*}
\frac{1}{n^{1/2-p/2}}\overset{n}{\underset{k=1}{\sum }}\frac{\left\Vert S_{k}f-f\right\Vert _{H_{p}(G)}^{p}}{k^{^{3/2-p}}}=0.
\end{equation*}
It follows the following equality
\begin{equation*}
\frac{1}{n^{1/2-p/2}}\overset{n}{\underset{k=1}{\sum }}\frac{\left\Vert
	S_{k}f\right\Vert _{H_{p}(G)}^{p}}{k^{3/2-p/2}}=\left\Vert f\right\Vert
_{H_{p}(G)}^{p}.
\end{equation*}

In the third chapter of this thesis we consider strong convergence results of Fej\'er means of the one-dimensional Walsh-Fourier series. According to Theorem We1 and Theorem Gog2 we only have to consider case $0<p\leq1/2$ (for details see \cite{tep5}, see also \cite{BT1},  \cite{BTT1}, \cite{BT2}, \cite{BT3}, \cite{BPT1}):

\textbf{Theorem \ref{theorem5.3.1}.} Let $0<p\leq 1/2$ and $f\in H_{p}(G)$. Then there exists an absolute constant $c_p$ depending only on $p,$ such that

\begin{equation*}
\frac{1}{\log ^{\left[ 1/2+p\right] }n}\overset{n}{\underset{m=1}{\sum }}%
\frac{\left\Vert \sigma _{m}f\right\Vert _{H_{p}(G)}^{p}}{m^{2-2p}}\leq
c_{p}\left\Vert f\right\Vert _{H_{p}(G)}^{p}.
\end{equation*}%

Moreover, let $0<p<1/2$  and $\Phi :\mathbb{N}_{+}\rightarrow
\lbrack 1, $\textit{\ }$\infty )$ be non-decreasing, non-negative function, such that $\Phi \left( n\right) \uparrow \infty $ and
\begin{equation*}
\overline{\underset{k\rightarrow \infty }{\lim }}\frac{k^{
		2-2p }}{\Phi \left({k}\right) }=\infty .
\end{equation*}
Then there exists a martingale $f\in H_{p}(G),$  such that
\begin{equation*}
\underset{m=1}{\overset{\infty }{\sum }}\frac{\left\Vert \sigma
	_{m}f\right\Vert _{weak-L_{p}(G)}^{p}}{\Phi \left( m\right) }=\infty .
\end{equation*}

When $ p=1/2 $ is was also proved that the following is true:

\textbf{Theorem \ref{theorem5.3.2}.} Let $f\in H_{1/2}(G).$ Then
\begin{equation*}
\underset{n\in \mathbf{\mathbb{N}}_{+}}{\sup }\underset{\left\Vert f\right\Vert _{H_{p}(G)}\leq 1}{\sup } \frac{1}{n}\underset{m=1}{\overset{n}{\sum }}\left\Vert \sigma_{m}f\right\Vert_{1/2}^{1/2}=\infty .
\end{equation*}

Theorem \ref{theorem5.3.1} follows that if $f\in H_{1/2}(G)$ then the following equalities are true:
\begin{equation*}
\lim_{n\rightarrow \infty }\frac{1}{\log n}\overset{n}{\underset{k=1}{\sum }}\frac{\left\Vert \sigma _{k}f-f\right\Vert _{H_{1/2}(G)}^{1/2}}{k}=0
\end{equation*}
and
\begin{equation*}
\lim_{n\rightarrow \infty }\frac{1}{\log n}\overset{n}{\underset{k=1}{\sum }}%
\frac{\left\Vert \sigma _{k}f\right\Vert _{H_{1/2}(G)}^{1/2}}{k}=\left\Vert
f\right\Vert _{H_{1/2}(G)}^{1/2}.
\end{equation*}

When $0<p<1/2$ and $f\in H_{p}(G),$ then Theorem \ref{theorem5.3.1} follows that there exists an absolute constant $c_p$ depending only on $ p, $ such that
\begin{equation*}
\frac{1}{n^{1/2-p}}\overset{n}{\underset{k=1}{\sum }}\frac{\left\Vert \sigma
	_{k}f\right\Vert _{H_{p}(G)}^{p}}{k^{3/2-p}}\leq c_{p}\left\Vert f\right\Vert
_{H_{p}(G)}^{p}.
\end{equation*}
Moreover,
\begin{equation*}
\frac{1}{n^{1/2-p}}\overset{n}{\underset{k=1}{\sum }}\frac{\left\Vert \sigma
	_{k}f-f\right\Vert _{H_{p}(G)}^{p}}{k^{^{3/2-p}}}=0.
\end{equation*}%
It follows that 
\begin{equation*}
\frac{1}{n^{1/2-p}}\overset{n}{\underset{k=1}{\sum }}\frac{\left\Vert \sigma
	_{k}f\right\Vert _{H_{p}(G)}^{p}}{k^{3/2-p}}=\left\Vert f\right\Vert
_{H_{p}(G)}^{p}.
\end{equation*}

\newpage

\section{Partial sums with respect to the one-dimensional Walsh-Fourier series on the martingale Hardy spaces}

\subsection{Basic notations}

\text{ \qquad \qquad \qquad \qquad \qquad \qquad \qquad \qquad \qquad \qquad \qquad \qquad \qquad \qquad \qquad  } Denote by $\mathbb{N}_{+}$ the set of the positive integers and by  $\mathbb{N}:=\mathbb{N}_{+}\cup \{0\}$ the set of non-negative integers. Denote by $Z_{2}$ the additive group of integers modulo-$2,$ which contains only two elements  $Z_{2}:=\{0,1\},$ group operation is modulo-$2$ sum and all sets are open. 

Define the group $G$ as the complete direct product of the groups $Z_{2}$ with the product of the discrete topologies  $Z_{2}$. The direct product $\mu $ of the measures $\mu _{n}\left( \{j\}\right)
:=1/2,\ (j\in Z_{2})$ is the Haar measure on $G$ with $\mu
\left( G\right) =1.$

The elements of $G$ are represented by sequences
\begin{equation*}
x:=(x_{0},x_{1},...,x_{j},...)\qquad \left( x_{k}=0,1\right) .
\end{equation*}

It is easy to give a base for the neighbourhood of $G$
\begin{equation*}
I_{0}\left( x\right) :=G,
\end{equation*}%
\begin{equation*}
I_{n}(x):=\{y\in G\mid y_{0}=x_{0},...,y_{n-1}=x_{n-1}\}\text{ }(x\in G,%
\text{ }n\in \mathbb{N}).
\end{equation*}

Set  $I_{n}:=I_{n}\left( 0\right) $ for any $n\in \mathbb{N}$ and $\overline{I_{n}}:=G$ $\backslash $ $I_{n}$.

It is evident that
\begin{equation} \label{2.1.2}
\overline{I_{M}}=\left( \overset{M-2}{\underset{k=0}{\bigcup }}\overset{M-1}{%
	\underset{l=k+1}{\bigcup }}I_{l+1}\left( e_{k}+e_{l}\right) \right) \bigcup \left( \underset{k=0}{\bigcup\limits^{M-1}}I_{M}\left( e_{k}\right) \right)= \underset{k=0}{\bigcup\limits^{M-1}}I_{k} \backslash I_{k+1}.
\end{equation}

If $n\in \mathbb{N},$ then it can be uniquely expressed as
$
n=\sum_{k=0}^{\infty }n_{j}2^{j}
$
where  $n_{j}\in Z_{2}$  $~(j\in \mathbb{N})$ and only a finite number of $n_{j}$s differ from zero. Set
\begin{equation*}
\left\langle n\right\rangle :=\min \{j\in \mathbb{N},n_{j}\neq 0\}\text{ \ \
	and \ \ \ }\left\vert n\right\vert :=\max \{j\in \mathbb{N},n_{j}\neq 0\},
\end{equation*}%
It is evident that $2^{\left\vert n\right\vert }\leq n\leq 2^{\left\vert n\right\vert
	+1}.$ Let
\begin{equation*}
d\left( n\right) :=\left\vert n\right\vert -\left\langle n\right\rangle,
\text{ \ for any \ }n\in \mathbb{N}.
\end{equation*}

Denote by $V(n)$ variation of natural number $n\in \mathbb{N}$ 

\begin{equation*}
V\left( n\right) =n_{0}+\overset{\infty }{\underset{k=1}{\sum }}\left|
n_{k}-n_{k-1}\right| .
\end{equation*}

Define $k$-th Rademacher functions by
\begin{equation*}
r_{k}\left( x\right) :=\left( -1\right) ^{x_{k}}\text{\qquad }\left( \text{ }%
x\in G,\text{ }k\in \mathbb{N}\right) .
\end{equation*}

By using Rademacher functions we define Walsh system $w:=(w_{n}:n\in \mathbb{N})$ $G$ as:
\begin{equation*}
w_{n}(x):=\underset{k=0}{\overset{\infty }{\prod }}r_{k}^{n_{k}}\left(x\right) =r_{\left\vert n\right\vert }\left( x\right) \left( -1\right) ^{%
	\underset{k=0}{\overset{\left\vert n\right\vert -1}{\sum }}n_{k}x_{k}}\text{%
	\qquad }\left( n\in \mathbb{N}\right).
\end{equation*}

The norm (quasi-norm) of space $L_{p}(G)$ and $weak-L_{p}(G) $ for $\ \left( 0<p<\infty \right)$ are respectively defined as
\begin{equation*}
\left\Vert f\right\Vert _{p}^p:=\left(\int_{G}\left\vert f(x)\right\vert
^{p}d\mu (x)\right),  \ \ \ \ \ \left\Vert f\right\Vert _{weak-L_{p}(G)}^p:=\underset{\lambda >0}{\sup}
\lambda^p \mu \left (x\in{G}:  |f|>\lambda \right).
\end{equation*}

Walsh system is orthonormal and complete in $L_{2}\left(G \right) $ (see  \cite{sws}).

For any $f\in L_{1}\left(G\right) $ the numbers
\begin{eqnarray*}
	\widehat{f}\left(n\right) :=\int\limits_{G}f(x)w_{n}(x)d\mu({x})
\end{eqnarray*}
are called $n$-th Walsh-Fourier coefficient of $f$.

$n$-th partial sum is denoted by
\begin{equation*}
S_{n}(f;{x}):=\sum\limits_{i=0}^{n-1}
\widehat{f}\left(i\right) w_{i}\left( x\right) .
\end{equation*}

Dirichlet kernels are defined by
\begin{equation*}
D_{n}\left( x\right) :=\sum\limits_{i=0}^{n-1} w_{i}\left( x\right) .
\end{equation*}

We also define the following maximal operators
\begin{eqnarray*}
	{S}^{\ast}f=\sup_{n\in \mathbb{N}}\left\vert S_{{n}}f\right\vert,
	\ \ \ \ \ \ \ \ \ \ 
	\widetilde{S}_{\#}^{\ast}f=\sup_{n\in \mathbb{N}}\left\vert S_{2^{n}}f\right\vert.
\end{eqnarray*}

The $\sigma $-algebra generated by the intervals $ I_{n}(x) $
with measure $2^{-n}$ is denoted by $\digamma_{n}\left(n\in \mathbb{N}\right).$ Conditional exponential operator with respect to
$\digamma _{n}\left( n\in \mathbb{N}\right) $ is denoted by $E_{n}$ and it is given by 

\begin{eqnarray*}
	E_{n}f(x)= S_{2^{n}}f\left( x\right)=\sum_{k=0}^{2^{n}-1}
	\widehat{f}\left(k\right) w_{k}(x) = \frac{1}{\left| I_{n}\left(x\right) \right| }\int_{I_{n}\left(x\right)
	}f(x)d\mu (x),
\end{eqnarray*}
where $\left| I_{n}\left(x\right) \right| =2^{-n}$  denotes length of set $I_{n}\left(x\right)$.

Sequence $f=\left(f_{n},\text{ }n\in \mathbb{N}\right) $ of functions $f_{n}\in L_{1}\left( G\right) $ is called dyadic martingale  (for details see \cite{Nev}, \cite{sws}) if

$\left( i\right) $ $f_{n}$ is measurable with respect to  $\sigma-$ algebras $\digamma _{n}$ for any $n\in
\mathbb{N}$,

$\left( ii\right) $ $E_{n}f_{m}=f_{n}$ for any $n\leq m$.

The maximal function of a martingale $f$ is defined by
\begin{equation*}
f^{\ast }=\sup_{n\in \mathbb{N}}\left\vert f_{n}\right\vert.
\end{equation*}

In the case $f\in L_{1}\left(G\right) ,$ the maximal functions are also be given by:
\begin{equation*}
f^{\ast }\left(x\right) =\sup\limits_{n\in \mathbb{N}}\frac{1}{\mu \left(
	I_{n}\left( x\right) \right) }\left\vert \int\limits_{I_{n}\left( x\right)
}f\left( u\right) d\mu \left( u\right) \right\vert.
\end{equation*}

For $0<p<\infty $ the Hardy martingale space $H_{p}\left( G\right)$
consists of all martingales, for which
\begin{equation*}
\left\Vert f\right\Vert _{H_{p}\left( G\right)}:=\left\Vert f^{\ast }\right\Vert
_{p}<\infty .
\end{equation*}

A bounded measurable function $a$ is said to be a $p$-atom if there exists an
dyadic interval $I$, such that%
\begin{equation*}
\int_{I}ad\mu =0,\text{ \ \ }\left\Vert a\right\Vert _{\infty }\leq \mu
\left( I\right) ^{-1/p},\text{ \ \ supp}\left( a\right) \subset I.\qquad
\end{equation*}%

It is easy to show that for martingale $f=(f_{n},n\in \mathbb{N}) $ and for any ${k}\in {\mathbb{N}}$ there exists a limit
\begin{equation*}
\widehat{f}\left(k\right):=\lim_{{n}\rightarrow \infty }\int_{{G}}f_{{{n}}}\left(
{x}\right) w_{{k}}\left( {x}\right) d\mu \left({x}\right)
\end{equation*}
and it is called  $k$-th Walsh-Fourier coefficients of $f$.

If $f_{0}\in L_{1}\left( G\right)$ and $f:=(E_{n}f_{0}:n\in \mathbb{N})$ is regular martingale then
\begin{eqnarray*}
	\widehat{f}\left(k\right)&=&\int_{{G}}f\left( {x}\right) w_{k}(x)
	d\mu \left({x}\right)=\widehat{f_{0}}\left(k\right),\ \ k\in \mathbb{N}.
\end{eqnarray*}

The modulus of continuity in  $H_{p}(G)$ space is defined by
\begin{equation*}
\omega _{H_{p}\left( G\right)}\left( \frac{1}{2^{n}},f\right) :=\left\Vert
f-S_{2^{n}}f\right\Vert _{H_{p}\left( G\right)}.
\end{equation*}

It is important to describe how can be understood difference $f-S_{2^{n}}f$, where $f$ be martingale $S_{2^{n}}f$ is a function:

\begin{remark} \label{lemma2.3.6666}
	Let $0<p\leq 1.$ Since
	\begin{equation*}
	S_{2^{n}}f=f^{\left( n\right) }\in L_1(G),\text{ \ where \ }f=\left( f^{\left( n\right)}:n\in \mathbb{N}\right) \in H_{p}(G)
	\end{equation*}%
	and
	\begin{eqnarray*}
		\left( S_{2^{k}}f^{\left( n\right) }:k \in \mathbb{N}\right)
		&=&\left( S_{2^{k}}S_{2^{n}},k \in \mathbb{N}\right) \\
		&=&\left( S_{2^{0}}f,\ldots ,S_{2^{n-1}}f,S_{2^{n}}f,S_{2^{n}}f,\ldots \right)=\left( f^{\left( 0\right) },\ldots ,f^{\left( n-1\right) },f^{\left(
			n\right) },f^{\left( n\right) },\ldots \right).
	\end{eqnarray*}
	Under the difference $  f-S_{2^{n}}f $ we mean the following martingale:
	$$ f:=(\left( f-S_{2^{n}}f\right) ^{\left( k\right) },  \ \ k\in \mathbb{N} )$$
	where
	\begin{equation*}
	\left( f-S_{2^{n}}f\right) ^{\left( k\right) }=\left\{
	\begin{array}{ll}
	0, & k=0,\ldots ,n, \\
	f^{\left( k\right) }-f^{\left( n\right) }, & k\geq n+1,\end{array}
	\right.
	\end{equation*}
\end{remark}

Consequently, the norm $ \left\Vert f-S_{2^{n}}f\right\Vert _{H_{p}\left( G\right)}  $ is understood as  $ H_p $-norm of 
$$  f-S_{2^{n}}f=(\left( f-S_{2^{n}}f\right) ^{\left( k\right) }, k\in \mathbb{N}) $$.

Watari \cite{wat} showed that there are strong connections between
\begin{equation*}
\omega _{p}\left( \frac{1}{2^{n}},f\right) ,\text{ \ }E_{2^{n}}\left(
L_{p},f\right) \text{ \ \ and\ \ \ }\left\Vert f-S_{2^{n}}f\right\Vert
_{p},\ \ p\geq 1,\text{ \ }n\in \mathbb{N}.
\end{equation*}

In particular,%
\begin{equation*}
\frac{1}{2}\omega _{p}\left( \frac{1}{2^{n}},f\right) \leq \left\Vert
f-S_{2^{n}}f\right\Vert _{p}\leq \omega _{p}\left( \frac{1}{2^{n}},f\right)
\label{eqvi}
\end{equation*}%
and%
\begin{equation*}
\frac{1}{2}\left\Vert f-S_{2^{n}}f\right\Vert _{p}\leq E_{2^{n}}\left(
L_{p},f\right) \leq \left\Vert f-S_{2^{n}}f\right\Vert _{p}.
\end{equation*}

\subsection{Auxiliary lemmas}

\text{ \qquad \qquad \qquad \qquad \qquad \qquad \qquad \qquad \qquad \qquad \qquad \qquad \qquad \qquad \qquad  } First we present and prove equalities and estimations of Dirichlet kernel and Lebesgue constants with respect to the one-dimensional Walsh-Fourier systems (see Lemmas \ref{lemma0}-\ref{lemma7}).

First equality of the following Lemma is proved in \cite{sws} and second identity is proved in  \cite{Ga2}:

\begin{lemma} \label{lemma0}
	Let $j, n\in \mathbb{N}$. Then
	\begin{equation*}
	D_{j+2^n}=D_{2^n}+w_{2^n}D_j,\text{ \ when \ } j\leq 2^n,
	\end{equation*}
	and
	\begin{equation*}
	D_{2^n-j}=D_{2^n}-\psi_{2^n-1}{D_j},\text{ \ when \ }j<2^n.
	\end{equation*}
\end{lemma}

The following estimation of Dirichlet kernel with respect to the one-dimensional Walsh-Fourier systems is proved in \cite{sws}:

\begin{lemma} \label{lemma1}
	
	Let $n\in \mathbb{N}$. Then
	
	\begin{equation*}
	D_{2^{n}}\left( x\right) =\left\{
	\begin{array}{l}
	2^{n},\text{ \ if  \  }x\in I_{n},
	\\ 0,\text{ \ if \ }x\notin I_{n},
	\end{array}
	\right.
	\end{equation*}
	and
	\begin{equation*}
	D_{n}=w_{n}\overset{\infty }{\underset{k=0}{\sum }}n_{k}r_{k}D_{2^{k}}=w_{n}%
	\overset{\infty }{\underset{k=0}{\sum }}n_{k}\left(
	D_{2^{k+1}}-D_{2^{k}}\right) ,\text{ \ for \ }n=\overset{\infty }{\underset{i=0
		}{\sum }}n_{i}2^{i}.
	\end{equation*}%
\end{lemma}

The following two-sided estimations of Lebesgue constants with respect to the one-dimensional Walsh-Fourier systems is proved in \cite{sws} and second equality is proved in \cite{fi}:

\begin{lemma} \label{lemma2}
	
	Let $n\in \mathbb{N}.$ Then
	
	\begin{equation*}
	\frac{1}{8}V\left( n\right) \leq \left\Vert D_{n}\right\Vert _{1}\leq
	V\left( n\right)
	\end{equation*}
	and
	\begin{equation*}
	\frac{1}{n\log n}\underset{k=1}{\overset{n}{\sum }}V\left( k\right) =\frac{1%
	}{4\log 2}+o\left( 1\right).
	\end{equation*}
	
\end{lemma}

Hardy martingale space $H_{p}\left( G\right) $ for any $0<p\leq 1$ can be characterize by simple functions which are called $ p $-atoms. The following is true (for details see \cite{S}, \cite{We1} and \cite{We5}):

\begin{lemma} \label{lemma3.2.4}
	A martingale $f=\left(f_{n},\text{ }n\in \mathbb{N}%
	\right) $ belongs to $H_{p}(G)\left( 0<p\leq 1\right) $ if and only if there exists a sequence of $ p $-atoms of
	$\left( a_{k},k\in \mathbb{N}\right) $ and sequence of real numbers 
	$\left( \mu _{k},\text{ }k\in \mathbb{N}\right) $ such that for all $n\in \mathbb{N}$,
	
	\begin{equation} \label{3.2.1.lemma3.2.4}
	\qquad \sum_{k=0}^{\infty }\mu _{k}S_{2^{n}}a_{k}=f_{n}
	\end{equation}%
	and
	\begin{equation*}
	\qquad \sum_{k=0}^{\infty }\left\vert \mu _{k}\right\vert ^{p}<\infty.
	\end{equation*}%
	Moreover,
	\begin{equation*}
	\left\Vert f\right\Vert _{H_{p}(G)}\backsim \inf \left(
	\sum_{k=0}^{\infty }\left\vert \mu _{k}\right\vert ^{p}\right) ^{1/p},
	\end{equation*}
	where the infimum is taken over all decomposition of  $f$ of the form (\ref{3.2.1.lemma3.2.4}).
\end{lemma}

The next five Examples of martingales will be used many times to prove sharpness of our main results. Such counterexamples first appear in the papers  of Goginava  \cite{Go} (see also \cite{gog4,GoPubl}). Such constructions of martingales are also used in the papers \cite{BPTW}, \cite{BT1}, \cite{BNPT}, \cite{BNT100}, \cite{MPT}, \cite{LPTT}, \cite{MST}, \cite{tep16}, \cite{tep18}, \cite{tep17}, \cite{NT4}, \cite{PTT100}, \cite{PTW1}, \cite{PTW2}, \cite{ptw}, \cite{tep1}, \cite{tep5}, \cite{tep12}, \cite{tep13}, \cite{tep14}, \cite{tep15},  \cite{tep20}, \cite{tep19}, \cite{tep22}, \cite{tep30}, \cite{tep23}, \cite{tep21}, \cite{tep31}, \cite{tetu1}.
For the one-dimensional case we use martingales which were used in \cite{tep_thesis}. So, we leave out the details of proof.

\begin{example}  \label{example2.2.1}
	Let $0<p\leq 1,$  $\left\{ \lambda _{k}:k\in \mathbb{N}\right\}$
	be sequence of real numbers
	\begin{equation} \label{3.3.2aa}
	\sum_{k=0}^{\infty }\left\vert \lambda _{k}\right\vert ^{p}\leq c_{p}<\infty
	\end{equation}%
	and $\left\{ a_{k}:k\in \mathbb{N}\right\} $ be sequence of $p$-atoms, given by 
	\begin{equation*}
	a_{k}(x):={2^{\left\vert \alpha _{k}\right\vert(1/p-1)}}\left(
	D_{2^{\left\vert \alpha _{k}\right\vert +1}}(x)-D_{2^{{\left\vert \alpha_{k}\right\vert }}}(x)\right),
	\end{equation*}%
	where $\left\vert \alpha _{k}\right\vert :=\max $ $\{j\in \mathbb{N}:$ $%
	\left(\alpha_{k}\right)_{j}\neq 0\}$ and $\left(\alpha_{k}\right)_{j}$
	denotes $j$-th binary coefficients of real number of $\alpha _{k}\in\mathbb{N}_{+}$. Then
	$\,f=\left( f_n:\text{ }n\in \mathbb{N}\right),$ where
	\begin{equation*}
	f_n(x):=\sum_{\left\{ k:\text{ }\left\vert \alpha
		_{k}\right\vert <n\right\} }\lambda _{k}a_{k}(x)
	\end{equation*}
	is martingale, which belongs to $H_{p}(G)$ for any $0<p\leq 1$.
	
	It is easy to show that
	\begin{equation}  \label{3.3.10AA}
	\widehat{f}(j)
	\end{equation}%
	\begin{equation*}
	=\left\{
	\begin{array}{ll}
	{\lambda _{k}2^{(1/p-1)\left\vert \alpha _{k}\right\vert }},
	& j\in \left\{ 2^{\left\vert \alpha _{k}\right\vert },...,
	2^{\left\vert \alpha _{k}\right\vert +1}-1\right\} ,\text{ }k\in \mathbb{N}%
	_{+}, \\
	0, & \text{\thinspace }j\notin \bigcup\limits_{k=1}^{\infty }\left\{
	2^{\left\vert \alpha _{k}\right\vert },...,2^{\left\vert \alpha_{k}\right\vert +1}-1\right\}.
	\end{array}%
	\right.
	\end{equation*}
	Let
	$2^{\left\vert \alpha _{l-1}\right\vert +1}\leq j\leq 2^{\left\vert\alpha_{l}\right\vert },$\ $l\in \mathbb{N}_{+}.$
	Then
	\begin{eqnarray}  \label{3.3.12AA}
	S_{j}f= S_{2^{\left\vert \alpha _{l-1}\right\vert +1}} 
	= \sum_{\eta =0}^{l-1}{\lambda _{\eta }2^{_{\left\vert \alpha _{\eta }\right\vert }(1/p-1)}}\left( D_{2^{\left\vert \alpha _{\eta }\right\vert
			+1}}-D_{2^{\left\vert \alpha _{\eta }\right\vert }}\right) .
	\end{eqnarray}
	
	Let $2^{\left\vert \alpha _{l}\right\vert }\leq j<2^{\left\vert \alpha_{l}\right\vert +1},$ $l\in \mathbb{N}_{+}.$ Then
	\begin{eqnarray}  \label{3.3.11AA}
	S_{j}f
	&=& S_{2^{\left\vert \alpha _{l}\right\vert }}+{\lambda
		_{l}2^{(1/p-1)\left\vert \alpha _{l}\right\vert}w_{2^{\left\vert
				\alpha _{l}\right\vert }}D_{j-2^{_{\left\vert \alpha _{l}\right\vert}}}} \\  \notag
	&=&\sum_{\eta =0}^{l-1}{\lambda _{\eta }2^{(1/p-1)\left\vert \alpha _{\eta}\right\vert}}\left( D_{2^{{\left\vert \alpha _{\eta
				}\right\vert +1}}}-D_{2^{{\left\vert \alpha _{\eta }\right\vert }}}\right) \\ \notag
	&+&{\lambda _{l}2^{(1/p-1)\left\vert \alpha _{l}\right\vert }w_{2^{\left\vert \alpha _{l}\right\vert }}D_{j-2^{\left\vert \alpha_{l}\right\vert }}}.
	\end{eqnarray}
	
	Moreover, \ for \ the \ modulus \ of \ continuity \ for \ $0<p\leq 1$ \ we \  have \ the \ following \ estimation:
	\begin{equation} \label{3.3.2aa0}
	\omega _{H_{p}}\left( \frac{1}{2^{n}},f\right) =O\left( \sum_{\left\{ k:%
		\text{ }\left\vert \alpha _{k}\right\vert \geq n\right\} }^{\infty
	}\left\vert \lambda _{k}\right\vert ^{p}\right) ^{1/p},\text{ \ \ as \ \ }n\rightarrow \infty.
	\end{equation}%
\end{example}

By applying Lemma \ref{lemma3.2.4} we easily obtain that the following is true (see \cite{We5}):

\begin{lemma} \label{lemma3.2.5}
	\label{lemma2.3} Let $0<p\leq 1$ and $T$ be $ \sigma  $-sub-linear operator, such that, for any $p$-atom  $a$,
	\begin{equation*}
	\int\limits_{G}\left\vert Ta({x})\right\vert ^{p}d\mu({x}) \leq
	c_{p}<\infty.
	\end{equation*}
	Then
	\begin{equation} \label{3.2.5000}
	\left\Vert Tf\right\Vert _{p}\leq c_{p}\left\Vert f\right\Vert
	_{H_{p}(G)}.
	\end{equation}
	
	In addition, if $T$  is bounded from  $L_{\infty}(G)$ to $L_{\infty}(G)$ then to prove (\ref{3.2.5000}) it is suffices to show that  
	\begin{equation*}
	\int\limits_{\overset{-}{I}}\left\vert Ta({x})\right\vert ^{p}d\mu({x}) \leq
	c_{p}<\infty,
	\end{equation*}
	for every $p$-atom  $a$, where $I$ denotes support of the atom $a$.
\end{lemma}

In the concrete cases the norm of Hardy martingale spaces can be calculated  by simpler formulas (for details see \cite{S}, \cite{We1} and \cite{We3}):

\begin{lemma}  \label{lemma3.2.3}
	If $g\in L_{1}\left( G\right)$ and $f:=(E_{n}g:n\in \mathbb{N})$ be regular martingale, then  $H_{p}\left( G\right)\text{ \ for \ }0<p\leq1$ norm can be calculated by
	\begin{equation*}
	\left\Vert f\right\Vert _{H_{p}(G)}=\left\Vert \sup\limits_{n\in \mathbb{N}}|S_{2^n}g|\right\Vert _{p}.
	\end{equation*}
\end{lemma}

The following lemmas are proved in \cite{tep5}, \cite{tep12}, \cite{tep13}.
\begin{lemma} \label{lemma3.2.8}
	\label{example002} Let $0<p\leq 1$, $2^{k}\leq n<2^{k+1}$ and $S_{n}f$ be $n$-th partial sum with respect to the one-dimensional Walsh-Fourier series, where $f\in H_{p}(G)$. Then for any fixed $n\in \mathbb{N}$,
	\begin{eqnarray*}
		\left\Vert S_{n}f\right\Vert _{H_{p}(G)}^p\leq\left\Vert \sup_{0\leq l\leq
			k}\left\vert S_{2^{l}}f\right\vert \right\Vert _{p}^p+\left\Vert
		S_{n}f\right\Vert _{p}^p 
		\leq \left\Vert \widetilde{S}_{\#}^{\ast }f\right\Vert _{p}^p+\left\Vert
		S_{n}f\right\Vert _{p}^p.
	\end{eqnarray*}
\end{lemma}
{\bf Proof}:
Let consider the following martingales
\begin{equation*}
f_{\#}:=\left( S_{2^{k}}S_{n}f,\text{ }k\in\mathbb{N}_+\right)
=\left(
S_{2^{0}},S_{2^{k}}f,\text{ }S_{n}f,...,S_{n}f,...\right),
\end{equation*}%
Hence, Lemma \ref{lemma3.2.3} immediately follows that
\begin{eqnarray*}
	\left\Vert S_{n}f\right\Vert _{H_{p}(G)}^p 
	\leq \left\Vert \sup_{0\leq l\leq k}\left\vert S_{2^{l}}f\right\vert \right\Vert _{p}^p+\left\Vert
	S_{n}f\right\Vert _{p}^p \leq \left\Vert \widetilde{S}_{\#}^{\ast }f\right\Vert
	_{p}^p+\left\Vert S_{n}f\right\Vert _{p}^p.
\end{eqnarray*}

Lemma is proved.

\subsection{Boundedness of subsequences of partial sums with respect to the one-dimensional Walsh-Fourier series on the martingale Hardy spaces} \text{ \qquad \qquad \qquad \qquad \qquad \qquad \qquad \qquad \qquad \qquad \qquad \qquad \qquad \qquad \qquad  }
In this section we consider boundedness of subsequences of partial sums with respect to the one-dimensional Walsh-Fourier series in the martingale Hardy spaces (for details see \cite{tep12}).

\begin{theorem} \label{th4.1.1}
	a) Let $0<p<1$ and $f\in H_{p}(G)$.  Then there exists an absolute constant $c_{p}$ depending only on $p,$ such that
	\begin{equation*}
	\text{ }\left\Vert S_{n}f\right\Vert _{H_{p}(G)}\leq c_{p}2^{d\left( n\right)
		\left( 1/p-1\right) }\left\Vert f\right\Vert _{H_{p}(G)}.
	\end{equation*}
	
	b) Let $0<p<1,$ $\left\{ m_{k}:\text{ }k \in\mathbb{N}_{+}\right\}$ be non-negative, increasing  sequence of natural numbers such that
	\begin{equation}
	\sup_{k\in \mathbb{N}}d\left( m_{k}\right) =\infty  \label{4.1.dnk}
	\end{equation}%
	and $\Phi :\mathbb{N}_{+}\rightarrow \lbrack 1,\infty )$ be non-decreasing function satisfying the condition
	\begin{equation}
	\overline{\underset{k\rightarrow \infty }{\lim }}\frac{2^{d\left(
			m_{k}\right) \left( 1/p-1\right) }}{\Phi \left( m_{k}\right) }=\infty.
	\label{4.1.1010}
	\end{equation}
	Then there exists a martingale $f\in H_{p}(G)$ such that
	\begin{equation*}
	\underset{k\in \mathbb{N}}{\sup }\left\Vert \frac{S_{m_{k}}f}{\Phi \left(
		m_{k}\right) }\right\Vert _{weak-L_{p}(G)}=\infty .
	\end{equation*}
\end{theorem}

{\bf Proof}:
Suppose that
\begin{equation}
\left\Vert 2^{\left( 1-1/p\right) d\left( n\right) }S_{n}f\right\Vert
_{p}\leq c_{p}\left\Vert f\right\Vert _{H_{p}(G)}.  \label{4.1.11.1}
\end{equation}

By combining Lemma \ref{lemma3.2.8} and inequalities (\ref{1.S2n}) and (\ref{4.1.11.1}), since $  2^{\left( 1-1/p\right) d\left( n\right) }\leq c_{p} $ we obtain that
\begin{eqnarray} \label{4.1.11.2}
\left\Vert 2^{\left( 1-1/p\right) d\left( n\right) }S_{n}f\right\Vert
_{H_{p}(G)}^p &\leq& \left\Vert  2^{\left( 1-1/p\right) d\left( n\right) }S_nf\right\Vert _{p}^p+2^{\left( 1-1/p\right) d\left( n\right) }\left\Vert\widetilde{S}_{\#}^{\ast}f\right\Vert_{p}^p \\ \notag
&\leq & c_{p}\left\Vert f\right\Vert _{H_{p}(G)}^p+c_{p}\left\Vert\widetilde{S}_{\#}^{\ast}f\right\Vert
_{p}^p \leq  c_{p}\left\Vert f\right\Vert _{H_{p}(G)}^p.
\end{eqnarray}
By combining Lemma \ref{lemma3.2.5} and (\ref{4.1.11.2}) it is suficies to show that 
\begin{equation}
\int\limits_{G}\left\vert 2^{\left( 1-1/p\right) d\left( n\right)
}S_{n}a\right\vert ^{p}d\mu \leq c_{p}<\infty ,  \label{4.1.25a}
\end{equation}%
for every $ p $-atom $a$, with support $I$, such that
$\mu \left(I\right) =2^{-M}$.

Without loss the generality we may assume that $p$-atom $a$ has support $I=I_{M}.$ Then it is easy to see that  $S_{n}a =0,$ where $2^{M}$ $\geq n$.
So, we may assume that $2^{M}<n$. Since $\left\Vert a\right\Vert_{\infty}\leq 2^{M/p}$ we can conclude that
\begin{eqnarray} \label{4.1.11a}
\left\vert 2^{\left( 1-1/p\right) d\left( n\right) }S_{n}a\left( x\right)
\right\vert &\leq& 2^{\left( 1-1/p\right) d\left( n\right) }\left\Vert
a\right\Vert _{\infty }\int_{I_{M}}\left\vert D_{n}\left( x+t\right)
\right\vert d\mu \left( t\right)  \\ \notag
&\leq & 2^{M/p}2^{\left( 1-1/p\right) d\left( n\right) }\int_{I_{M}}\left\vert
D_{n}\left( x+t\right) \right\vert d\mu \left(t\right).
\end{eqnarray}

Let $x\in I_{M}$. Since $V\left( n\right) \leq 2d\left( n\right) ,$  by using first estimations of Lemma \ref{lemma2} we can conclude that
\begin{eqnarray*}
	\left\vert 2^{\left( 1-1/p\right) d\left( n\right)}S_{n}a\right\vert 
	\leq2^{M/p}2^{\left( 1-1/p\right) d\left( n\right) }V\left( n\right) 
	\leq 2^{M/p}d\left( n\right) 2^{\left( 1-1/p\right) d\left( n\right) }
\end{eqnarray*}%
and
\begin{eqnarray} \label{4.1.11b}
\int_{I_{M}}\left\vert 2^{\left( 1-1/p\right) d\left( n\right)
}S_{n}a\right\vert ^{p}d\mu \leq d\left( n\right) 2^{\left( 1-1/p\right)
	d\left( n\right) }<c_{p}<\infty .
\end{eqnarray}

Let $t\in I_{M}$ and $x\in I_{s}\backslash I_{s+1},\,$ where $0\leq s\leq
M-1<\left\langle n\right\rangle $ or $0\leq s<\left\langle n\right\rangle
\leq M-1.$ Then $x+t$ $\in I_{s}\backslash I_{s+1}$ and if we use both equality of Lemma \ref{lemma1} we get that $D_{n}\left( x+t\right) =0$ and it follows that
\begin{equation}
\left\vert 2^{\left( 1-1/p\right) d\left( n\right) }S_{n}a\left( x\right)
\right\vert =0. \label{4.1.11bbb}
\end{equation}

Let $x\in I_{s}\backslash I_{s+1},$ $\,\left\langle n\right\rangle \leq
s\leq M-1.$ Then $x+t\in I_{s}\backslash I_{s+1},$ where $t\in I_{M}$. Then by using again both equality of Lemma \ref{lemma1} we have that
\begin{equation*}
\left\vert D_{n}\left( x+t\right) \right\vert \leq
\sum_{j=0}^{s}n_{j}2^{j}\leq c2^{s}.
\end{equation*}%
If we apply again (\ref{4.1.11a}) we can conclude that
\begin{eqnarray} \label{4.1.12}
\left\vert 2^{\left( 1-1/p\right) d\left( n\right) }S_{n}a\left( x\right)
\right\vert
&\leq & 2^{\left( 1-1/p\right) d\left( n\right) }2^{M/p}\frac{2^{s}%
}{2^{M}} \\ \notag
&\leq & 2^{\left< n\right> \left(1/p-1\right) }2^{M(1/p-1)}\frac{2^{s}%
}{2^{\left| n\right| \left(1/p-1\right) }} \leq 2^{\left\langle n\right\rangle \left( 1/p-1\right) }2^{s}.
\end{eqnarray}

By identity (\ref{2.1.2})  and inequalities  (\ref{4.1.11bbb}) and (\ref{4.1.12}) we find that
\begin{eqnarray*}
	&&\int_{\overline{I_{M}}}\left\vert 2^{\left( 1-1/p\right) d\left( n\right)
	}S_{n}a\left( x\right) \right\vert ^{p}d\mu \left( x\right) \\
	&=&\overset{M-1}{\underset{s=\left\langle n\right\rangle }{\sum }}%
	\int_{I_{s}\backslash I_{s+1}}\left\vert 2^{\left\langle n\right\rangle
		\left( 1/p-1\right) }2^{s}\right\vert ^{p}d\mu \left( x\right)\leq  c\overset%
	{M-1}{\underset{s=\left\langle n\right\rangle }{\sum }}\frac{2^{\left\langle
			n\right\rangle \left( 1-p\right) }}{2^{s\left( 1-p\right) }}\leq
	c_{p}<\infty .
\end{eqnarray*}

Now, we prove part b) of Theorem \ref{th4.1.1}. By using condition (\ref{4.1.1010}) there exists sequence of natural numbers $\left\{ \alpha _{k}:\text{ }k\in\mathbb{N}_+\right\} \subset \left\{
m_{k}:\text{ }k\in\mathbb{N}_+\right\},$ such that
\begin{equation}
\sum_{\eta =0}^{\infty }\frac{\Phi ^{p/2}\left( \alpha _{\eta }\right) }{%
	2^{d\left( \alpha _{\eta }\right) \left( 1-p\right) /2}}<\infty ,
\label{4.1.12f}
\end{equation}

Let $f=\left(f_{n},n\in\mathbb{N}_+\right) \in H_{p}(G)$ be a martingale from the Example \ref{example2.2.1}, where
\begin{equation} \label{4.1.100}
\lambda _{k}=\frac{\Phi ^{1/2}\left( \alpha _{k}\right) }{2^{d\left( \alpha
		_{k}\right) \left( 1/p-1\right) /2}}.
\end{equation}

Then, if we use (\ref{4.1.12f}) we obtain that condition (\ref{3.3.2aa})  is fulfilled and it follows that $f=\left(f_{n},n\in\mathbb{N}_+\right) \in H_{p}(G).$

If we apply (\ref{3.3.10AA}) when $\lambda _{k}$ are given by the formula (\ref{4.1.100}) then we get that
\begin{equation} \label{4.1.6}
\widehat{f}(j)
\end{equation}
\begin{equation*}
=\left\{
\begin{array}{c}
\Phi ^{1/2}\left( \alpha _{k}\right) 2^{\left( \left\vert \alpha
	_{k}\right\vert +\left\langle \alpha _{k}\right\rangle \right) \left(
	1/p-1\right) /2},\text{ \ if \thinspace \thinspace }j\in \left\{
2^{_{\left\vert \alpha _{k}\right\vert }},...,2^{_{\left\vert \alpha
		_{k}\right\vert +1}}-1\right\} ,\text{ }k\in\mathbb{N}_+ \\
0,\text{ \ if \ }j\notin \bigcup\limits_{k=0}^{\infty }\left\{ 2^{_{\left\vert \alpha _{k}\right\vert}},...,2^{_{\left\vert \alpha _{k}\right\vert +1}}-1\right\} .\text{ }%
\end{array}\right.
\end{equation*}

In the view of (\ref{3.3.11AA}) when $\lambda _{k}$ are given by (\ref{4.1.100}) we get that
\begin{eqnarray} \label{4.1.6aaa}
\frac{S_{\alpha _{k}}f}{\Phi \left( \alpha _{k}\right) } 
&=&\frac{1}{\Phi \left( \alpha _{k}\right) }\sum_{\eta =0}^{k-1}\Phi
^{1/2}\left(\alpha _{\eta }\right) 2^{\left( \left\vert \alpha _{\eta
	}\right\vert +\left\langle \alpha _{\eta }\right\rangle \right) \left(
	1/p-1\right) /2}\left( D_{2^{\left\vert \alpha _{\eta }\right\vert
		+1}}-D_{2^{\left\vert \alpha _{\eta }\right\vert}}\right) \\ \notag
&+&\frac{2^{\left( \left\vert \alpha _{k}\right\vert +\left\langle \alpha
		_{k}\right\rangle \right) \left( 1/p-1\right) /2}w_{2^{\left\vert \alpha
			_{k}\right\vert }}D_{\alpha _{k}-2^{\left\vert \alpha _{k}\right\vert }}}{%
	\Phi ^{1/2}\left( \alpha _{k}\right) }:=I+II.
\end{eqnarray}

by using (\ref{4.1.12f}) for $I$ we have that
\begin{eqnarray} \label{4.1.8a}
&&\left\Vert I\right\Vert _{weak-L_{p}(G)}^{p}  \\ \notag
&\leq &\frac{1}{\Phi ^{p}\left( \alpha _{k}\right) }\sum_{\eta =0}^{k-1}\frac{%
	\Phi ^{p/2}\left( \alpha _{\eta }\right) }{2^{d\left( \alpha _{\eta }\right)
		\left( 1-p\right) /2}}\left\Vert 2^{\left\vert \alpha _{\eta }\right\vert
	\left( 1/p-1\right) }\left( D_{2^{\left\vert \alpha _{\eta }\right\vert
		+1}}-D_{2^{\left\vert \alpha _{\eta }\right\vert}}\right) \right\Vert
_{weak-L_{p}(G)}^{p} \\ \notag
&\leq & \frac{1}{\Phi ^{p}\left( \alpha _{k}\right) }\sum_{\eta =0}^{\infty }%
\frac{\Phi ^{p/2}\left( \alpha _{\eta }\right) }{2^{d\left( \alpha _{\eta
		}\right) \left( 1-p\right) /2}}\leq c<\infty.
\end{eqnarray}

Let $x\in I_{\left\langle \alpha _{k}\right\rangle }\backslash
I_{\left\langle \alpha _{k}\right\rangle +1}.$  Since
$\left\vert \alpha _{k}\right\vert \neq\left\langle \alpha _{k}\right\rangle$ and
$\left\langle
\alpha _{k}-2^{\left\vert \alpha _{k}\right\vert }\right\rangle
=\left\langle \alpha _{k}\right\rangle.$

By using both inequalities of Lemma  \ref{lemma1} we get that
\begin{eqnarray} \label{4.1.77}
&&\left\vert D_{\alpha _{k}-2^{\left\vert \alpha _{k}\right\vert }}(x)\right\vert \\ \notag
&=&\left\vert \left( D_{2^{\left\langle \alpha _{k}\right\rangle
		+1}}(x)-D_{2^{\left\langle \alpha _{k}\right\rangle }}(x)\right) +\overset{\left\vert \alpha _{k}\right\vert -1}{\underset{j=\left\langle \alpha_{k}\right\rangle +1}{\sum }}\left( \alpha _{k}\right)_{j}\left(
D_{2^{i+1}}(x)-D_{2^{i}}(x)\right) \right\vert \\ \notag
&=&\left\vert -D_{2^{\left\langle
		\alpha _{k}\right\rangle }}(x)\right\vert =2^{\left\langle \alpha
	_{k}\right\rangle }
\end{eqnarray}%
and
\begin{eqnarray} \label{4.1.12aa}
\left\vert II\right\vert =\frac{2^{\left( \left\vert \alpha _{k}\right\vert
		+\left\langle \alpha _{k}\right\rangle \right) \left( 1/p-1\right) /2}}{\Phi
	^{1/2}\left( \alpha _{k}\right) }\left\vert D_{\alpha _{k}-2^{\left\vert
		\alpha _{k}\right\vert }}\left( x\right) \right\vert =\frac{2^{\left\vert \alpha _{k}\right\vert \left( 1/p-1\right)
		/2}2^{\left\langle \alpha _{k}\right\rangle \left( 1/p+1\right) /2}}{\Phi
	^{1/2}\left( \alpha _{k}\right) }.
\end{eqnarray}

By combining (\ref{4.1.8a}) and (\ref{4.1.12aa}) we get that
\begin{eqnarray*}
	&&\left\Vert \frac{S_{\alpha _{k}}f}{\Phi \left( \alpha _{k}\right) }%
	\right\Vert _{weak-L_{p}(G)}^{p} \\ \notag
	&\geq & \left\Vert II\right\Vert _{weak-L_{p}(G)}^{p}-\left\Vert I\right\Vert _{weak-L_{p}(G)}^{p} \\
	&\geq &\frac{2^{\left( \left\vert \alpha _{k}\right\vert \right) \left(
			1/p-1\right) /2}2^{\left\langle \alpha _{k}\right\rangle \left( 1/p+1\right)
			/2}}{\Phi ^{1/2}\left( \alpha _{k}\right) }\mu \left\{ x\in G:\text{ }%
	\left\vert II\right\vert \geq \frac{2^{\left( \left\vert \alpha
			_{k}\right\vert \right) \left( 1/p-1\right) /2}2^{\left\langle \alpha
			_{k}\right\rangle \left( 1/p+1\right) /2}}{\Phi ^{1/2}\left( \alpha
		_{k}\right) }\right\} ^{1/p} \\
	&\geq & \frac{2^{\left( \left\vert \alpha _{k}\right\vert \right) \left(
			1/p-1\right) /2}2^{\left\langle \alpha _{k}\right\rangle \left( 1/p+1\right)
			/2}}{\Phi ^{1/2}\left( \alpha _{k}\right) }\left( \mu \left\{
	I_{\left\langle \alpha _{k}\right\rangle }\backslash I_{\left\langle \alpha
		_{k}\right\rangle +1}\right\} \right) ^{1/p} \\
	&\geq & c\frac{2^{d\left( \alpha _{k}\right) \left( 1/p-1\right) /2}}{\Phi
		^{1/2}\left( \alpha _{k}\right) }\rightarrow \infty, \ \text{ \ as \  }
	k\rightarrow \infty.
\end{eqnarray*}

The proof of Theorem \ref{th4.1.1} is complete.

\begin{corollary} \label{cor4.1.1}
	a) Let $n\in \mathbb{N}_{+}$, $0<p<1$ and $f\in H_{p}(G)$. Then there exists an absolute constant  $c_{p},$ depending only on $p$ such that
	\begin{equation*}
	\text{ }\left\Vert S_{n}f\right\Vert _{H_{p}(G)}\leq c_{p}\left( n\mu \left\{
	\text{supp}\left( D_{n}\right) \right\} \right) ^{1/p-1}\left\Vert
	f\right\Vert _{H_{p}(G)}.
	\end{equation*}
	
	\textit{b) Let} $0<p<1,$ $\left\{ m_{k}:\text{ }k\in \mathbb{N}_+\right\} $  be increasing sequence of natural numbers, such that
	\begin{equation}\label{suppdnk}
	\sup_{k\in \mathbb{N}}m_{k}\mu \left\{ \text{supp}\left( D_{m_{k}}\right)
	\right\} =\infty
	\end{equation}%
	and
	$\Phi :\mathbb{N}_{+}\rightarrow \lbrack 1,\infty )$ be non-decreasing function satisfying the condition
	\begin{equation}
	\overline{\underset{k\rightarrow \infty }{\lim }}\frac{\left( m_{k}\mu
		\left\{ \text{supp}\left( D_{m_{k}}\right) \right\} \right) ^{1/p-1}}{\Phi
		\left(m_{k}\right) }=\infty .  \label{4.1.12e}
	\end{equation}%
	\textit{Then there exists a martingale} $f\in H_{p}(G)$ \textit{such that}
	\begin{equation*}
	\underset{k\in \mathbb{N}}{\sup }\left\Vert \frac{S_{m_{k}}f}{\Phi \left(
		m_{k}\right) }\right\Vert _{weak-L_{p}(G)}=\infty.
	\end{equation*}
\end{corollary}

{\bf Proof}:
By applying both inequalities of Lemma \ref{lemma1} we get that
\begin{equation*}
I_{\left\langle n\right\rangle }\backslash I_{\left\langle n\right\rangle
	+1}\subset \text{supp}\left\{ D_{n}\right\} \subset I_{\left\langle
	n\right\rangle }\text{ \ and \ }2^{-\left\langle n\right\rangle -1}\leq \mu
\left\{ \text{supp}\left( D_{n}\right) \right\} \leq 2^{-\left\langle
	n\right\rangle }.
\end{equation*}

Hence,
\begin{equation*}
\frac{2^{d\left( n\right) \left( 1/p-1\right) }}{4}\leq \left( n\mu \left\{
\text{supp}\left( D_{n}\right) \right\} \right) ^{1/p-1}\leq 2^{d\left(
	n\right) \left( 1/p-1\right) }.
\end{equation*}%

Corollary \ref{cor4.1.1} is proved.

\begin{theorem} \label{th4.1.2}
	a) Let $n\in \mathbb{N}_{+}$ and $f\in H_{1}(G).$ Then there exists an absolute constant $c,$ such that
	\begin{equation*}
	\left\Vert S_{n}f\right\Vert _{H_{1}(G)}\leq cV\left( n\right) \left\Vert f\right\Vert _{H_{1}(G)}.
	\end{equation*}
	
	\textit{b) Let} $\left\{ m_{k}:\text{ }k\in \mathbb{N}_{+}\right\} $ be non-negative increasing  sequence of natural numbers such that
	\begin{equation} \label{4.1.vnk}
	\sup_{k\in \mathbb{N}}V\left( m_{k}\right) =\infty
	\end{equation}%
	and $\Phi :\mathbb{N}_{+}\rightarrow \lbrack 1,\infty )$ be non-decreasing function satisfying the condition
	\begin{equation}  \label{4.1.17aa}
	\overline{\underset{k\rightarrow \infty }{\lim }}\frac{V\left( m_{k}\right)
	}{\Phi \left( m_{k}\right) }=\infty .
	\end{equation}
	
	Then there exists a martingale $f\in H_{1}(G),$ such that
	\begin{equation*}
	\underset{{k\in \mathbb{N}}}{\sup }\left\Vert \frac{S_{m_{k}}f}{\Phi \left(
		m_{k}\right) }\right\Vert _{1}=\infty.
	\end{equation*}
\end{theorem}

{\bf Proof}:
Since
\begin{equation} \label{4.1.12k}
\left\Vert \frac{S_{n}f}{V\left( n\right)}\right\Vert _{1}\leq \left\Vert
f\right\Vert _{1}\leq \left\Vert f\right\Vert _{H_{1}(G)}
\end{equation}
by combining Lemmas \ref{lemma3.2.8} and (\ref{4.1.12k}) we can conclude that
\begin{eqnarray} \label{4.1.50}
\left\Vert \frac{S_{n}f}{V\left( n\right) }\right\Vert _{H_{1}(G)}
&\leq & \left\Vert \frac{S_{n}f}{V\left( n\right) }\right\Vert_{1}+\frac{1}{V\left( n\right) }\left\Vert\widetilde{S}_{\#}^{\ast}f\right\Vert_{1} \\ \notag
&\leq & c\left\Vert f\right\Vert _{H_{1}(G)}+c\left\Vert\widetilde{S}_{\#}^{\ast}f\right\Vert
_{1}\leq c\left\Vert f\right\Vert _{H_{1}(G)}.
\end{eqnarray}

Now prove second part of Theorem \ref{th4.1.2}. Let $\left\{ m_{k}:\text{ }k\in \mathbb{N}_{+}\right\} $ be increasing sequence of natural numbers and function $\ \Phi :\mathbb{N}_{+}\rightarrow \lbrack 1,\infty )$ satisfies conditions (\ref{4.1.vnk}) and (\ref{4.1.17aa}). Then there exists  non-negative, increasing sequence  $\left\{ \alpha _{k}:k\in \mathbb{N}_{+} \right\} \subset\left\{ m_{k}:k \in\mathbb{N}_{+}\right\} $ such that
\begin{equation} \label{4.1.2aaa}
\sum_{k=1}^{\infty}\frac{\Phi ^{1/2}\left(\alpha_{k}\right)}{
	V^{1/2}\left(\alpha_{k}\right)}\leq \beta <\infty.
\end{equation}

Let  $f=\left(f_{n},n\in \mathbb{N}_+\right)$ be martingale from Example \ref{example2.2.1}, where
\begin{equation} \label{4.1.101}
\text{\ }\lambda _{k}=\frac{\Phi ^{1/2}\left( \alpha _{k}\right) }{%
	V^{1/2}\left( \alpha _{k}\right) }.
\end{equation}

By applying condition (\ref{4.1.2aaa}) we can conclude that condition (\ref{3.3.2aa}) is fulfilled and it follows that $f=\left(f_{n},n\in \mathbb{N}_+\right) \in H_{1}(G).$

In the view of (\ref{3.3.10AA}) when $\lambda _{k}$ are given by (\ref{4.1.101}) we get that
\begin{equation} \label{4.1.5aa}
\widehat{f}(j)=\left\{
\begin{array}{ll}
\frac{\Phi ^{1/2}\left( \alpha _{k}\right) }{V^{1/2}\left( \alpha
	_{k}\right) }, & \text{ \  if \ }j\in \left\{
2^{_{\left\vert \alpha _{k}\right\vert }},...,2^{_{\left\vert \alpha
		_{k}\right\vert +1}}-1\right\} ,\text{ }k=0,1,... \\
0, & \text{ \ if \ }j\notin
\bigcup\limits_{k=0}^{\infty }\left\{ 2^{_{\left\vert \alpha _{k}\right\vert
}},...,2^{_{\left\vert \alpha _{k}\right\vert +1}}-1\right\} .\text{ }%
\end{array}%
\right.
\end{equation}

Analogously to (\ref{4.1.6aaa}) if we apply (\ref{3.3.11AA}) when $\lambda _{k}$ are given by (\ref{4.1.101}) we get that
\begin{eqnarray*}
	S_{\alpha _{k}}f = \sum_{\eta =0}^{k-1}\frac{\Phi ^{1/2}\left( \alpha _{\eta
		}\right) }{V^{1/2}\left( \alpha _{\eta }\right) }\left( D_{2^{\left\vert
			\alpha _{\eta }\right\vert +1}}-D_{2^{\left\vert \alpha _{\eta }\right\vert
	}}\right)+\frac{\Phi ^{1/2}\left( \alpha _{k}\right) }{V^{1/2}\left( \alpha
		_{k}\right) }w_{2^{\left\vert \alpha _{k}\right\vert }}D_{\alpha
		_{k}-2^{\left\vert \alpha _{k}\right\vert }}.
\end{eqnarray*}

By applying first estimation of Lemma \ref{lemma2} and (\ref{4.1.2aaa}) we can conclude that
\begin{eqnarray*}
	\left\Vert \frac{S_{\alpha _{k}}f}{\Phi \left( \alpha _{k}\right) }%
	\right\Vert _{1}
	&\geq &\frac{\Phi ^{1/2}\left( \alpha _{k}\right) }{\Phi
		\left( \alpha _{k}\right) V^{1/2}\left( \alpha _{k}\right) }\left\Vert
	D_{\alpha _{k}-2^{\left\vert \alpha _{k}\right\vert }}\right\Vert _{1} \\
	&-&\frac{1}{\Phi \left( \alpha _{k}\right) }\sum_{\eta =0}^{k-1}\frac{\Phi
		^{1/2}\left( \alpha _{\eta }\right) }{V^{1/2}\left( \alpha _{\eta }\right) }\left\Vert  D_{2^{\left\vert
			\alpha _{\eta }\right\vert +1}}-D_{2^{\left\vert \alpha _{\eta }\right\vert}}\right\Vert _{1} \\
	&\geq &\frac{V\left( \alpha _{k}-2^{\left\vert \alpha _{k}\right\vert }\right)
		\Phi ^{1/2}\left( \alpha _{k}\right) }{8\Phi \left( \alpha _{k}\right)
		V^{1/2}\left( \alpha _{k}\right) } \\
	&-&\frac{1}{\Phi \left( \alpha _{k}\right) }%
	\sum_{\eta =0}^{\infty }\frac{\Phi ^{1/2}\left( \alpha _{\eta }\right) }{%
		V^{1/2}\left( \alpha _{\eta }\right) } \\
	&\geq & \frac{cV^{1/2}\left( \alpha _{k}\right) }{\Phi ^{1/2}\left( \alpha
		_{k}\right) }\rightarrow \infty ,\text{ \ as \ }k\rightarrow \infty.
\end{eqnarray*}

Theorem \ref{th4.1.2} is proved.

\begin{corollary} \label{cor4.1.2}
	Let $n\in \mathbb{N}$, $0<p\leq 1$ and $f\in H_{p}(G)$. Then there exists an absolute constant $ c_{p}, $ depending only on $ p, $ such that
	\begin{equation} \label{4.1.sn2n20}
	\left\Vert S_{2^{n}}f\right\Vert _{H_{p}(G)}\leq c_{p}\left\Vert
	f\right\Vert _{H_{p}(G)}.
	\end{equation}
\end{corollary}

{\bf Proof}:
To prove Theorem \ref{cor4.1.2} we only have to show that
\begin{equation*}
\left\vert 2^{n}\right\vert =n,\text{ \  \ } \left\langle
2^{n}\right\rangle =n-1\text{ \ and \ }d\left(
2^{n}\right)=0.
\end{equation*}
By applying first part of Theorem \ref{th4.1.1} we immediately get that (\ref{4.1.sn2n20}) for any $0<p\leq 1$ and proof of Corollary \ref{cor4.1.2} is proved.

\begin{corollary} \label{cor4.1.3}
	Let $n\in \mathbb{N}$, $0<p\leq 1$ and $f\in H_{p}(G)$. Then there exists an absolute constant $ c_p, $  depending only on $p$ such that 
	\begin{equation} \label{4.1.sn2n2}
	\left\Vert S_{2^{n}+2^{n-1}}f\right\Vert _{H_{p}(G)}\leq c_{p}\left\Vert
	f\right\Vert _{H_{p}(G)}.
	\end{equation}
\end{corollary}

{\bf Proof}:
Since
\begin{equation*}
\left\vert 2^{n}+2^{n-1}\right\vert =n,\left\langle
2^{n}+2^{n-1}\right\rangle =n-1\text{ \ and \ }d\left(
2^{n}+2^{n-1}\right)=1
\end{equation*}
by first part of Theorem \ref{th4.1.1} we get that  (\ref{4.1.sn2n2}) holds, for any $0<p\leq 1$ and proof of Corollary \ref{cor4.1.3} is complete.

\begin{corollary} \label{cor4.1.4}
	Let $n\in \mathbb{N}$ and $0<p<1.$ Then there exists a martingale $f\in
	H_{p}(G),$ such that
	\begin{equation} \label{4.1.sn2n1}
	\underset{n\in \mathbb{N}}{\sup }\left\Vert S_{2^{n}+1}f\right\Vert
	_{weak-L_{p}(G)}=\infty.
	\end{equation}
	
	On the other hand, there exists an absolute constant $ c $, such that
	
	\begin{equation} \label{4.1.sn2n22}
	\left\Vert S_{2^{n}+1}f\right\Vert _{H_{1}(G)}\leq c\left\Vert
	f\right\Vert _{H_{1}(G)}.
	\end{equation}
	
\end{corollary}

{\bf Proof}:
Since
\begin{equation}
\left\vert 2^{n}+1\right\vert =n,\left\langle 2^{n}+1\right\rangle =0\text{
	\ and \ }d\left( 2^{n}+1\right) =n.  \label{4.1.cor1}
\end{equation}%
by applying second part of Theorem \ref{th4.1.1} we get that there exists a martingale $f=\left(f_{n},n\in \mathbb{N}_+\right) \in H_{p}(G),$ for $0<p<1,$ such that (\ref{4.1.sn2n1}) holds.

On the other hand, proof of  (\ref{4.1.sn2n22}) follows simple observation that $$ V( 2^{n}+1)=4<\infty. $$

Corollary \ref{cor4.1.4} is proved.

\bigskip

\subsection{Modulus of continuity and convergence in norm of subsequences of partial sums with respect to the one-dimensional Walsh-Fourier series on the martingale Hardy spaces} 
\text{ \qquad \qquad \qquad \qquad \qquad \qquad \qquad \qquad \qquad \qquad \qquad \qquad \qquad \qquad \qquad  }
 In this section we apply Theorem \ref{th4.1.1} and  Theorem \ref{th4.1.2} to find necessary and sufficient conditions for modulus of continuity, for which subsequences of partial sums with respect to the one-dimensional Walsh-Fourier series are bounded in the martingale Hardy spaces.

First, we prove the following estimation:

\begin{theorem} \label{theorem4.2.1}
	Let $ n\in \mathbb{N}_+ $ and $2^{k}<n\leq 2^{k+1}.$ Then there exists an absolute constant $c_{p},$ depending only on $p$ such that
	\begin{equation} \label{4.2.sn1}
	\left\Vert S_{n}f-f\right\Vert _{H_{p}(G)}\leq c_{p}2^{d\left( n\right) \left(1/p-1\right)}\omega _{H_{p}(G)}\left( \frac{1}{2^{k}},f\right),
	\text{ \ \ \ } (f\in H_p(G))\text{ \ \ \ }\left( 0<p<1\right)
	\end{equation}%
	and
	\begin{equation} \label{4.2.sn2}
	\left\Vert S_{n}f-f\right\Vert _{H_{1}(G)}\leq c_{1}V\left( n\right) \omega
	_{H_{1}(G)}\left( \frac{1}{2^{k}},f\right),\text{ \ \ \ } (f\in H_1(G)) .
	\end{equation}
\end{theorem}

{\bf Proof}:
Let $0<p<1$ and $2^{k}<n\leq 2^{k+1}.$ By applying first part of Theorem \ref{th4.1.1} we get that
\begin{eqnarray} \label{4.2.1000}
&&\left\Vert S_{n}f-f\right\Vert _{H_{p}(G)}^p \\ \notag
&\leq& c_{p}\left\Vert
S_{n}f-S_{2^{k}}f\right\Vert _{H_{p}(G)}^p+c_{p}\left\Vert
S_{2^{k}}f-f\right\Vert _{H_{p}(G)}^p  \\ \notag
&=&c_{p}\left\Vert S_{n}\left( S_{2^{k}}f-f\right) \right\Vert
_{H_{p}(G)}^p+c_{p}\left\Vert S_{2^{k}}f-f\right\Vert _{H_{p}(G)}^p \\ \notag
&\leq & c_{p}\left( 1+2^{d\left( n\right) \left( 1-p\right) }\right) \omega
_{H_{p}(G)}^p\left( \frac{1}{2^{k}},f\right) \\ \notag
&\leq & c_{p}2^{d\left( n\right) \left(
	1-p\right) }\omega _{H_{p}(G)}^p\left( \frac{1}{2^{k}},f\right).
\end{eqnarray}

The proof of (\ref{4.2.sn2}) is analogical to (\ref{4.2.sn1}). Analogously to (\ref{4.2.sn1}) we can also prove estimation (\ref{4.2.sn2}). So, we leave out the details.

Theorem \ref{theorem4.2.1} is proved.

\begin{theorem} \label{th4.2.2}
	a) Let $k\in\mathbb{N}_+$, $0<p<1,$ $f\in H_{p}(G)$ and $\{m_{k}:k\in\mathbb{N}_+\}$ be increasing sequence of natural numbers, such that
	\begin{equation}  \label{4.2.18a}
	\omega _{H_{p}(G)}\left( \frac{1}{2^{\left\vert m_{k}\right\vert}},f\right)
	=o\left(\frac{1}{2^{d\left(m_{k}\right)\left( 1/p-1\right)}}\right)
	\text{ \ as \ }k\rightarrow \infty.
	\end{equation}%
	Then
	\begin{equation} \label{4.2.con1}
	\left\Vert S_{m_{k}}f-f\right\Vert _{H_{p}(G)}\rightarrow 0\text{ \ as \ }%
	k\rightarrow \infty .
	\end{equation}
	
	b) Let $\{m_{k}:k\in\mathbb{N}_+\}$ be increasing sequence of natural numbers,such that condition (\ref{4.1.dnk}) is fulfilled. Then there exists a martingale  $f\in H_{p}(G)$ and increasing sequence of natural numbers  $\{\alpha _{k}:k\in\mathbb{N}_+\}\subset \{m_{k}:k\in\mathbb{N}_+\},$ such that
	\begin{equation*}
	\omega _{H_{p}(G)}\left( \frac{1}{2^{\left\vert \alpha _{k}\right\vert }}%
	,f\right) =O\left( \frac{1}{2^{d\left( \alpha _{k}\right) \left(
			1/p-1\right) }}\right) \text{ \ as \ }k\rightarrow \infty \text{\ }
	\end{equation*}%
	and
	\begin{equation} \label{4.2.con11}
	\limsup\limits_{k\rightarrow \infty }\left\Vert S_{\alpha
		_{k}}f-f\right\Vert _{weak-L_{p}(G)}>c_{p}>0,\text{ \ as  \ }k\rightarrow \infty,
	\end{equation}%
	where $c_{p}$ is an absolute constant depending only on  $p$.
\end{theorem}

{\bf Proof}:
Let $0<p<1$, $f\in H_{p}(G)$ and $\{m_{k}:k\in\mathbb{N}_+\}$ be increasing sequence of natural numbers, such that condition  (\ref{4.2.18a}) is fulfilled. By combining Theorem \ref{theorem4.2.1} and estimation (\ref{4.2.sn1}) we get that (\ref{4.2.con1}) holds true.

Now, prove second part of Theorem \ref{th4.2.2}. In the view of  (\ref{4.1.dnk}) we simply get that there exists sequence $\{\alpha _{k}:k\in\mathbb{N}_+\}\subset \{m_{k}:k\in\mathbb{N}_+\},$ such that
\begin{equation} \label{4.2.4.18}
\ 2^{d\left( \alpha _{k}\right) }\uparrow \infty ,\,\,\,\,\text{ \ as\ }\ \
k\rightarrow \infty \text{, \ \ \ }2^{2\left( 1/p-1\right) d\left( \alpha
	_{k}\right) }\leq 2^{\left( 1/p-1\right) d\left( \alpha _{k+1}\right) }.
\end{equation}%

Let  $f=\left(f_{n},n\in \mathbb{N}\right)$ be a martingale from Example \ref{example2.2.1}, such that
\begin{equation} \label{4.2.4.188}
\lambda_{i}={2^{-(1/p-1)d\left(\alpha_{i}\right)}}.
\end{equation}

By applying (\ref{4.2.4.18}) we obtain that condition (\ref{3.3.2aa}) is fulfilled and it follows that  $\ f\in H_{p}(G).$

By applying (\ref{3.3.10AA}), whene $\lambda _{k}$ are given by (\ref{4.2.4.188}), then
\begin{equation} \label{4.2.4.22}
\widehat{f}(j)=\left\{
\begin{array}{ll}
2^{\left( 1/p-1\right) \left\langle \alpha _{k}\right\rangle }, & \text{ \ if
	\ }j\in \left\{ 2^{_{\left\vert \alpha _{k}\right\vert
}},...,2^{_{\left\vert \alpha _{k}\right\vert +1}}-1\right\} ,\text{ }%
k\in\mathbb{N}_+, \\
0, & \text{\ if \ }j\notin
\bigcup\limits_{k=0}^{\infty }\left\{ 2^{_{\left\vert \alpha _{k}\right\vert
}},...,2^{_{\left\vert \alpha _{k}\right\vert +1}}-1\right\} .\text{ }%
\end{array}
\right.
\end{equation}

By combining (\ref{4.2.4.18}) and (\ref{3.3.2aa0}) we have that
\begin{eqnarray} \label{4.2.4.21}
\omega _{H_{p}(G)}(\frac{1}{2^{\left\vert \alpha _{k}\right\vert }},f) \leq\sum\limits_{i=k}^{\infty }\frac{1}{2^{\left( 1/p-1\right) d\left( \alpha_{i}\right) }} = O\left( \frac{1}{2^{\left( 1/p-1\right) d\left( \alpha_{k}\right) }}\right) ,\text{ \ as \ }k\rightarrow \infty.
\end{eqnarray}

By using (\ref{4.1.77}) we get that
\begin{equation*}
\left\vert D_{\alpha _{k}-2^{\left\langle \alpha
		_{k}\right\rangle }}\right\vert \geq 2^{\left\langle \alpha
	_{k}\right\rangle },\text{ \ \ \ \ where \ \ \ }I_{\left\langle \alpha_{k}\right\rangle }\backslash I_{\left\langle \alpha _{k}\right\rangle +1}.
\end{equation*}

In the view of (\ref{3.3.11AA}) we can conclude that
\begin{equation*}
S_{\alpha _{k}}f=S_{2^{|\alpha _{k}|}}f+2^{\left( 1/p-1\right) \left\langle \alpha _{k}\right\rangle }{w_{2^{|\alpha _{k}|}}D_{{\alpha _{k}}-2^{|\alpha _{k}|}}},
\end{equation*}

Since
\begin{eqnarray*}
	\Vert D_{\alpha _{k}}\Vert _{weak-L_{p}(G)}
	&\geq & 2^{\left\langle \alpha
		_{k}\right\rangle }\mu \left\{ x\in I_{\left\langle \alpha _{k}\right\rangle} \backslash I_{\left\langle \alpha _{k}\right\rangle +1}:\text{ }\left \vert D_{\alpha _{k}}\right\vert \geq 2^{\left\langle \alpha _{k}\right\rangle}\right\} ^{1/p} \\
	&\geq & 2^{\left\langle \alpha _{k}\right\rangle }\left( \mu \left\{
	I_{\left\langle \alpha _{k}\right\rangle }\backslash I_{\left\langle \alpha_{k}\right\rangle +1}\right\} \right) ^{1/p}\geq 2^{\left\langle \alpha_{k}\right\rangle \left(1-1/p\right) },
\end{eqnarray*}
if we apply (\ref{1.S2n000}) (see also Theorem T2) we obtain that
\begin{eqnarray*}
	\Vert f-S_{\alpha _{k}}f\Vert_{weak-L_{p}(G)}^p 
	&\geq& 2^{\left( 1-p\right) \left\langle \alpha _{k}\right\rangle }\Vert {w_{2^{|\alpha _{k}|}}D_{{\alpha _{k}}-2^{|\alpha _{k}|}}}\Vert _{weak-L_{p}(G)}^p \\
	&-&\Vert f-S_{2^{|\alpha _{k}|}}f\Vert_{weak-L_{p}(G)}^p \\
	&\geq & c-o(1)>c>0, \text{\ \ \ as  \ \ }k\rightarrow \infty .
\end{eqnarray*}
Proof of Theorem \ref{th4.2.2} is complete.

\begin{corollary} \label{cor4.2.1}
	a) Let $0<p<1,$ $f\in H_{p}(G)$ and $\{m_{k}:k\in \mathbb{N}_+\}$ be increasing sequence of natural numbers, such that
	\begin{equation} \label{4.2.cond2}
	\omega_{H_{p}(G)}\left(\frac{1}{2^{\left\vert m_{k}\right\vert }},f\right)
	=o\left( \frac{1}{\left(m_{k}\mu \left( \text{supp}D_{m_{k}}\right) \right)
		^{1/p-1}}\right) \text{ \ as\ }k\rightarrow \infty .
	\end{equation}%
	Then (\ref{4.2.con1}) holds.
	
	b) Let $\{m_{k}:k\in \mathbb{N}_+\}$ be increasing sequence of natural numbers, such that
	\begin{equation} \label{4.2.con1ab}
	\sup_{k\in \mathbb{N}_+}m_{k}\mu \left\{ \text{supp}\left( D_{m_{k}}\right)
	\right\} =\infty.
	\end{equation}
	Then there exist a martingale $f\in H_{p}(G)$ and sequence $\{\alpha _{k}:k\in \mathbb{N}_+\}\subset \{m_{k}:k\in \mathbb{N}_+\},$ such that
	\begin{equation*}
	\omega _{H_{p}(G)}\left( \frac{1}{2^{\left\vert \alpha _{k}\right\vert }},f\right) =O\left( \frac{1}{\left( \alpha _{k}\mu \left( \text{supp}D_{\alpha _{k}}\right) \right)^{1/p-1}}\right) \text{\ \ as \ }k\rightarrow\infty
	\end{equation*}%
	and (\ref{4.2.con11}) holds.
\end{corollary}

\begin{theorem} \label{th4.2.3}
	a) Let $f\in H_{1}(G)$ and  $\{m_{k}:k\in \mathbb{N}_+\}$ be increasing sequence of natural numbers, such that
	\begin{equation} \label{4.2.cond1}
	\omega _{H_{1}(G)}\left( \frac{1}{2^{\left\vert m_{k}\right\vert }},f\right)=o\left( \frac{1}{V\left( m_{k}\right) }\right) \text{ \ as \ }k\rightarrow \infty.
	\end{equation}%
	Then
	\begin{equation} \label{4.2.cond1a}
	\left\Vert S_{m_{k}}f-f\right\Vert _{H_{1}(G)}\rightarrow 0\text{ \ as \ }%
	k\rightarrow \infty .
	\end{equation}
	
	b) Let $\{m_{k}:k\in \mathbb{N}_+\}$ be increasing sequence of natural numbers, such that condition (\ref{4.1.vnk}) is fulfilled. Then there exists a martingale $f\in H_{1}(G)$ and increasing sequence of natural numbers
	$\{\alpha _{k}:k\in \mathbb{N}_+\}\subset \{m_{k}:k\in \mathbb{N}_+\}$ such that
	\begin{equation*}
	\omega _{H_{1}(G)}\left( \frac{1}{2^{\left\vert \alpha _{k}\right\vert }},f\right) =O\left( \frac{1}{V\left( \alpha _{k}\right) }\right) \text{ \ as
		\ }k\rightarrow \infty
	\end{equation*}%
	and
	\begin{equation} \label{cond10}
	\limsup\limits_{k\rightarrow \infty }\left\Vert S_{\alpha
		_{k}}f-f\right\Vert _{1}>c>0\,\,\,\text{\ as \ }
	k\rightarrow \infty,
	\end{equation}
	where $c$ is an absolute constant.
\end{theorem}

{\bf Proof}:
Let $f\in H_{1}(G)$ and $\{m_{k}:k\in \mathbb{N}_+\}$  be increasing sequence of natural numbers, such that (\ref{4.2.cond1}). By applying Theorem \ref{theorem4.2.1} we get that condition (\ref{4.2.cond1a}) is fulfilled.

Now, we prove second part of Theorem \ref{th4.2.3}. By applying (\ref{4.1.vnk}) we conclude that there exists sequence  $\{\alpha _{k}:k\in \mathbb{N}_+\}\subset
\{m_{k}:k\in \mathbb{N}_+\}$, such that
\begin{equation}
V(\alpha _{k})\uparrow \infty , \text{ \ \ as \ \ }k\rightarrow \infty \text{ \ \ and \ \ }V^{2}(\alpha _{k})\leq V(\alpha _{k+1}) \text{ \ \ }k\in \mathbb{N}_+.  \label{4.2.4.7}
\end{equation}

Let  $f=\left(f_{n},n\in \mathbb{N}_+\right)$ be a martingale from the Example \ref{example2.2.1}, where
\begin{equation*}
\lambda_{k}=\frac{1}{V(\alpha _{k})}.
\end{equation*}%

By applying (\ref{4.2.4.7}) we conclude that (\ref{3.3.2aa}) is fulfilled and we conclude that $f=\left( f_{n},n\in
\mathbb{N}_+\right) \in H_{1}(G).$

In the view of (\ref{3.3.10AA}) we have that
\begin{equation}
\widehat{f}(j)=\left\{
\begin{array}{ll}
\frac{1}{V(\alpha _{k})}, & \text{ if \thinspace \thinspace }j\in \left\{
2^{_{\left\vert \alpha _{k}\right\vert }},...,2^{_{\left\vert \alpha
		_{k}\right\vert +1}}-1\right\} ,\text{ }k=0,1,... \\
0, & \text{\ if \thinspace \thinspace \thinspace }j\notin
\bigcup\limits_{k=0}^{\infty }\left\{ 2^{_{\left\vert \alpha _{k}\right\vert}},...,2^{_{\left\vert \alpha _{k}\right\vert +1}}-1\right\} .\text{ }%
\end{array}%
\right.  \label{4.13}
\end{equation}

According to (\ref{3.3.2aa0}) we get that
\begin{eqnarray*} 
	w_{H_{1}(G)}(1/2^n,f) 
	&=&\left\Vert f-S_{2^{n}}f\right\Vert _{H_{1}(G)}
	\leq \sum\limits_{i=n+1}^{\infty }%
	\frac{1}{V(\alpha _{i})} \\ \notag
	&=& O\left( \frac{1}{V(\alpha _{n})}\right), \text{\ \ as \ \ }n\rightarrow \infty.  
\end{eqnarray*}

By applying (\ref{3.3.11AA}) we can conclude that
\begin{equation*} \label{4.2.4.1222222}
S_{\alpha _{k}}f=S_{2^{|\alpha _{k}|}}f+\frac{w_{2^{|\alpha _{k}|}}D_{{\alpha _{k}}-2^{|\alpha _{k}|}}}{V(\alpha _{k})},
\end{equation*}

If we use (\ref{1.S2n000}) and Theorem T2  we get that
\begin{eqnarray*}
	\Vert f-S_{\alpha _{k}}f\Vert _{1}
	&\geq & \Vert \frac{w_{2^{|\alpha _{k}|}}D_{{\alpha _{k}}-2^{|\alpha _{k}|}}}{V(\alpha _{k})}\Vert _{1}-\Vert f-S_{2^{|\alpha _{k}|}}f\Vert _{1} \\
	&\geq & \frac{V(\alpha _{k}-2^{|\alpha _{k}|})}{8V(\alpha _{k})}-o(1)>c>0, \text{\ \ as \ \ }k\rightarrow \infty .
\end{eqnarray*}

The proof of Theorem \ref{th4.2.3} is proved.

Theorem \ref{theorem5.2.2} follows the following corollaries which are \cite{tep7}:

\begin{corollary}\label{corollary5.2.2snsn}
	a) Let $0<p<1,$ $f\in H_{p}(G)$ and
	\begin{equation*}
	\omega _{H_{p}(G)}\left( 1/2^{k},f\right) =o\left(
	1/2^{k( 1/p-1)}\right) ,\text{ as  }
	k\rightarrow \infty .
	\end{equation*}%
	Then
	\begin{equation*}
	\left\Vert S_{k}f-f\right\Vert _{H_{p}(G)}\rightarrow
	0, \text{ \ as \ }k\rightarrow \infty .
	\end{equation*}
	
	b) There exists a martingale $f\in H_{p}(G)$ $\left( 0<p<1\right),$ such that
	\begin{equation*}
	\omega _{H_{p}(G)}\left( 1/2^{k},f\right) =O\left(
	1/2^{k\left( 1/p-1\right) }\right),\text{ \ as \ }%
	k\rightarrow \infty
	\end{equation*}%
	and
	\begin{equation*}
	\left\Vert S_{k}f-f\right\Vert _{weak-L_{p}(G)}\nrightarrow
	0, \text{ \ as \ }k\rightarrow\infty .
	\end{equation*}
\end{corollary}

\begin{corollary}\label{corollary5.2.2snsn2}
	a) Let  $f\in H_{1}(G)$ and
	\begin{equation*}
	\omega _{H_{1}(G)}\left( 1/2^{k},f\right) =o\left(\frac{1}{ k}\right) ,\text{ as \ }
	k\rightarrow \infty .
	\end{equation*}%
	Then
	\begin{equation*}
	\left\Vert S_{k}f-f\right\Vert _{H_{1}(G)}\rightarrow
	0,\text{\ as \ }k\rightarrow\infty .
	\end{equation*}
	
	b) There exists a martingale $f\in H_{1}(G),$ such that
	\begin{equation*}
	\omega _{H_{1}(G)}\left( 1/2^{k},f\right)=O\left(\frac{1}{k}\right),\text{ \ as \ }k\rightarrow \infty
	\end{equation*}
	and
	\begin{equation*}
	\left\Vert S_kf-f\right\Vert_1\nrightarrow 0,\text{ \ as \ }k\rightarrow \infty.
	\end{equation*}
\end{corollary}
\newpage

\section{Fej\'er means with respect to the one-dimensional Walsh-Fourier series on the martingale Hardy spaces}
\subsection{Basic notations} 
\text{ \qquad \qquad \qquad \qquad \qquad \qquad \qquad \qquad \qquad \qquad \qquad \qquad \qquad \qquad \qquad  } For the one-dimensional case Fej\'er means with respect to the one-dimensional Walsh-Fourier series $\sigma _{n}$ is defined by:

\begin{eqnarray*}
	\qquad \sigma _{n}f(x) &:&=\frac{1}{n}\sum_{k=1}^{n}S_{k}f(x)\text{\qquad\ \ \ }\left( \text{ }n\in \mathbb{N}_{+}\right).
\end{eqnarray*}

The following equality is true (for details see \cite{1} and \cite{sws}):

\begin{equation*}
\sigma _{n}f\left( x\right) =\frac{1}{n}\overset{n-1}{\underset{k=0}{\sum }}%
\left( D_{k}\ast f\right) \left( x\right)
\end{equation*}%
\begin{equation*}
=\left( f\ast K_{n}\right) \left( x\right) =\int_{G_{m}}f\left( t\right)
K_{n}\left( x-t\right) d\mu \left( t\right) .
\end{equation*}
where
\begin{eqnarray*}
	K_{n}(x) &:&=\frac{1}{n}\overset{n}{\underset{k=1}{\sum }}D_{k}(x)\text{ \qquad
		\thinspace }\left( \text{ }n\in \mathbb{N}_{+}\text{ }\right)
\end{eqnarray*}

In the literature $ K_{n} $ is called $ n $-th Fej\'er kernel. 

We also define the following maximal operators
\begin{eqnarray*}
	{\sigma }^{\ast}f&=&\sup_{n\in \mathbb{N}}\left\vert \sigma_{{n}}f\right\vert  \\
	\widetilde{\sigma }_{\#}^{\ast}f&=&\sup_{n\in \mathbb{N}}\left\vert \sigma_{2^{n}}f\right\vert.
\end{eqnarray*}.

For any natural number $n\in \mathbb{N}$ we also need the following expression
\begin{equation*}
n=\sum_{i=1}^{s}2^{n_{i}}, \text{ \qquad
	\thinspace } n_{1}<n_{2}<...<n_{s}.
\end{equation*}

Set
\begin{equation*}
n^{\left( i\right) }:=2^{n_{1}}+...+2^{n_{i-1}},\text{ \ }i=2,...,s
\end{equation*}%
and
\begin{equation*}
\mathbb{A}_{0,2}:=\left\{ n\in \mathbb{N}:\text{ }n=2^{0}+2^{2}+%
\sum_{i=3}^{s_{n}}2^{n_{i}}\right\} .
\end{equation*}

Then, for any natural number $n\in \mathbb{N}$ there exists numbers
\begin{equation*}
0 \leq l_{1}\leq m_{1}\leq l_{2}-2<l_{2}\leq m_{2}\leq ...\leq l_{s}-2<l_{s}\leq m_{s}
\end{equation*}
such that it can be written as
\begin{equation*}
n=\sum_{i=1}^{s}\sum_{k=l_{i}}^{m_{i}}2^{k},
\end{equation*}
where $s$ is depending on $n$.

It is evident that
\begin{equation*}
s\leq V\left(n\right) \leq 2s+1.
\end{equation*}

\subsection{Auxiliary lemmas}

\text{ \qquad } The following equality and estimation of Fej\'er kernels with respect to the one-dimensional Walsh-Fourier series is proved in \cite{sws}:

\begin{lemma} \label{lemma3}
	
	Let $n\in \mathbb{N}$ and $n=\sum_{i=1}^{s}2^{n_{i}},$ $n_{1}<n_{2}<...<n_{s}$. Then
	\begin{equation*}
	nK_{n}=\sum_{r=1}^{s}\left( \underset{j=r+1}{\overset{s}{\prod }}%
	w_{2^{n_{j}}}\right) 2^{n_{r}}K_{2^{n_{r}}}+\sum_{t=2}^{s}\left( \underset{%
		j=t+1}{\overset{s}{\prod }}w_{2^{n_{j}}}\right) n^{\left( t\right)
	}D_{2^{n_{t}}},
	\end{equation*}
	and
	\begin{equation*}
	\sup_{n\in \mathbb{N}}\int_{G}\left\vert K_{n}\left( x\right) \right\vert d\mu \left(
	x\right) \leq c<\infty ,
	\end{equation*}%
	where $ c $ is an absolute constant.
\end{lemma}

The following equality is proved in \cite{sws} (see also \cite{gat}):

\begin{lemma} \label{lemma4}
	Let $n>t$ and $t,n\in \mathbb{N}$. Then we have the following expression for $2^{n}$-th Fej\'er kernels with respect to the one-dimensional Walsh-Fourier series:
	\begin{equation*}
	K_{2^{n}}\left( x\right) =\left\{
	\begin{array}{c}
	\text{ }2^{t-1},\text{\  if \ }x\in I_{n}\left( e_{t}\right) , \\
	\frac{2^{n}+1}{2},\text{ \  if \ }x\in I_{n},\text{\ } \\
	0,\text{ \  otherwise.\  }
	\end{array}%
	\right.
	\end{equation*}%
	
\end{lemma}

The following estimation is proved by Goginava \cite{GoSzeged}:
\begin{lemma} \label{lemma5}
	Let  $x\in I_{l+1}\left(e_{k}+e_{l}\right), \ \ k=0,...,M-2, \ \ l=0,...,M-1.$
	Then
	\begin{equation*}
	\int_{I_{M}}\left\vert K_{n}\left( x+t\right) \right\vert d\mu \left(
	t\right) \leq \frac{c2^{l+k}}{n2^{M}},\text{ \  where \ }n>2^{M}.
	\end{equation*}
	
	Let $x\in I_{M}\left(e_{k}\right),m=0,...,M-1.$ Then
	\begin{equation*}
	\int_{I_{M}}\left\vert K_{n}\left( x+t\right) \right\vert d\mu \left(t\right) \leq \frac{c2^{k}}{2^{M}},\text{ \ \ for \ \ }n>2^{M},
	\end{equation*}
	where $ c $ is an absolute constant.
\end{lemma}

The following estimations of Fej\'er kernels with respect to the one-dimensional Walsh-Fourier series is proved in \cite{tep13}:

\begin{lemma} \label{lemma6}
	Let
	$$n=\sum_{i=1}^{r}\sum_{k=l_{i}}^{m_{i}}2^{k},$$
	where
	$$m_{1}\geq l_{1}>l_{1}-2\geq m_{2}\geq l_{2}>l_{2}-2>...>m_{s}\geq l_{s}\geq0.$$
	Then
	\begin{equation*}
	\left\vert nK_{n}\right\vert \leq c\sum_{A=1}^{r}\left( 2^{l_{A}}\left\vert
	K_{2^{l_{A}}}\right\vert +2^{m_{A}}\left\vert K_{2^{m_{A}}}\right\vert
	+2^{l_{A}}\sum_{k=l_{A}}^{m_{A}}D_{2^{k}}\right) +cV\left( n\right),
	\end{equation*}
	where $ c $ is an absolute constant.
\end{lemma}

{\bf Proof}:
Let $$n=\sum_{i=1}^{r}2^{n_{i}}, n_{1}>n_{2}>...>n_{r}\geq 0.$$ By using Lemma \ref{lemma3} for $n$-th Fej\'er kernels we can conclude that
\begin{eqnarray*}
	nK_{n}&=&\sum_{A=1}^{r}\left( \underset{j=1}{\overset{A-1}{\prod }}%
	w_{2^{n_{j}}}\right) \left( \left( 2^{n_{A}}K_{2^{n_{A}}}+\left(
	2^{n_{A}}-1\right) D_{2^{n_{A}}}\right) \right) \\
	&-&\sum_{A=1}^{r}\left( \underset{j=1}{\overset{A-1}{\prod }}%
	w_{2^{n_{j}}}\right) \left( 2^{n_{A}}-1-n^{\left( A\right) }\right)
	D_{2^{n_{A}}}=I_{1}-I_{2}.
\end{eqnarray*}

For $I_{1}$ we have the following equality
\begin{eqnarray*}
	I_{1}&=&\sum_{v=1}^{r}\left( \underset{j=1}{\overset{v-1}{\prod }}\underset{%
		i=l_{j}}{\overset{m_{j}}{\prod }}w_{2^{i}}\right) \left(
	\sum_{k=l_{v}}^{m_{v}}\left( \underset{j=k+1}{\overset{m_{v}}{\prod }}%
	w_{2^{j}}\right) \left( 2^{k}K_{2^{k}}-\left( 2^{k}-1\right)
	D_{2^{k}}\right) \right) \\
	&=&\sum_{v=1}^{r}\left( \underset{j=1}{\overset{v-1}{\prod }}\underset{i=l_{j}}%
	{\overset{m_{j}}{\prod }}w_{2^{i}}\right) \left(
	\sum_{k=0}^{m_{v}}-\sum_{k=0}^{l_{v}-1}\right) \left( \underset{j=k+1}{%
		\overset{m_{v}}{\prod }}w_{2^{j}}\right) \left( 2^{k}K_{2^{k}}-\left(
	2^{k}-1\right) D_{2^{k}}\right) \\
	&=&\sum_{v=1}^{r}\left( \underset{j=1}{\overset{v-1}{\prod }}\underset{i=l_{j}}%
	{\overset{m_{j}}{\prod }}w_{2^{i}}\right) \left( \sum_{k=0}^{m_{v}}\left(
	\underset{j=k+1}{\overset{m_{v}}{\prod }}w_{2^{j}}\right) \left(
	2^{k}K_{2^{k}}-\left( 2^{k}-1\right) D_{2^{k}}\right) \right) \\
	&-&\sum_{v=1}^{r}\left( \underset{j=1}{\overset{v}{\prod }}\underset{i=l_{j}}{%
		\overset{m_{j}}{\prod }}w_{2^{i}}\right) \left( \sum_{k=0}^{l_{v}-1}\left(
	\underset{j=k+1}{\overset{l_{v}-1}{\prod }}w_{2^{j}}\right) \left(
	2^{k}K_{2^{k}}-\left( 2^{k}-1\right) D_{2^{k}}\right) \right) .
\end{eqnarray*}

Since
$$2^{n}-1=\sum_{k=0}^{n-1}2^{k}$$
and
\begin{equation*}
\left( 2^{n}-1\right) K_{2^{n}-1}=\sum_{k=0}^{n-1}\left(
\prod_{j=k+1}^{n-1}w_{2^{j}}\right) \left( 2^{k}K_{2^{k}}-\left(
2^{k}-1\right) D_{2^{k}}\right),
\end{equation*}%
we obtain that
\begin{eqnarray*}
	I_{1}&=&\sum_{v=1}^{r}\left( \underset{j=1}{\overset{v-1}{\prod }}\underset{%
		i=l_{j}}{\overset{m_{j}}{\prod }}w_{2^{i}}\right) \left(
	2^{m_{v}+1}-1\right) K_{2^{m_{v}+1}-1} \\
	&-&\sum_{v=1}^{r}\left( \underset{j=1}{\overset{v}{\prod }}\underset{i=l_{j}}{%
		\overset{m_{j}}{\prod }}w_{2^{i}}\right) \left( 2^{l_{v}}-1\right)
	K_{2^{l_{v}}-1}.
\end{eqnarray*}

If we apply estimations
$$\left\vert K_{2^{n}}\right\vert \leq c\left\vert
K_{2^{n-1}}\right\vert $$
and
$$ \left\vert K_{2^{n}-1}\right\vert \leq c\left\vert
K_{2^{n}}\right\vert +c\  $$
we get that
\begin{equation}
\left\vert I_{1}\right\vert \leq c\sum_{v=1}^{r}\left( 2^{l_{v}}\left\vert
K_{2^{l_{v}}}\right\vert +2^{m_{v}}\left\vert K_{2^{m_{v}}}\right\vert
+cr\right) .  \label{1.12c}
\end{equation}

Let $l_{j}<n_{A}\leq m_{j},$ where $j=1,...,s.$ Then
$$n^{\left( A\right)
}\geq \sum_{v=l_{j}}^{n_{A}-1}2^{v}\geq 2^{n_{A}}-2^{l_{j}}$$ 
and
$$
2^{n_{A}}-1-n^{\left( A\right) }\leq 2^{l_{j}}.$$

If $l_{j}=n_{A}$ where $j=1,...,s,$ then
$$n^{\left( A\right) }\leq 2^{m_{j-1}+1}<2^{l_{j}}.$$

By using these estimations we can conclude that
\begin{equation}
\left\vert I_{2}\right\vert \leq
c\sum_{v=1}^{r}2^{l_{v}}\sum_{k=l_{v}}^{m_{v}}D_{2^{k}}.  \label{1.12d}
\end{equation}

By combining (\ref{1.12c})-(\ref{1.12d}) we get the proof of Lemma \ref{lemma6}.

The following estimations of Fej\'er kernels with respect to the one-dimensional Walsh-Fourier series is proved in \cite{tep13}:

\begin{lemma} \label{lemma7}
	\label{lemma(Tephnadze)} Let $$n=\sum_{i=1}^{s}\sum_{k=l_{i}}^{m_{i}}2^{k},$$
	where
	$$0\leq l_{1}\leq m_{1}\leq l_{2}-2<l_{2}\leq m_{2}\leq ...\leq
	l_{s}-2<l_{s}\leq m_{s}.$$
	
	Then
	\begin{equation*}
	n\left\vert K_{n}\left( x\right) \right\vert \geq \frac{2^{2l_{i}}}{16},%
	\text{ \ \ for \ \ }x\in I_{l_{i}+1}\left( e_{l_{i}-1}+e_{l_{i}}\right).
	\end{equation*}
\end{lemma}

{\bf Proof}: 
If we apply Lemma \ref{lemma3} for $n=\sum_{i=1}^{s}\sum_{k=l_{i}}^{m_{i}}2^{k}$ we can write that
\begin{eqnarray*}
	nK_{n} &=&\sum_{r=1}^{s}\sum_{k=l_{r}}^{m_{r}}\left( \underset{j=r+1}{%
		\overset{s}{\prod }}\underset{q=l_{j}}{\overset{m_{j}}{\prod }}w_{2^{q}}%
	\underset{j=k+1}{\overset{m_{r}}{\prod }}w_{2^{j}}\right) 2^{k}K_{2^{k}} \\
	&+&\sum_{r=1}^{s}\sum_{k=l_{r}}^{m_{r}}\left( \underset{j=r+1}{\overset{s}{%
			\prod }}\underset{q=l_{j}}{\overset{m_{j}}{\prod }}w_{2^{q}}\underset{j=k+1}{%
		\overset{m_{r}}{\prod }}w_{2^{j}}\right) \left(
	\sum_{t=1}^{r-1}\sum_{q=l_{t}}^{m_{t}}2^{q}+\sum_{q=l_{r}}^{k-1}2^{q}\right)
	D_{2^{k}}.
\end{eqnarray*}

Let $x\in I_{l_{i}+1}\left( e_{l_{i}-1}+e_{l_{i}}\right) .$ Then

\begin{eqnarray*}
	n\left\vert K_{n}\right\vert \geq\left\vert
	2^{l_{i}}K_{2^{l_{i}}}\right\vert
	-\sum_{r=1}^{i-1}\sum_{k=l_{r}}^{m_{r}}\left\vert 2^{k}K_{2^{k}}\right\vert
	-\sum_{r=1}^{i-1}\sum_{k=l_{r}}^{m_{r}}\left\vert 2^{k}D_{2^{k}}\right\vert= I-II-III.
\end{eqnarray*}

Lemma \ref{lemma4} follows that

\begin{equation}
I=\left\vert 2^{l_{i}}K_{2^{l_{i}}}\left( x\right) \right\vert =\frac{%
	2^{2l_{i}}}{4}.  \label{2.2.10.0}
\end{equation}

Since $m_{i-1}\leq l_{i}-2$, we easily obtain that the following estimation is true:
\begin{eqnarray} \label{2.2.10.1}
II&\leq &\sum_{n=0}^{l_{i}-2}\left\vert 2^{n}K_{2^{n}}\left( x\right)
\right\vert \\ \notag
&\leq & \sum_{n=0}^{l_{i}-2}2^{n}\frac{\left( 2^{n}+1\right) }{2} \\ \notag
&\leq & \frac{2^{2l_{i}}}{24}+\frac{2^{l_{i}}}{4}-\frac{2}{3}.
\end{eqnarray}

For $III$ we get that
\begin{equation} \label{2.2.10.2}
III\leq \sum_{k=0}^{l_{i}-2}\left\vert 2^{k}D_{2^{k}}\left( x\right)
\right\vert \leq \sum_{k=0}^{l_{i}-2}4^{k}=\frac{2^{2l_{i}}}{12}-\frac{1}{3}.
\end{equation}

By combining (\ref{2.2.10.0}-\ref{2.2.10.2}) we can conclude that
\begin{equation}
n\left\vert K_{n}\left( x\right) \right\vert \geq I-II-III\geq \frac{%
	2^{2l_{i}}}{8}-\frac{2^{l_{i}}}{4}+1.  \label{10.3}
\end{equation}

Suppose that $l_{i}\geq 2$. Then
\begin{equation*}
n\left\vert K_{n}\left( x\right) \right\vert \geq \frac{2^{2l_{i}}}{8}-\frac{%
	2^{2l_{i}}}{16}\geq \frac{2^{2l_{i}}}{16}.
\end{equation*}

If $l_{i}=0$ or $l_{i}=1$, then by applying (\ref{10.3}) we get that
\begin{equation*}
n\left\vert K_{n}\left( x\right) \right\vert \geq \frac{7}{8}\geq \frac{%
	2^{2l_{i}}}{16},
\end{equation*}
Lemma is proved.

The following estimations of Fej\'er kernels with respect to the one-dimensional Walsh-Fourier series is proved in \cite{tep13} (see also \cite{tep_thesis}):

\begin{lemma} \label{lemma3.2.10}
	Let $0<p\leq1$, $2^{k}\leq n<2^{k+1}$ and $\sigma_{n}f$ be Fej\'er means with respect to the one-dimensional Walsh-Fourier series, where $f\in H_{p}(G)$. Then, for any fixed $n\in \mathbb{N}$,
	\begin{eqnarray*}
		&&\left\Vert \sigma_{n}f\right\Vert _{H_{p}(G)} \\
		&&\leq \left\Vert \sup_{0\leq l \leq k}\left\vert \sigma_{2^{l}}f\right\vert \right\Vert _{p}+\left\Vert \sup_{0\leq l\leq k}\left\vert S_{2^{l}}f\right\vert \right\Vert _{p}+\left\Vert \sigma_{n}f\right\Vert _{p} \\
		&&\leq \left\Vert \widetilde{\sigma}_{\#}^{\ast
		}f\right\Vert _{p}+\left\Vert \widetilde{S}_{\#}^{\ast
		}f\right\Vert _{p}+\left\Vert \sigma_{n}f\right\Vert _{p}.
	\end{eqnarray*}
\end{lemma}

{\bf Proof}:
Let consider the following martingale
\begin{eqnarray*}
	f_{\#}&=&\left( S_{2^{k}}\sigma_{n}f,\text{ }k \in \mathbb{N}\right) \\
	&=&\left( \frac{2^{0}\sigma_{2^{0}}}{n}+\frac{(n-2^{0})S_{2^{0}}f}{n},...,\frac{2^{k}\sigma_{2^{k}}f}{n}+\frac{(n-2^{k})S_{2^{k}}f}{n},
	\sigma_{n}f,...,\sigma_{n}f,...\right).
\end{eqnarray*}
By using Lemma \ref{lemma3.2.3} we immediately get
\begin{eqnarray*}
	&&\left\Vert \sigma_{n}f\right\Vert _{H_{p}(G^2)}^{p} \\
	&&\leq \left\Vert \sup_{0\leq l\leq k}\left\vert \sigma_{2^{l}}f\right\vert \right\Vert _{p}^{p}+\left\Vert \sup_{0\leq l\leq k}\left\vert S_{2^{l}}f\right\vert \right\Vert _{p}^{p}+\left\Vert
	S_{n}f\right\Vert _{p}^{p} \\
	&&\leq \left\Vert \widetilde{\sigma}_{\#}^{\ast }f\right\Vert
	_{p}^{p}+\left\Vert \widetilde{S}_{\#}^{\ast }f\right\Vert
	_{p}^{p}+\left\Vert \sigma_{n}f\right\Vert _{p}^{p}.
\end{eqnarray*}

Lemma is proved.

\subsection{Boundedness of subsequences of Fej\'er means with respect to the one-dimensional Walsh-Fourier series on the martingale Hardy spaces}
\text{ \qquad \qquad \qquad \qquad \qquad \qquad \qquad \qquad \qquad \qquad \qquad } In this section we study boundedness of subsequences of Fej\'er means with respect to the one-dimensional Walsh-Fourier series in the martingale Hardy spaces (For details see \cite{tep13}).

First, we consider case $ p=1/2 $. The following estimation is true:
\begin{theorem} \label{th5.1.1}
	\label{theorem1}a) Let $f\in H_{1/2}(G).$ Then there exists an absolute constant $c,$ such that
	\begin{equation*}
	\left\Vert \sigma _{n}f\right\Vert _{H_{1/2}(G)}\leq cV^{2}\left( n\right)
	\left\Vert f\right\Vert _{H_{1/2}(G)}.
	\end{equation*}
	
	b) Let $\left\{ n_{k}:k\in\mathbb{N}_{+}\right\} $  be increasing sequence of natural numbers, such that $\sup_{k \in \mathbb{N}_{+}}V\left(
	n_{k}\right) =\infty $ and $\Phi :\mathbb{N}_{+}\rightarrow \lbrack
	1,\infty )$ be non-decreasing function satisfying the conditions
	$\Phi \left( n\right) \uparrow \infty $
	and
	\begin{equation}
	\overline{\underset{k\rightarrow \infty }{\lim }}\frac{V^{2}\left(
		n_{k}\right) }{\Phi \left( n_{k}\right) }=\infty .  \label{5.1.30}
	\end{equation}%
	Then there exists a martingale $f\in H_{1/2}(G),$ such that
	\begin{equation*}
	\underset{k\in \mathbb{N}}{\sup }\left\Vert \frac{\sigma _{n_{k}}f}{\Phi
		\left(n_{k}\right) }\right\Vert _{1/2}=\infty .
	\end{equation*}
\end{theorem}

{\bf Proof}:
Suppose that
\begin{equation} \label{5.1.12k}
\left\Vert \frac{\sigma _{n}f}{V^{2}\left( n\right) }\right\Vert _{1/2}\leq
c\left\Vert f\right\Vert_{H_{1/2}(G)}.
\end{equation}%

By combining estimations (\ref{1.S2n}), (\ref{sigmamax}) and Lemma \ref{lemma3.2.10} we can conclude that
\begin{eqnarray} \label{5.1.12l}
&&\left\Vert \frac{\sigma _{n}f}{V^{2}\left( n\right) }\right\Vert
_{H_{1/2}(G)}^{1/2} \\ \notag
&\leq &\left\Vert \frac{\sigma _{n}f}{V^{2}\left( n\right) }\right\Vert_{1/2}^{1/2}+\frac{1}{V^{2}\left( n\right) }\left\Vert \sigma_{\#}^{\ast }f \right\Vert_{1/2}^{1/2} +\frac{1}{V^{2}\left( n\right) }\left\Vert\widetilde{S}_{\#}^{\ast }\right\Vert_{1/2}^{1/2}\\ \notag
&\leq & \left\Vert \frac{\widetilde{\sigma} _{n}f}{V^{2}\left( n\right) }\right\Vert_{1/2}^{1/2}+\left\Vert \widetilde{\sigma}_{\#}^{\ast }f \right\Vert_{1/2}^{1/2}+\left\Vert \widetilde{S}_{\#}^{\ast }f \right\Vert_{1/2}^{1/2} \leq  c\left\Vert f\right\Vert _{H_{1/2}(G)}^{1/2}.
\end{eqnarray}

By combining Lemma \ref{lemma3.2.5} and (\ref{5.1.12l}), Theorem \ref{theorem1} will be proved if we show that
\begin{equation*}
\int_{\overline{I_{M}}}\left( \frac{\left\vert \sigma _{n}a\right\vert }{V^{2}\left( n\right) }\right) ^{1/2}d\mu  \leq c<\infty ,
\end{equation*}%
for any $ 1/2 $-atom $a$.

Without loss the generality we may assume that $a$ is $1/2$-atom, with support $I,$ for which $\mu \left( I\right) =2^{-M},$  $I=I_{M}.$ It is easy to check that $\sigma _{n}\left( a\right) =0,$ when $n\leq 2^{M}.$ Therefore, we may assume that $n>2^{M}.$ Set
\begin{eqnarray*}
	&& II_{\alpha _{A}}^{1}\left( x\right):=2^{M}\int_{I_{M}}2^{\alpha
		_{A}}\left\vert K_{2^{\alpha _{A}}}\left( x+t\right) \right\vert d\mu \left(t\right) ,\text{\ } \\
	&& II_{l_{A}}^{2}\left( x\right)
	=2^{M}\int_{I_{M}}2^{l_{A}}\sum_{k=l_{A}}^{m_{A}}D_{2^{k}}\left(
	x+t\right) d\mu \left( t\right).
\end{eqnarray*}

Let $x\in I_{M}.$ Since $\sigma _{n}$ is bounded from $L_{\infty }(G)$ to $L_{\infty }(G)$, for $n>2^{M}$ and $\left\Vert a\right\Vert _{\infty }\leq 2^{2M},$ by using Lemma \ref{lemma5} we can conclude that
\begin{eqnarray*}
	&&\frac{\left\vert \sigma _{n}a\left( x\right) \right\vert }{V^{2}\left(n\right)} \\
	&& \leq \frac{c}{V^{2}\left( n\right) }\int_{I_{M}}\left\vert a\left(
	x\right) \right\vert \left\vert K_{n}\left( x+t\right) \right\vert d\mu\left( t\right) \\
	&& \leq \frac{c\left\Vert a\right\Vert _{\infty }}{V^{2}\left(n\right) }\int_{I_{M}}\left\vert K_{n}\left( x+t\right) \right\vert d\mu \left(t\right) \\
	&& \leq \frac{c2^{2M}}{V^{2}\left( n\right) }\int_{I_{M}}\left\vert
	K_{n}\left( x+t\right) \right\vert d\mu \left( t\right)\\
	&&\leq \frac{c2^{M}}{V^{2}\left( n\right) }\left\{
	\sum_{A=1}^{s}\int_{I_{M}}2^{l_{A}}\left\vert K_{2^{l_{A}}}\left( x+t\right)
	\right\vert d\mu \left( t\right) +\int_{I_{M}}2^{m_{A}}\left\vert
	K_{2^{m_{A}}}\left( x+t\right) \right\vert d\mu \left( t\right) \right\}\\
	&&+\frac{c2^{M}}{V^{2}\left( n\right) }\sum_{A=1}^{s}\int_{I_{M}}2^{l_{A}}%
	\sum_{k=l_{A}}^{m_{A}}D_{2^{k}}\left( x+t\right) d\mu \left( t\right) +\frac{c2^{M}}{V^{2}\left( n\right) }\int_{I_{M}}V\left( n\right) d\mu \left(t\right)\\
	&&=\frac{c}{V^{2}\left( n\right) }\sum_{A=1}^{s}\left( II_{^{l_{A}}}^{1}\left(
	x\right) +II_{^{m_{A}}}^{1}\left( x\right) +II_{l_{A}}^{2}\left( x\right)
	\right) +c.
\end{eqnarray*}

Hence,
\begin{eqnarray*} \\
	&&\int_{\overline{I_{M}}}\left\vert \frac{\sigma _{n}a\left( x\right) }{%
		V^{2}\left( n\right) }\right\vert ^{1/2}d\mu \left( x\right)
	\\
	&\leq &\frac{c}{V\left( n\right) }\left( \sum_{A=1}^{s}\int_{\overline{I_{M}}%
	}\left\vert II_{l_{A}}^{1}\left( x\right) \right\vert ^{1/2}d\mu \left(
	x\right) \right. \\
	&&\left. +\int_{\overline{I_{M}}}\left\vert II_{m_{A}}^{1}\left( x\right)
	\right\vert ^{1/2}d\mu \left( x\right) +\int_{\overline{I_{M}}}\left\vert
	II_{l_{A}}^{2}\left( x\right) \right\vert ^{1/2}d\mu \left( x\right) \right)
	+c.
\end{eqnarray*}

Since  $s\leq 4V\left( n\right) $ we obtain that Theorem \ref{th5.1.1} will be proved if we show that
\begin{equation} \label{5.1.11.1}
\int_{\overline{I_{M}}}\left\vert II_{\alpha _{A}}^{1}\left( x\right)
\right\vert ^{1/2}d\mu \left( x\right) \leq c<\infty ,\text{\ }\int_{%
	\overline{I_{M}}}\left\vert II_{l_{A}}^{2}\left( x\right) \right\vert
^{1/2}d\mu \left( x\right) \leq c<\infty ,\text{\ }
\end{equation}%
where $\alpha _{A}=l_{A}$ or $\alpha _{A}=m_{A},$ $A=1,...,s$.

Let $t\in I_{M}$ and $x\in I_{l+1}\left( e_{k}+e_{l}\right) ,$ $0\leq
k<l<\alpha _{A}\leq M$ or $0\leq k<l\leq M\leq \alpha _{A}.$ Since $x+t\in
I_{l+1}\left( e_{k}+e_{l}\right),$ by applying Lemma \ref{lemma4} we can conclude that
\begin{equation} \label{5.1.10a}
K_{2^{\alpha _{A}}}\left( x+t\right) =0\text{ \ \ and \ \ }II_{\alpha
	_{A}}^{1}\left( x\right) =0.
\end{equation}

Let $x\in I_{l+1}\left( e_{k}+e_{l}\right),$ $0\leq k<\alpha _{A}\leq l\leq M.$ Then $x+t\in I_{l+1}\left( e_{k}+e_{l}\right),$ where $t\in I_{M}$ and if we apply again Lemma \ref{lemma4} we get that
\begin{equation} \label{5.1.10b}
2^{\alpha _{A}}\left\vert K_{2^{\alpha _{A}}}\left( x+t\right) \right\vert\leq 2^{\alpha _{A}+k}\text{ \ \ and \ \ }II_{\alpha _{A}}^{1}\left(x\right) \leq 2^{\alpha_{A}+k}.
\end{equation}

Analogously to (\ref{5.1.10b}) for $0\leq \alpha _{A}\leq k<l\leq M$ we can prove that
\begin{equation} \label{5.1.10c}
2^{\alpha_{A}}\left\vert K_{2^{\alpha _{A}}}\left( x+t\right) \right\vert\leq 2^{2\alpha_{A}},\text{ \ }II_{\alpha _{A}}^{1}\left( x\right)\leq2^{2\alpha_{A}},\text{\ }t\in I_{M},\text{\ }x\in I_{l+1}\left(
e_{k}+e_{l}\right).
\end{equation}

Let $0\leq \alpha _{A}\leq M-1,$ where $A=1,...,s.$ According to (\ref{2.1.2}) and (\ref{5.1.10a}-\ref{5.1.10c}) we find that
\begin{eqnarray*}
	&&\int_{\overline{I_{M}}}\left\vert II_{\alpha _{A}}^{1}\left( x\right)
	\right\vert ^{1/2}d\mu \left( x\right) \\
	&=&\overset{M-2}{\underset{k=0}{\sum }}\overset{M-1}{\underset{l=k+1}{\sum }}%
	\int_{I_{l+1}\left( e_{k}+e_{l}\right) }\left\vert II_{\alpha
		_{A}}^{1}\left( x\right) \right\vert ^{1/2}d\mu \left( x\right)\\
	&+&\overset{M-1%
	}{\underset{k=0}{\sum }}\int_{I_{M}\left( e_{k}\right) }\left\vert
	II_{\alpha _{A}}^{1}\left( x\right) \right\vert ^{1/2}d\mu \left( x\right) \\
	&\leq & c\overset{\alpha _{A}-1}{\underset{k=0}{\sum }}\overset{M-1}{\underset{%
			l=\alpha _{A}+1}{\sum }}\int_{I_{l+1}\left( e_{k}+e_{l}\right) }2^{\left(
		\alpha _{A}+k\right) /2}d\mu \left( x\right) \\
	& +& c\overset{M-2}{\underset{%
			k=\alpha _{A}}{\sum }}\overset{M-1}{\underset{l=k+1}{\sum }}%
	\int_{I_{l+1}\left( e_{k}+e_{l}\right) }2^{\alpha _{A}}d\mu \left( x\right) \\
	&+& c\overset{\alpha _{A}-1}{\underset{k=0}{\sum }}\int_{I_{M}\left(
		e_{k}\right) }2^{\left( \alpha _{A}+k\right) /2}d\mu \left( x\right)
	+c\overset{M-1}{\underset{k=\alpha _{A}}{\sum }}\int_{I_{M}\left( e_{k}\right)
	}2^{\alpha _{A}}d\mu \left( x\right) \\
	&\leq & c\overset{\alpha _{A}-1}{\underset{k=0}{\sum }}\overset{M-1}{\underset{%
			l=\alpha _{A}+1}{\sum }}\frac{2^{\left( \alpha _{A}+k\right) /2}}{2^{l}}+c%
	\overset{M-2}{\underset{k=\alpha _{A}}{\sum }}\overset{M-1}{\underset{l=k+1}{%
			\sum }}\frac{2^{\alpha _{A}}}{2^{l}} \\
	&+& c\overset{\alpha _{A}-1}{\underset{k=0}{%
			\sum }}\frac{2^{\left( \alpha _{A}+k\right) /2}}{2^{M}}+c\overset{M-1}{%
		\underset{k=\alpha _{A}}{\sum }}\frac{2^{\alpha _{A}}}{2^{M}}
	\leq c<\infty.
\end{eqnarray*}

Let $\alpha _{A}\geq M.$ Analogously to $II_{\alpha _{A}}^{1}\left( x\right)$ we can prove (\ref{5.1.11.1}), for $A=1,...,s.$

Now, prove boundedness of $II_{l_{A}}^{2}$. Let $t\in I_{M}$ and $x\in I_{i}\backslash I_{i+1},$ $i\leq l_{A}-1.$ Since $x+t\in I_{i}\backslash
I_{i+1},$ if we apply first equality of Lemma \ref{lemma1} we get that
\begin{equation} \label{5.1.13a}
II_{l_{A}}^{2}\left( x\right) =0.
\end{equation}

Let $x\in I_{i}\backslash I_{i+1},$ $l_{A}\leq i\leq m_{A}.$ Since  $n\geq
2^{M}$ and $t\in I_{M},$ if we apply first equality of Lemma \ref{lemma1} we get that
\begin{equation} \label{5.1.13b}
II_{l_{A}}^{2}\left( x\right) \leq
2^{M}\int_{I_{M}}2^{l_{A}}\sum_{k=l_{A}}^{i}D_{2^{k}}\left( x+t\right) d\mu\left( t\right) \leq c2^{l_{A}+i}.
\end{equation}

Let $x\in I_{i}\backslash I_{i+1},$ $m_{A}<i\leq M-1.$ Then  $x+t\in I_{i}\backslash I_{i+1},$ for any $t\in I_{M}$ and by first equality of Lemma \ref{lemma1} we have that
\begin{equation} \label{5.1.13c}
II_{l_{A}}^{2}\left( x\right) \leq c2^{M}\int_{I_{M}}2^{l_{A}+m_{A}}\leq
c2^{l_{A}+m_{A}}.
\end{equation}

Let  $0\leq l_{A}\leq m_{A}\leq M.$ Then, in the view of (\ref{2.1.2})  and (\ref{5.1.13a}-\ref{5.1.13c}) we can conclude that
\begin{eqnarray*}
	&&\int_{\overline{I_{M}}}\left\vert II_{l_{A}}^{2}\left( x\right) \right\vert
	^{1/2}d\mu \left( x\right) \\
	&=&\left(
	\sum_{i=0}^{l_{A}-1}+\sum_{i=l_{A}}^{m_{A}}+\sum_{i=m_{A}+1}^{M-1}\right)
	\int_{I_{i}\backslash I_{i+1}}\left\vert II_{l_{A}}^{2}\left( x\right)
	\right\vert ^{1/2}d\mu \left( x\right) \\
	&\leq& c\sum_{i=l_{A}}^{m_{A}}\int_{I_{i}\backslash I_{i+1}}2^{\left(
		l_{A}+i\right) /2}d\mu \left( x\right) \\
	&+& c\sum_{i=m_{A}+1}^{M-1}\int_{I_{i}\backslash I_{i+1}}2^{\left(
		l_{A}+m_{A}\right) /2}d\mu \left( x\right) \\
	&\leq & c\sum_{i=l_{A}}^{m_{A}}2^{\left( l_{A}+i\right) /2}\frac{1}{2^{i}} \\
	&+& c\sum_{i=m_{A}+1}^{M-1}2^{\left(l_{A}+m_{A}\right) /2}\frac{1}{2^{i}}\leq c<\infty .
\end{eqnarray*}

Analogously, we can prove same estimations in the cases
$0\leq l_{A}\leq M<m_{A}$ and $M\leq l_{A}\leq m_{A}.$

Now, we prove part b) of Theorem \ref{th5.1.1}. According to (\ref{5.1.30}), there exists increasing sequence $\left\{ \alpha _{k}:\text{ }k\in\mathbb{N}_{+}\right\} \subset \left\{ n_{k}\text{ }:k\in\mathbb{N}_{+}\right\} $ of natural numbers such that
\begin{equation} \label{5.1.2aaa}
\sum_{k=1}^{\infty }\frac{\Phi ^{1/4}\left( \alpha _{k}\right)}{V^{1/2}\left(
	\alpha _{k}\right) }\leq c<\infty.
\end{equation}

Let $f=\left(f_{n},n\in\mathbb{N}_{+}\right)$ be martingale form Example \ref{example2.2.1}, where
\begin{equation*}
\lambda _{k}:=\Phi ^{1/2}\left( \alpha _{k}\right) /V\left( \alpha
_{k}\right).
\end{equation*}

According to (\ref{5.1.2aaa}) we get that condition (\ref{3.3.2aa}) is fulfilled and it follows that $f=\left(f_{n},n\in\mathbb{N}_{+}\right)$.

By applying (\ref{3.3.10AA}) we get that
\begin{equation} \label{5.1.5aa}
\widehat{f}(j)
\end{equation}
\begin{equation*}
=\left\{
\begin{array}{ll}
2^{\left\vert \alpha _{k}\right\vert }\Phi ^{1/2}\left( \alpha _{k}\right)/V\left( \alpha _{k}\right) , & \text{\thinspace \thinspace }j\in \left\{2^{\left\vert \alpha _{k}\right\vert },...,2^{_{\left\vert \alpha_{k}\right\vert +1}}-1\right\} ,\text{ }k \in\mathbb{N}_{+} \\
0\,, & \text{\thinspace }j\notin \bigcup\limits_{k=0}^{\infty }\left\{2^{_{\left\vert \alpha _{k}\right\vert }},..., 2^{_{\left\vert \alpha_{k}\right\vert +1}}-1\right\}.
\end{array}
\right.
\end{equation*}

Let $2^{\left\vert \alpha _{k}\right\vert }<j<\alpha _{k}.$ If we apply (\ref{3.3.11AA}) we get that
\begin{equation} \label{5.1.sn}
S_{j}f=S_{2^{\left\vert \alpha_{k}\right\vert }}f+\frac{w_{2^{\left\vert \alpha_{k}\right\vert }}D_{j-2^{\left\vert \alpha _{k}\right\vert }}\Phi^{1/2}\left( \alpha _{k}\right) }{V\left( \alpha _{k}\right)}
\end{equation}

Hence,
\begin{eqnarray} \label{5.1.7aaa}
&&\frac{\sigma _{_{\alpha _{k}}}f}{\Phi \left( \alpha _{k}\right) } \\ \notag
&=&\frac{1}{\Phi \left( \alpha _{k}\right) \alpha_{k}}\sum_{j=1}^{2^{\left\vert \alpha_{k}\right\vert }}S_{j}f \\ \notag
&+&\frac{1}{\Phi \left( \alpha _{k}\right) \alpha _{k}}\sum_{j=2^{\left\vert \alpha _{k}\right\vert }+1}^{\alpha_{k}}S_{j}f
\\ \notag
&=&\frac{\sigma _{_{2^{\left\vert \alpha _{k}\right\vert }}}f}{\Phi \left(\alpha _{k}\right) \alpha _{k}} \\ \notag
&+&\frac{\left( \alpha _{k}-2^{\left\vert\alpha_{k}\right\vert }\right) S_{2^{\left\vert \alpha _{k}\right\vert }}f}{\Phi \left(\alpha _{k}\right) \alpha _{k}} \\ \notag
&+&\frac{w_{2^{\left\vert \alpha_{k}\right\vert }}2^{\left\vert \alpha_{k}\right\vert }\Phi^{1/2}\left(\alpha _{k}\right)}{\Phi \left( \alpha_{k}\right) V\left( \alpha_{k}\right) \alpha _{k}}\sum_{j=2^{_{\left\vert\alpha_{k}\right\vert }}+1}^{\alpha _{k}}D_{j-2^{\left\vert \alpha_{k}\right\vert }} \\ \notag
&=& III_{1}+III_{2}+III_{3}.
\end{eqnarray}%

For $ III_{3} $ we can conclude that
\begin{eqnarray} \label{5.1.9aaa}
&&\left\vert III_{3}\right\vert \\ \notag
&=&\frac{2^{\left\vert \alpha _{k}\right\vert }\Phi ^{1/2}\left( \alpha
	_{k}\right) }{\Phi \left( \alpha _{k}\right) V\left( \alpha _{k}\right)
	\alpha _{k}}\left\vert \sum_{j=1}^{\alpha _{k}-2^{_{\left\vert \alpha
			_{k}\right\vert }}}D_{j}\right\vert \\ \notag
&=&\frac{2^{\left\vert \alpha
		_{k}\right\vert }\Phi ^{1/2}\left( \alpha _{k}\right) }{\Phi \left( \alpha
	_{k}\right) V\left( \alpha _{k}\right) \alpha _{k}}\left( \alpha
_{k}-2^{_{\left\vert \alpha _{k}\right\vert }}\right) \left\vert K_{\alpha
	_{k}-2^{_{\left\vert \alpha _{k}\right\vert }}}\right\vert \\ \notag
&\geq& \frac{c\left( \alpha _{k}-2^{_{\left\vert \alpha _{k}\right\vert }}\right)\left\vert K_{\alpha _{k}-2^{_{\left\vert \alpha _{k}\right\vert
	}}}\right\vert }{\Phi ^{1/2}\left( \alpha _{k}\right) V\left( \alpha
	_{k}\right).}
\end{eqnarray}

Let
\begin{equation*}
\alpha_{k}=\sum_{i=1}^{r_{k}}\sum_{k=l_{i}^{k}}^{m_{i}^{k}}2^{k},
\end{equation*}
where
\begin{equation*}
m_{1}^{k}\geq l_{1}^{k}>l_{1}^{k}-2\geq m_{2}^{k}\geq
l_{2}^{k}>l_{2}^{k}-2\geq ...\geq m_{s}^{k}\geq l_{s}^{k}\geq 0.
\end{equation*}

Since (see theorems \ref{th4.1.1} and \ref{th5.1.1})
\begin{equation*}
\left\Vert III_{1}\right\Vert
_{1/2}\leq c, \left\Vert III_{2}\right\Vert _{1/2}\leq c,
\end{equation*}
and
\begin{equation*}
\mu \left\{E_{l_{i}^{k}}\right\} \geq 1/2^{l_{i}^{k}-1},
\end{equation*}
By combining (\ref{5.1.7aaa}), (\ref{5.1.9aaa}) and Lemma \ref{lemma7} we get that
\begin{eqnarray*}
	&&\int_{G}\left\vert \sigma _{\alpha _{k}}f(x)/\Phi \left( \alpha _{k}\right)
	\right\vert ^{1/2}d\mu \left( x\right) \\
	&\geq & \left\Vert III_{3}\right\Vert
	_{1/2}^{1/2}-\left\Vert III_{2}\right\Vert _{1/2}^{1/2}-\left\Vert
	III_{1}\right\Vert _{1/2}^{1/2} \\
	&\geq & c\text{ }\underset{i=2}{\overset{r_{k}-2}{\sum }}\int_{E_{l_{i}^{k}}}%
	\left\vert 2^{2l_{i}^{k}}/\left( \Phi ^{1/2}\left( \alpha _{k}\right)
	V\left( \alpha _{k}\right) \right) \right\vert ^{1/2}d\mu \left( x\right) -2c \\
	&\geq & c\overset{r_{k}-2}{\underset{i=2}{\sum }}1/\left( V^{1/2}\left( \alpha
	_{k}\right) \Phi ^{1/4}\left( \alpha _{k}\right) \right)-2c \\
	&\geq&
	cr_{k}/\left( V^{1/2}\left( \alpha _{k}\right) \Phi ^{1/4}\left( \alpha
	_{k}\right) \right) \\
	&\geq & cV^{1/2}\left( \alpha _{k}\right) /\Phi ^{1/4}\left( \alpha _{k}\right)
	\rightarrow \infty ,\text{ as }k\rightarrow \infty .
\end{eqnarray*}

Theorem \ref{th5.1.1} is proved.

\begin{theorem} \label{th5.1.2}
	a) Let $0<p<1/2,$ $f\in H_{p}(G).$ Then there exists an absolute constant $c_{p}$ depending only on $p$ such that
	\begin{equation*}
	\text{ }\left\Vert \sigma _{n}f\right\Vert _{H_{p}(G)}\leq c_{p}2^{d\left(
		n\right) \left( 1/p-2\right) }\left\Vert f\right\Vert _{H_{p}(G)}.
	\end{equation*}
	
	b) Let $0<p<1/2$ and $\Phi \left( n\right):\mathbb{N}_{+}\rightarrow \lbrack
	1,\infty ) $ be non-decreasing function such that
	\begin{equation}  \label{5.1.31aaa}
	\sup_{k\in \mathbb{N}_+}d\left( n_{k}\right) =\infty ,\text{ \ \ }\overline{\underset{k\rightarrow \infty }{\lim }}\frac{2^{d\left( n_{k}\right) \left(1/p-2\right) }}{\Phi \left( n_{k}\right) }=\infty .
	\end{equation}%
	Then there exist a martingale $f\in H_{p}(G),$ such that
	\begin{equation*}
	\underset{k\in\mathbb{N}_+}{\sup}\left\Vert\frac{\sigma_{n_{k}}f}{\Phi \left(n_{k}\right)}\right\Vert_{weak-L_{p}(G)}=\infty.
	\end{equation*}
\end{theorem}
{\bf Proof}:
Let $n\in \mathbb{N}.$ Analogously to (\ref{5.1.12l}) it is sufficient to prove that
\begin{equation*}
\int\limits_{\overline{I}_{M}}\left( 2^{d\left( n\right) \left( 2-1/p\right)
}\left\vert \sigma _{n}\left( a\right) \right\vert \right) ^{p}d\mu \leq
c_{p}<\infty ,
\end{equation*}%
for every $ p $-atom $a$, where $I$ denotes support of the atom.

Analogously to Theorem \ref{th5.1.1} we may assume that $a$
is $ p $-atom with support $\ I=I_{M}$, $\mu \left( I_{M}\right)=2^{-M}$ and  $n>2^{M}.$
Since $\left\Vert a\right\Vert _{\infty }\leq
2^{M/p} $ we can conclude that
\begin{eqnarray*}
	&& 2^{d\left( n\right) \left( 2-1/p\right) }\left\vert \sigma _{n}a\right\vert \\
	&\leq & 2^{d\left( n\right) \left( 2-1/p\right)}\left\Vert a\right\Vert _{\infty }\int_{I_{M}}\left\vert K_{n}\left(
	x+t\right) \right\vert d\mu \left( t\right) \\
	&\leq & 2^{d\left( n\right) \left( 2-1/p\right) }2^{M/p}\int_{I_{M}}\left\vert
	K_{n}\left( x+t\right) \right\vert d\mu \left( t\right).
\end{eqnarray*}

Let $x\in I_{l+1}\left( e_{k}+e_{l}\right) ,\,0\leq k,l\leq \left[ n\right]
\leq M.$ Then, by applying Lemma \ref{lemma4} we get that $K_{n}\left( x+t\right)=0, $ where $t\in I_{M}$ and hence,
\begin{equation} \label{5.1.12a}
2^{d\left( n\right) \left( 2-1/p\right) }\left\vert \sigma _{n}a \right\vert =0.
\end{equation}

Let $x\in I_{l+1}\left( e_{k}+e_{l}\right) ,\,\left[ n\right] \leq k,l\leq M$
or  $k\leq \left[ n\right] \leq l\leq M.$ Then Lemma \ref{lemma6} follows that
\begin{eqnarray} \label{5.1.12}
2^{d\left( n\right) \left( 2-1/p\right) }\left\vert \sigma _{n}a \right\vert &\leq &  2^{d\left( n\right) \left( 2-1/p\right)}2^{M\left( 1/p-2\right) +k+l}  \\ \notag
&\leq & c_{p}2^{\left[ n\right] \left( 1/p-2\right) +k+l}.
\end{eqnarray}

By combining (\ref{2.1.2}), (\ref{5.1.12a}) and (\ref{5.1.12}) we can conclude that
\begin{eqnarray*}
	&&\int_{\overline{I_{M}}}\left\vert 2^{d\left( n\right) \left( 2-1/p\right)
	}\sigma _{n}a\left( x\right) \right\vert ^{p}d\mu \left( x\right) \\
	&\leq& \left( \overset{\left[ n\right] -2}{\underset{k=0}{\sum }}\overset{%
		\left[ n\right] -1}{\underset{l=k+1}{\sum }}+\overset{\left[ n\right] -1}{%
		\underset{k=0}{\sum }}\overset{M-1}{\underset{l=\left[ n\right] }{\sum }}+%
	\overset{M-2}{\underset{k=\left[ n\right] }{\sum }}\overset{M-1}{\underset{%
			l=k+1}{\sum }}\right) \int_{I_{l+1}\left( e_{k}+e_{l}\right) }\left\vert
	2^{d\left( n\right) \left( 2-1/p\right) }\sigma _{n}a\left( x\right)
	\right\vert ^{p}d\mu \left( x\right) \\
	&+&\overset{M-1}{\underset{k=0}{\sum }}\int_{I_{M}\left( e_{k}\right)
	}\left\vert 2^{d\left( n\right) \left( 2-1/p\right) }\sigma _{n}a\left(
	x\right) \right\vert ^{p}d\mu \left( x\right) \\
	&\leq & c_{p}\overset{M-2}{%
		\underset{k=\left[ n\right] }{\sum }}\overset{M-1}{\underset{l=k+1}{\sum }}%
	\frac{1}{2^{l}}2^{\left[ n\right] \left( 2p-1\right) }2^{p\left( k+l\right) } \\
	&+& c_{p}\overset{\left[ n\right] }{\underset{k=0}{\sum }}\overset{M-1}{%
		\underset{l=\left[ n\right] +1}{\sum }}\frac{1}{2^{l}}2^{\left[ n\right]
		\left( 2p-1\right) }2^{p\left( k+l\right) } \\
	&+&\frac{c_{p}2^{\left[ n\right]
			\left( 2p-1\right) }}{2^{M}}\overset{\left[ n\right] }{\underset{k=0}{\sum }}%
	2^{p\left( k+M\right) }<c_{p}<\infty .
\end{eqnarray*}

Now, we prove part b) of Theorem \ref{th5.1.2}. According to   (\ref{5.1.31aaa}) there exists an increasing sequence of natural numbers $\left\{ \alpha _{k}:\text{ }%
k\in \mathbb{N}_+\right\} \subset \left\{n_{k}:\text{ }k\in \mathbb{N}_+\right\},$ such that
$\alpha _{0}\geq 3$ and
\begin{equation} \label{5.1.121}
\sum_{\eta =0}^{\infty }u^{-p}\left( \alpha _{\eta }\right) <c_{p}<\infty, \text{ \ \ }
\text{\ }u\left( \alpha _{k}\right) =2^{d\left( \alpha _{k}\right) \left(
	1/p-2\right) /2}/\Phi ^{1/2}\left( \alpha _{k}\right) .
\end{equation}

Let $f$ be martingale from Example \ref{example2.2.1}, where
\begin{equation*}
\lambda_{k}=u^{-1}\left( \alpha _{k}\right),
\end{equation*}

If we apply  (\ref{5.1.121}) we get that (\ref{3.3.2aa}) is fulfilled and it follows that $f\in H_{p}(G).$ According to (\ref{3.3.10AA}) we have that
\begin{equation} \label{5.1.6aa}
\widehat{f}(j)
=\left\{
\begin{array}{ll}
2^{\left\vert \alpha _{k}\right\vert \left( 1/p-1\right) }/u\left( \alpha
_{k}\right) , & j\in \left\{ 2^{\left\vert \alpha _{k}\right\vert
},...,2^{_{\left\vert \alpha _{k}\right\vert +1}}-1\right\} ,\text{ }%
k\in\mathbb{N}_+, \\
0\,, & j\notin \bigcup\limits_{k=0}^{\infty }\left\{ 2^{\left\vert \alpha
	_{k}\right\vert },...,2^{_{\left\vert \alpha _{k}\right\vert +1}}-1\right\} .
\end{array}
\right.
\end{equation}

Let $2^{\left\vert \alpha _{k}\right\vert }<j<\alpha _{k}.$ Then, analogously to (\ref{5.1.sn}) and (\ref{5.1.7aaa}), if we apply (\ref{5.1.6aa}) we get that
\begin{eqnarray*}
	&&\frac{\sigma _{_{\alpha _{k}}}f}{\Phi \left( \alpha _{k}\right)} \\ &=&\frac{%
		\sigma _{_{2^{\left\vert \alpha _{k}\right\vert }}}f}{\Phi \left( \alpha
		_{k}\right) \alpha _{k}}+\frac{\left( \alpha _{k}-2^{\left\vert \alpha
			_{k}\right\vert }\right) S_{2^{\left\vert \alpha _{k}\right\vert }}f}{\Phi
		\left( \alpha _{k}\right) \alpha _{k}} \\
	&+&\frac{2^{\left\vert \alpha _{k}\right\vert \left( 1/p-1\right) }}{\Phi
		\left( \alpha _{k}\right) u\left( \alpha _{k}\right) \alpha _{k}}%
	\sum_{j=2^{_{\left\vert \alpha _{k}\right\vert }}}^{\alpha _{k}-1}\left(
	D_{_{j}}-D_{2^{\left\vert \alpha _{k}\right\vert }}\right) \\
	&=& IV_{1}+IV_{2}+IV_{3}.
\end{eqnarray*}

Let $\alpha _{k}\in \mathbb{N}$ and $E_{\left[ \alpha _{k}\right] }:=I_{_{%
		\left[ \alpha _{k}\right] +1}}\left( e_{\left[ \alpha _{k}\right] -1}+e_{%
	\left[ \alpha _{k}\right] }\right) .$ Since $\left[ \alpha
_{k}-2^{\left\vert \alpha _{k}\right\vert }\right] =\left[ \alpha _{k}\right],$
analogously to (\ref{5.1.9aaa}), if we apply Lemma \ref{lemma7} for $IV_{3}$ we have the following estimation
\begin{eqnarray*}
	\left\vert IV_{3}\right\vert &=&\frac{2^{\left\vert \alpha _{k}\right\vert
			\left( 1/p-1\right) }}{\Phi \left( \alpha _{k}\right) u\left( \alpha
		_{k}\right) \alpha _{k}}\left( \alpha _{k}-2^{\left\vert \alpha
		_{k}\right\vert }\right) \left\vert K_{\alpha _{k}-2^{\left\vert \alpha
			_{k}\right\vert }}\right\vert \\
	&=&\frac{2^{\left\vert \alpha _{k}\right\vert \left( 1/p-1\right) }}{\Phi
		\left( \alpha _{k}\right) u\left( \alpha _{k}\right) \alpha _{k}}\left\vert
	2^{\left[ \alpha _{k}\right] }K_{\left[ \alpha _{k}\right] }\right\vert \\
	&\geq &
	\frac{2^{\left\vert \alpha _{k}\right\vert \left( 1/p-2\right) }2^{2\left[
			\alpha _{k}\right] -4}}{\Phi \left( \alpha _{k}\right) u\left( \alpha
		_{k}\right) } \\
	&\geq& 2^{\left\vert \alpha _{k}\right\vert \left( 1/p-2\right) /2}2^{2\left[
		\alpha _{k}\right] -4}/\Phi ^{1/2}\left( \alpha _{k}\right) .
\end{eqnarray*}
Hence,
\begin{eqnarray*}
	&&\left\Vert IV_{3}\right\Vert _{weak-L_{p}(G)}^{p} \\
	&\geq & \left( \frac{2^{\left\vert \alpha _{k}\right\vert \left( 1/p-2\right)
			/2}2^{2\left[ \alpha _{k}\right] -4}}{\Phi ^{1/2}\left( \alpha _{k}\right) }%
	\right) ^{p}\mu \left\{ x\in G:\text{ }\left\vert IV_{3}\right\vert \geq
	\frac{2^{\left\vert \alpha _{k}\right\vert \left( 1/p-2\right) /2}2^{2\left[
			\alpha _{k}\right] -4}}{\Phi ^{1/2}\left( \alpha _{k}\right) }\right\} \\
	&\geq & c_{p}\left( 2^{2\left[ \alpha _{k}\right] +\left\vert \alpha
		_{k}\right\vert \left( 1/p-2\right) /2}/\Phi ^{1/2}\left( \alpha _{k}\right)
	\right) ^{p}\mu (E_{\left[ \alpha _{k}\right] })\\
	&\geq & c_{p}\left( 2^{\left( \left\vert \alpha _{k}\right\vert -\left[ \alpha
		_{k}\right] \right) \left( 1/p-2\right) }/\Phi \left( \alpha _{k}\right)
	\right) ^{p/2} \\
	&=& c_{p}\left( 2^{d\left( \alpha _{k}\right) \left( 1/p-2\right)
	}/\Phi \left( \alpha _{k}\right) \right) ^{p/2}\rightarrow \infty ,\text{\ as }k\rightarrow \infty .
\end{eqnarray*}

By combining Corollary \ref{cor4.1.2} and first part of Theorem \ref{th5.1.2} we find that
\begin{equation*}
\left\Vert IV_{1}\right\Vert _{weak-L_{p}(G)}\leq
c_{p}<\infty , \text{ \ \ } \left\Vert IV_{2}\right\Vert _{weak-L_{p}(G)}\leq
c_{p}<\infty.
\end{equation*}

On the other hand, for sufficiently large $ n $ we can conclude that
\begin{eqnarray*}
	&&\left\Vert \sigma _{\alpha _{k}}f\right\Vert _{weak-L_{p}(G)}^{p} \\
	&\geq&\left\Vert IV_{3}\right\Vert _{weak-L_{p}(G)}^{p}-\left\Vert
	IV_{2}\right\Vert _{weak-L_{p}(G)}^{p}-\left\Vert IV_{1}\right\Vert
	_{weak-L_{p}(G)}^{p} \\
	&\geq & \frac{1}{2}\left\Vert IV_{3}\right\Vert _{weak-L_{p}(G)}^{p}\rightarrow \infty ,\text{ as }k\rightarrow \infty .
\end{eqnarray*}

Theorem \ref{th5.1.2} is proved.

The proofs of Corollaries \ref{cor5.1.1}-\ref{cor5.1.3} are similar to the proofs of Corollaries \ref{cor4.1.2}-\ref{cor4.1.4}. So, we leave out the details of proofs and just present these results:

\begin{corollary}\label{cor5.1.1}
	Let $p>0$ and $f\in H_{p}(G)$. Then
	\begin{equation*}
	\left\Vert \sigma_{2^{k}}f-f\right\Vert _{H_p(G)}\rightarrow 0,\text{ as }k\rightarrow\infty.
	\end{equation*}
\end{corollary}

\begin{corollary}\label{cor5.1.2}
	Let $p>0$ and $f\in H_{p}(G)$. Then
	\begin{equation*}
	\left\Vert \sigma_{2^{k}+2^{k-1}}f-f\right\Vert _{H_p(G)}\rightarrow 0,\text{ as }k\rightarrow\infty.
	\end{equation*}
\end{corollary}

\begin{corollary}\label{cor5.1.3}
	Let $0<p<1/2$. Then there exists a martingale $f\in H_{p}(G)$, such that
	\begin{equation*}
	\left\Vert \sigma_{2^{k}+1}f-f\right\Vert _{weak-L_{p}(G)}\nrightarrow 0,\text{ as }k\rightarrow\infty.
	\end{equation*}

	On the other hand, for any $f\in H_{1/2}(G)$ the following is true:
	\begin{equation*}
	\left\Vert \sigma_{2^{k}+1}f-f\right\Vert _{H_{1/2}(G)}\rightarrow 0,\text{ as }k\rightarrow\infty.
	\end{equation*}
	
\end{corollary}

\subsection{Modulus of continuity and convergence in norm of subsequences of  Fej\'er means with respect to the one-dimensional Walsh-Fourier series on the martingale Hardy spaces}

\text{ \qquad \qquad \qquad \qquad \qquad \qquad \qquad \qquad \qquad \qquad \qquad \qquad \qquad \qquad \qquad  } In this section we apply Theorem \ref{th5.1.1} and Theorem \ref{th5.1.2} to find necessary and sufficient conditions for modulus of continuity of martingale $ f\in H_p $, for which subsequences of  Fej\'er means with respect to the one-dimensional Walsh-Fourier series converge in $ H_p $-norm.

First, we prove the following result:

\begin{theorem} \label{theorem5.2.1}
	a) Let $f\in H_{1/2}(G), \ \ \sup_{k\in \mathbb{N}_+}V\left( n_{k}\right) =\infty$ and
	\begin{equation} \label{5.2.cond}
	\omega _{H_{p}(G)}\left( 1/2^{\left\vert n_{k}\right\vert },f\right) =o\left(
	1/V^{2}\left( n_{k}\right) \right) ,\text{ \ as \ }k\rightarrow \infty .
	\end{equation}%
	Then
	\begin{equation*}
	\left\Vert \sigma_{n_{k}}f-f\right\Vert _{H_{1/2}(G)}\rightarrow 0,\text{ as }k\rightarrow\infty.
	\end{equation*}
	
	b) Let $\sup_{k\in \mathbb{N}_+}V\left( n_{k}\right) =\infty .$ Then there exists a martingale $f\in H_{1/2}(G),$ such that
	\begin{equation} \label{5.2.cond2}
	\omega _{H_{1/2}(G)}\left( 1/2^{\left\vert n_{k}\right\vert },f\right) =O\left(
	1/V^{2}\left( n_{k}\right) \right) ,\text{ \ as \ }k\rightarrow \infty
	\end{equation}%
	and
	\begin{equation} \label{5.2.kn3}
	\left\Vert \sigma _{n_{k}}f-f\right\Vert _{H_{1/2}(G)}\nrightarrow 0, \text{ \ as \ }k\rightarrow \infty .
	\end{equation}
\end{theorem}

{\bf Proof}:
Let $f\in H_{1/2}(G)$ and $2^{k}<n\leq 2^{k+1}.$ Then
\begin{eqnarray*} 
	&&\left\Vert \sigma _{n}f-f\right\Vert _{H_{1/2}(G)}^{1/2} \\
	&\leq & \left\Vert \sigma _{n}f-\sigma _{n}S_{2^{k}}f\right\Vert
	_{H_{1/2}(G)}^{1/2}
	+\left\Vert \sigma _{n}S_{2^{k}}f-S_{2^{k}}f\right\Vert
	_{H_{1/2}(G)}^{1/2}
	+\left\Vert S_{2^{k}}f-f\right\Vert _{H_{1/2}(G)}^{1/2}\\
	&=&\left\Vert \sigma _{n}\left( S_{2^{k}}f-f\right) \right\Vert
	_{H_{1/2}(G)}^{1/2}
	+\left\Vert S_{2^{k}}f-f\right\Vert _{H_{1/2}(G)}^{1/2}+\left\Vert
	\sigma _{n}S_{2^{k}}f-S_{2^{k}}f\right\Vert _{H_{1/2}(G)}^{1/2}\\
	&\leq& c\left( V\left( n\right) +1\right) \omega
	_{H_{1/2}(G)}^{1/2}\left( 1/2^{k},f\right)+\left\Vert \sigma
	_{n}S_{2^{k}}f-S_{2^{k}}f\right\Vert _{H_{1/2}(G)}^{1/2}.\\
\end{eqnarray*}

It is evident that
\begin{eqnarray*}
	\sigma _{n}S_{2^{k}}f-S_{2^{k}}f
	=\frac{2^{k}}{n}\left( S_{2^{k}}\sigma
	_{2^{k}}f-S_{2^{k}}f\right)
	=\frac{2^{k}}{n}S_{2^{k}}\left( \sigma
	_{2^{k}}f-f\right).
\end{eqnarray*}

Let $p>0.$ By combining Corollaries \ref{cor4.1.2} and \ref{cor5.1.1} we can conclude that
\begin{eqnarray*}
	\left\Vert \sigma _{n}S_{2^{k}}f-S_{2^{k}}f\right\Vert _{H_{1/2}(G)}^{1/2}
	&\leq &\frac{2^{k/2}}{n^{1/2}}\left\Vert S_{2^{k}}\left( \sigma
	_{2^{k}}f-f\right) \right\Vert _{H_{1/2}(G)}^{1/2}
	\\
	&&
	\leq \left\Vert \sigma
	_{2^{k}}f-f\right\Vert _{H_{1/2}(G)}^{1/2}\rightarrow 0,\text{\ as \ } k\rightarrow
	\infty.
\end{eqnarray*}

Now, we prove part b) of Theorem \ref{theorem5.2.1}. Since $\sup_{k\in \mathbb{N}_+}V(\alpha _{k})=\infty,$ then there exists a martingale  $\{\alpha _{k}:k\in \mathbb{N}_+\}\subset \{n_{k}:k\in \mathbb{N}_+\}$ such that $V(\alpha _{k})\uparrow \infty,$ as $k\rightarrow \infty $ and
\begin{equation} \label{5.2.4.7}
V^{2}(\alpha _{k})\leq V(\alpha _{k+1}).
\end{equation}

Let $f$ be martingale from Example \ref{example2.2.1}, where
\begin{equation*}
\lambda_{k}=V^{-2}(\alpha _{k}),
\end{equation*}

If we apply (\ref{5.2.4.7}) we get that condition (\ref{3.3.2aa}) is fulfilled and it follows that $f\in H_{p}(G).$
By using (\ref{3.3.10AA}) we find that

\begin{equation} \label{5.2.4.13}
\widehat{f}(j)=\left\{
\begin{array}{ll}
2^{_{\left\vert \alpha _{k}\right\vert }}/V^{2}(\alpha _{k}), & \text{%
	\thinspace }j\in \left\{ 2^{\left\vert \alpha _{k}\right\vert
},...,2^{_{\left\vert \alpha _{k}\right\vert +1}}-1\right\} ,\text{ }%
k\in\mathbb{N}_+ \\
0\,, & j\notin \bigcup\limits_{k=0}^{\infty }\left\{ 2^{_{\left\vert \alpha
		_{k}\right\vert }},...,2^{_{\left\vert \alpha _{k}\right\vert +1}}-1\right\}
.\text{ }%
\end{array}%
\right.
\end{equation}

By combining (\ref{3.3.2aa0}) and (\ref{5.2.4.7}) we can conclude that
\begin{eqnarray} \label{4.2.4.12}
&& w_{H_{1/2}(G)}(1/2^n,f)
=\left\Vert f-S_{2^{n}}f\right\Vert_{H_{1/2}(G)} \\ \notag
&&\leq
\sum\limits_{i=n+1}^{\infty }1/V^{2}(\alpha _{i})= O\left( 1/V^{2}(\alpha
_{n})\right),\text{ \ as \ } n\rightarrow \infty.
\end{eqnarray}

Let $2^{_{\left\vert \alpha _{k}\right\vert }}<j\leq \alpha _{k}.$ By using (\ref{3.3.11AA}) we get that
\begin{equation*}
S_{j}f=S_{2^{_{\left\vert \alpha _{k}\right\vert }}}f+\frac{%
	2^{_{\left\vert \alpha _{k}\right\vert }}w_{2^{_{\left\vert \alpha
				_{k}\right\vert }}}D_{j-2^{_{\left\vert \alpha _{k}\right\vert }}}}{%
	V^{2}(\alpha _{k})}.
\end{equation*}%
Hence,
\begin{eqnarray}  \label{5.2.nn}
&&\sigma _{\alpha _{k}}f-f \\ \notag
&=&\frac{2^{_{\left\vert \alpha _{k}\right\vert }}}{%
	\alpha _{k}}\left( \sigma _{2^{_{\left\vert \alpha _{k}\right\vert
}}}f-f\right)+\frac{\alpha _{k}-2^{_{\left\vert \alpha _{k}\right\vert }}}{\alpha _{k}}%
\left(S_{2^{_{\left\vert \alpha _{k}\right\vert }}}f-f\right) \\ \notag
&+&\frac{%
	2^{_{\left\vert \alpha _{k}\right\vert }}w_{2^{_{\left\vert \alpha
				_{k}\right\vert }}}\left( \alpha _{k}-2^{_{\left\vert \alpha _{k}\right\vert
	}}\right) K_{\alpha _{k}-2^{_{\left\vert \alpha _{k}\right\vert }}}}{\alpha_{k}V^{2}(\alpha _{k})}.
\end{eqnarray}

According to (\ref{1.S2n000}),  (\ref{fe22222}) and (\ref{5.2.nn}) we have that
\begin{eqnarray} \label{5.2.a11}
&&\Vert \sigma _{\alpha _{k}}f-f\Vert _{1/2}^{1/2} \\ \notag
&\geq & \frac{c}{V(\alpha _{k})}%
\Vert \left( \alpha _{k}-2^{_{\left\vert \alpha _{k}\right\vert }}\right)
K_{\alpha _{k}-2^{_{\left\vert \alpha _{k}\right\vert }}}\Vert _{1/2}^{1/2} \\ \notag
&-&\left( \frac{2^{_{\left\vert \alpha _{k}\right\vert }}}{\alpha _{k}}\right)
^{1/2}\Vert \sigma _{2^{_{\left\vert \alpha _{k}\right\vert }}}f-f\Vert
_{1/2}^{1/2} \\ \notag
&-&\left( \frac{\alpha _{k}-2^{_{\left\vert \alpha _{k}\right\vert
}}}{\alpha _{k}}\right) ^{1/2}\Vert S_{2^{_{\left\vert \alpha
			_{k}\right\vert }}}f-f\Vert _{1/2}^{1/2}.
\end{eqnarray}

Let
$$\alpha _{k}=\sum_{i=1}^{r_{k}}\sum_{k=l_{i}^{k}}^{m_{i}^{k}}2^{k},$$
where
$$m_{1}^{k}\geq l_{1}^{k}>l_{1}^{k}-2\geq m_{2}^{k}\geq
l_{2}^{k}>l_{2}^{k}-2>...>m_{s}^{k}\geq l_{s}^{k}\geq 0$$
and
$$E_{l_{i}^{k}}:=I_{_{l_{i}^{k}+1}}\left(e_{l_{i}^{k}-1}+e_{l_{i}^{k}}\right).$$

By using Lemma \ref{lemma7} we get that
\begin{eqnarray} \label{5.2.33}
&&\int_{G}\left\vert \left( \alpha _{k}-2^{_{\left\vert \alpha _{k}\right\vert
}}\right) K_{\alpha _{k}-2^{_{\left\vert \alpha _{k}\right\vert }}}\left(
x\right) \right\vert ^{1/2}d\mu   \\ \notag
&\geq &\frac{1}{16}\underset{i=2}{\overset{r_{k}-2}{\sum }}%
\int_{E_{l_{i}^{k}}}\left\vert \left( \alpha _{k}-2^{_{\left\vert \alpha
		_{k}\right\vert }}\right) K_{\alpha _{k}-2^{_{\left\vert \alpha
			_{k}\right\vert }}}\left( x\right) \right\vert ^{1/2}d\mu \left( x\right)
\\ \notag
&\geq & \frac{1}{16}\underset{i=2}{\overset{r_{k}-2}{\sum }}\frac{1}{%
	2^{l_{i}^{k}}}2^{l_{i}^{k}}
\geq cr_{k}\geq cV(\alpha _{k}).
\end{eqnarray}

By combining estimations (\ref{5.2.a11}-\ref{5.2.33}), Corollaries \ref{cor4.1.2} and \ref{cor5.1.1} we get that (\ref{5.2.kn3}) holds true and Theorem \ref{theorem5.2.1} is proved.

\begin{theorem}\label{theorem5.2.2}
	a) Let $0<p<1/2,$ $f\in H_{p}(G)$, $\ \sup_{k\in \mathbb{N}_+}d\left(
	n_{k}\right) =\infty $ and
	\begin{equation}
	\omega _{H_{p}(G)}\left( 1/2^{\left\vert n_{k}\right\vert },f\right) =o\left(
	1/2^{d\left( n_{k}\right) \left( 1/p-2\right) }\right) ,\text{ as \ }%
	k\rightarrow \infty .  \label{5.2.cond3}
	\end{equation}%
	Then
	\begin{equation} \label{5.2.fe2}
	\left\Vert\sigma_{n_{k}}f-f\right\Vert_{H_{p}(G)}\rightarrow0,\text{\ as \ } k\rightarrow \infty.
	\end{equation}
	
	b) Let $\sup_{k\in \mathbb{N}_+}d\left( n_{k}\right) =\infty .$ Then there exists a martingale $f\in H_{p}(G)$ $\left( 0<p<1/2\right) ,$\ \ such that
	\begin{eqnarray} \label{cond4}
	&&\omega _{H_{p}(G)}\left( 1/2^{\left\vert n_{k}\right\vert },f\right)= O\left(
	1/2^{d\left( n_{k}\right) \left( 1/p-2\right) }\right) ,\text{ \ as \ }%
	k\rightarrow \infty
	\end{eqnarray}%
	and
	\begin{equation}
	\left\Vert\sigma_{n_{k}}f-f\right\Vert_{weak-L_{p}(G)}\nrightarrow
	0,\text{\ as \ }k\rightarrow \infty.
	\label{kn4}
	\end{equation}
\end{theorem}

{\bf Proof}:
\textbf{\ }Let $0<p<1/2.$ Then under condition (\ref{5.2.cond3}) if we repeat steps of the proof of Theorem \ref{theorem5.2.1}, we immediately get that (\ref{5.2.fe2}) holds.

Let prove part b) of Theorem \ref{theorem5.2.2}. Since  $\sup_{k}d\left(
n_{k}\right) =\infty ,$ there exists $\{\alpha _{k}:k\in \mathbb{N}_+\}\subset\{n_{k}:k\in \mathbb{N}_+\}$ such that $\sup_{k\in \mathbb{N}_+}d\left( \alpha _{k}\right) =\infty $
and
\begin{equation} \label{5.2.4.18}
2^{2d\left( \alpha _{k}\right) \left( 1/p-2\right) }\leq 2^{d\left( \alpha
	_{k+1}\right) \left( 1/p-2\right) }.
\end{equation}

Let $f$ be a martingale from Lemma \ref{example2.2.1}, where
\begin{equation*}
\lambda_{k}=2^{-\left( 1/p-2\right) d\left( \alpha _{i}\right)}.
\end{equation*}

If we use (\ref{5.2.4.18}) we conclude that condition (\ref{3.3.2aa}) is fulfilled and it follows that $f\in H_{p}(G).$

According to (\ref{3.3.10AA}) we get that

\begin{equation} \label{5.2.4.22}
\widehat{f}(j)=\left\{
\begin{array}{ll}
2^{\left( 1/p-2\right) \left[ \alpha _{k}\right] }, & \text{\thinspace
	\thinspace }j\in \left\{ 2^{\left\vert \alpha _{k}\right\vert
},...,2^{_{\left\vert \alpha _{k}\right\vert +1}}-1\right\} ,\text{ }%
k\in \mathbb{N}_+ \\
0\,, & \text{\thinspace }j\notin \bigcup\limits_{n=0}^{\infty }\left\{
2^{_{\left\vert \alpha _{n}\right\vert }},...,2^{_{\left\vert \alpha
		_{n}\right\vert +1}}-1\right\} .\text{ }%
\end{array}%
\right.
\end{equation}

By combining (\ref{3.3.2aa0}) and (\ref{5.2.4.18}) we have that
\begin{eqnarray} \label{4.21}
\omega _{H_{p}(G)}(1/2^{\left\vert \alpha _{k}\right\vert },f) 
&\leq & \sum\limits_{i=k}^{\infty }1/2^{d\left( \alpha _{i}\right) \left(
	1/p-2\right) } \\ \notag
&=& O\left( 1/2^{d\left( \alpha _{k}\right) \left( 1/p-2\right)
}\right)\text{ \ as \ } k \rightarrow \infty.
\end{eqnarray}

Analogously to the proof of previous theorem, if we use also Corollaries \ref{cor4.1.2} and \ref{cor5.1.1}, for the sufficiently large $ k $ we can conclude that
\begin{eqnarray} \label{5.2.nnmm00}
&&\Vert \sigma _{\alpha _{k}}f-f\Vert _{weak-L_{p}(G)}^{p} \\ \notag
&\geq& 2^{\left(
	1-2p\right) \left[ \alpha _{k}\right] }\Vert \left( \alpha
_{k}-2^{\left\vert \alpha _{k}\right\vert }\right) K_{\alpha
	_{k}-2^{_{\left\vert \alpha _{k}\right\vert }}}\Vert _{weak-L_{p}(G)}^{p} \\  \notag
&-&\left( \frac{2^{\left\vert \alpha _{k}\right\vert }}{\alpha _{k}}\right)
^{p}\Vert \sigma _{2^{_{\left\vert \alpha _{k}\right\vert }}}f-f\Vert
_{weak-L_{p}(G)}^{p} \\ \notag
&-&\left( \frac{\alpha _{k}-2^{_{\left\vert \alpha
			_{k}\right\vert }}}{\alpha _{k}}\right) ^{p}\Vert S_{2^{_{\left\vert \alpha
			_{k}\right\vert }}}f-f\Vert _{weak-L_{p}(G)}^{p} \\ \notag
&\geq & 2^{\left(
	1-2p\right) \left[ \alpha _{k}\right]-1 }\Vert \left( \alpha
_{k}-2^{\left\vert \alpha _{k}\right\vert }\right) K_{\alpha
	_{k}-2^{_{\left\vert \alpha _{k}\right\vert }}}\Vert _{weak-L_{p}(G)}^{p}
\end{eqnarray}

Let $x\in E_{\left[ \alpha _{k}\right] }.$ Lemma \ref{lemma7} follows that
\begin{eqnarray*}
	\mu \left( x\in G:\left( \alpha _{k}-2^{_{\left\vert \alpha _{k}\right\vert
	}}\right) \left\vert K_{\alpha _{k}-2^{_{\left\vert \alpha _{k}\right\vert
	}}}\right\vert \geq 2^{2\left[ \alpha _{k}\right] -4}\right)\geq \mu \left(
	E_{\left[ \alpha _{k}\right] }\right) \geq 1/2^{\left[ \alpha _{k}\right]
		-4}
\end{eqnarray*}
and
\begin{eqnarray} \label{5.2.nnmm}
2^{2p\left[ \alpha _{k}\right] -4}\mu \left( x\in G:\left( \alpha
_{k}-2^{_{\left\vert \alpha _{k}\right\vert }}\right) \left\vert K_{\alpha
	_{k}-2^{_{\left\vert \alpha _{k}\right\vert }}}\right\vert \geq 2^{2\left[
	\alpha _{k}\right] -4}\right) \geq 2^{\left( 2p-1\right) \left[ \alpha _{k}%
	\right] -4}.
\end{eqnarray}

Hence, by combining (\ref{1.S2n000}),  (\ref{fe22222}), (\ref{5.2.nnmm00}) and (\ref{5.2.nnmm}) we get that
\begin{equation*}
\left\Vert \sigma _{n_{k}}f-f\right\Vert _{weak-L_p(G)}\nrightarrow
0,\,\,\,\text{as}\,\,\,k\rightarrow \infty.
\end{equation*}%
The proof of Theorem \ref{theorem5.2.2} is complete.

By using Theorem \ref{theorem5.2.2} we easily get an important result which was proved in \cite{tep6}:

\begin{corollary}\label{corollary5.2.2ss2}
	a) Let  $f\in H_{1/2}(G)$ and
	\begin{equation*}
	\omega _{H_{1/2}(G)}\left( 1/2^{k},f\right) =o\left(\frac{1}{k^2}\right) ,\text{ as \ }%
	k\rightarrow \infty .
	\end{equation*}%
	Then
	\begin{equation*}
	\left\Vert \sigma _{k}f-f\right\Vert _{H_{1/2}(G)}\rightarrow
	0,\,\,\,\text{ \ as \ }k\rightarrow \infty.
	\end{equation*}
	
	b) There exists a martingale $f\in H_{1/2}(G),$ for which
	\begin{equation*}
	\omega _{H_{1/2}(G)}\left( 1/2^{k},f\right) =O\left(\frac{1}{k^2}\right),\text{ \ as \ }
	k\rightarrow \infty
	\end{equation*}%
	and
	\begin{equation*}
	\left\Vert \sigma_{k}f-f\right\Vert_{1/2}\nrightarrow0,\,\,\,\text{ \ as \ }k\rightarrow\infty .
	\end{equation*}
\end{corollary}

\begin{corollary}\label{corollary5.2.2ss}
	a) Let $0<p<1/2,$ $f\in H_{p}(G)$ and
	
	\begin{equation*}
	\omega _{H_{p}(G)}\left( 1/2^{k},f\right) =o\left(
	1/2^{k( 1/p-2)}\right) ,\text{ as \ }%
	k\rightarrow \infty .
	\end{equation*}%
	Then
	\begin{equation*}
	\left\Vert \sigma _{k}f-f\right\Vert _{H_{p}(G)}\rightarrow
	0,\,\,\,\text{as\thinspace \thinspace \thinspace }k\rightarrow \infty .
	\end{equation*}
	
	b) Then there exists a martingale $f\in H_{p}(G)$ $\left( 0<p<1/2\right) ,$\ \ for which
	
	\begin{equation*}
	\omega _{H_{p}(G)}\left( 1/2^{k},f\right) =O\left(
	1/2^{k\left( 1/p-2\right) }\right) ,\text{ \ as \ }%
	k\rightarrow \infty
	\end{equation*}%
	and
	\begin{equation*}
	\left\Vert \sigma _{k}f-f\right\Vert _{weak-L_{p}(G)}\nrightarrow
	0,\,\,\,\text{ \ as \ }k\rightarrow \infty .
	\end{equation*}
\end{corollary}

\subsection{Strong convergence of  Fej\'er means with respect to the one-dimensional Walsh-Fourier series on the martingale Hardy spaces}

\text{ \qquad \qquad \qquad \qquad \qquad \qquad \qquad \qquad \qquad \qquad \qquad \qquad \qquad \qquad \qquad  } In this section we consider strong convergence results of Fej\'er means with respect to the one-dimensional Walsh-Fourier series in the martingale Hardy spaces, when $ 0<p\leq 1/2 $ (for details see \cite{tep5}).

The following is true:

\begin{theorem} \label{theorem5.3.1}
	a) Let $0<p\leq 1/2$ and $f\in H_{p}(G)$. Then there exists a constant $c_{p},$ depending only on $p$, such that
	
	\begin{equation*}
	\frac{1}{\log ^{\left[ 1/2+p\right] }n}\overset{n}{\underset{m=1}{\sum }}%
	\frac{\left\Vert \sigma _{m}f\right\Vert _{H_{p}(G)}^{p}}{m^{2-2p}}\leq
	c_{p}\left\Vert f\right\Vert _{H_{p}(G)}^{p}.
	\end{equation*}%
	b) Let $0<p<1/2,$ $\Phi :\mathbb{N}_{+}\rightarrow\lbrack 1,\infty )$ be non-decreasing function, such that $\Phi \left( n\right) \uparrow \infty $ and
	\begin{equation*}
	\overline{\underset{k\rightarrow \infty }{\lim }}\frac{k^{
			2-2p }}{\Phi \left({k}\right) }=\infty .
	\end{equation*}
	
	Then there exists a martingale $f\in H_{p}(G),$ such that
	\begin{equation*}
	\underset{m=1}{\overset{\infty }{\sum }}\frac{\left\Vert \sigma
		_{m}f\right\Vert _{weak-L_{p}(G)}^{p}}{\Phi \left( m\right) }=\infty .
	\end{equation*}
\end{theorem}

{\bf Proof}:
Suppose that
\begin{equation*}
\frac{1}{\log ^{\left[ 1/2+p\right] }n}\overset{n}{\underset{m=1}{\sum }}%
\frac{\left\Vert \sigma _{m}f\right\Vert _{p}^{p}}{m^{2-2p}}\leq
c_{p}\left\Vert f\right\Vert _{H_{p}(G)}^{p}.
\end{equation*}

By combining (\ref{1.S2n}), (\ref{sigmamax}) and Lemma \ref{lemma3.2.10} we can conclude that
\begin{eqnarray} \label{5.3.5.3}
&&\frac{1}{\log ^{\left[ 1/2+p\right] }n}\overset{n}{\underset{m=1}{\sum }}%
\frac{\left\Vert \sigma _{m}f\right\Vert _{H_{p}(G)}^{p}}{m^{2-2p}} \\
&& \leq\frac{1}{\log ^{\left[ 1/2+p\right] }n}\overset{n}{\underset{m=1}{\sum }}%
\frac{\left\Vert \sigma _{m}f\right\Vert _{p}^{p}}{m^{2-2p}}+\|\widetilde{\sigma}_{\#}^{\ast }f\|_{H_p(G)}+\|\widetilde{S}_{\#}^{\ast }f\|_{H_p(G)}  \notag \\
&&  \leq c_{p}\left\Vert f\right\Vert _{H_{p}(G)}^{p}. \notag
\end{eqnarray}

According to Lemma \ref{lemma3.2.5} and (\ref{5.3.5.3}) Theorem \ref{theorem5.3.1} will be proved if we show that

\begin{equation*}
\frac{1}{\log ^{\left[ 1/2+p\right] }n}\overset{n}{\underset{m=1}{\sum }}%
\frac{\left\Vert \sigma _{m}a\right\Vert _{p}^{p}}{m^{2-2p}}\leq c<\infty ,%
\text{ \ \ \ }m=2,3,...
\end{equation*}%
for any  $ p $-atom $a$. We may assume that $a$ is $ p $-atom, with support  $ I$, $\mu \left( I\right) =2^{-M}$ and $I=I_{M}.$ It is evident that $\sigma _{n}\left( a\right) =0,$ when $n\leq 2^{M}.$ Therefore, we may assume that $n>2^{M}.$

Let $x\in I_{M}.$ Since $\sigma _{n}$ is bounded from $L_{\infty }(G)$ to $L_{\infty }(G)$ (The boundedness follows fact that Fej\'er kernels are uniformly bounded in the space $L_{1}(G)$, which is proved in Lemma \ref{lemma3}) and $\left\Vert
a\right\Vert _{\infty }\leq 2^{M/p}$ we can conclude that
\begin{equation*}
\int_{I_{M}}\left\vert \sigma _{m}a\left( x\right) \right\vert ^{p}d\mu
\left( x\right) \leq \left\Vert \sigma _{m}a\right\Vert _{\infty
}^{p}/2^{M}
\end{equation*}
\begin{equation*}
\leq \left\Vert a\right\Vert _{\infty
}^{p}/2^{M}\leq c<\infty ,\text{ }0<p\leq 1/2.
\end{equation*}

Let $0<p\leq 1/2.$ Then
\begin{equation*}
\frac{1}{\log ^{\left[ 1/2+p\right] }n}\overset{n}{\underset{m=1}{\sum }}%
\frac{\int_{I_{M}}\left\vert \sigma _{m}a\left( x\right) \right\vert
	^{p}d\mu \left( x\right) }{m^{2-2p}}
\end{equation*}
\begin{equation*}
\leq \frac{c}{\log ^{\left[ 1/2+p\right]
	}n}\overset{n}{\underset{m=1}{\sum }}\frac{1}{m^{2-2p}}\leq c<\infty .
\end{equation*}

It is evident that
\begin{equation*}
\left\vert \sigma _{m}a\left( x\right) \right\vert \leq
\int_{I_{M}}\left\vert a\left( t\right) \right\vert \left\vert K_{m}\left(
x+t\right) \right\vert d\mu \left( t\right)
\end{equation*}%
\begin{equation*}
\leq
2^{M/p}\int_{I_{M}}\left\vert K_{m}\left( x+t\right) \right\vert d\mu \left(
t\right) .
\end{equation*}

Lemma \ref{lemma4} follows that
\begin{equation} \label{5.3.12}
\left\vert \sigma _{m}a\left( x\right) \right\vert \leq \frac{%
	c2^{k+l}2^{M\left( 1/p-1\right) }}{m},\text{ \ \ }x\in I_{l+1}\left(
e_{k}+e_{l}\right) ,\,0\leq k<l<M
\end{equation}%
and
\begin{equation} \label{5.3.12a}
\left\vert \sigma _{m}a\left( x\right) \right\vert \leq c2^{M\left(
	1/p-1\right) }2^{k},\text{ \ }x\in I_{M}\left( e_{k}\right) ,\,0\leq k<M.
\end{equation}

If we use identity (\ref{2.1.2}) and (\ref{5.3.12}-\ref{5.3.12a}) we get that
\begin{eqnarray} \label{7.2}
&&\int_{\overline{I_{M}}}\left\vert \sigma _{m}a\left( x\right) \right\vert
^{p}d\mu \left( x\right)  \\
&=&\overset{M-2}{\underset{k=0}{\sum }}\overset{M-1}{\underset{l=k+1}{\sum }}%
\int_{I_{l+1}\left( e_{k}+e_{l}\right) }\left\vert \sigma _{m}a\left(
x\right) \right\vert ^{p}d\mu \left( x\right)
\notag \\
& &
+\overset{M-1}{\underset{k=0}{%
		\sum }}\int_{I_{M}\left( e_{k}\right) }\left\vert \sigma _{m}a\left(
x\right) \right\vert ^{p}d\mu \left( x\right)  \notag \\
&\leq &c\overset{M-2}{\underset{k=0}{\sum }}\overset{M-1}{\underset{l=k+1}{%
		\sum }}\frac{1}{2^{l}}\frac{2^{p\left( k+l\right) }2^{M\left( 1-p\right) }}{%
	m^{p}}+c\overset{M-1}{\underset{k=0}{\sum }}\frac{1}{2^{M}}2^{M\left(
	1-p\right) }2^{pk}  \notag \\
&\leq &\frac{c2^{M\left( 1-p\right) }}{m^{p}}\overset{M-2}{\underset{k=0}{%
		\sum }}\overset{M-1}{\underset{l=k+1}{\sum }}\frac{2^{p\left( k+l\right) }}{%
	2^{l}}+c\overset{M-1}{\underset{k=0}{\sum }}\frac{2^{pk}}{2^{pM}}  \notag \\
&\leq &\frac{c2^{M\left( 1-p\right) }M^{\left[ 1/2+p\right] }}{m^{p}}+c.
\notag
\end{eqnarray}

Hence,
\begin{eqnarray*}
	&&\frac{1}{\log ^{\left[ 1/2+p\right] }n}\overset{n}{\underset{m=2^{M}+1}{%
			\sum }}\frac{\int_{\overline{I_{M}}}\left\vert \sigma _{m}a\left( x\right)
		\right\vert ^{p}d\mu \left( x\right) }{m^{2-2p}} \\
	&\leq &\frac{1}{\log ^{\left[ 1/2+p\right] }n}\left( \overset{n}{\underset{%
			m=2^{M}+1}{\sum }}\frac{c2^{M\left( 1-p\right) }M^{\left[ 1/2+p\right] }}{%
		m^{2-p}}+\overset{n}{\underset{m=2^{M}+1}{\sum }}\frac{c}{m^{2-2p}}\right)
	<c<\infty .
\end{eqnarray*}

The proof of part a) of theorem \ref{theorem5.3.1} is complete.

\bigskip Now, we prove part b) of Theorem \ref{theorem5.3.1}.  Let $\Phi \left( n\right) $ non-decreasing function satisfying the condition
\begin{equation} \label{5.3.12j}
\underset{k\rightarrow \infty }{\lim }\frac{2^{\left( \left\vert
		n_{k}\right\vert +1\right) \left( 2-2p\right) }}{\Phi \left( 2^{\left\vert
		n_{k}\right\vert +1}\right) }=\infty .
\end{equation}

According to (\ref{5.3.12j}), there exists an increasing sequence  $\left\{ \alpha _{k}:%
\text{ }k\in\mathbb{N}_+\right\} \subset \left\{ n_{k}:\text{ }k\in\mathbb{N}_+\right\} $ such that
\begin{equation} \label{5.3.122}
\left\vert \alpha _{k}\right\vert \geq 2,\text{ \ \ \ where \ }k\in\mathbb{N}_+
\end{equation}%
and
\begin{eqnarray} \label{5.3.121}
\sum_{\eta =0}^{\infty }\frac{\Phi ^{1/2}\left( 2^{\left\vert \alpha _{\eta
		}\right\vert +1}\right) }{2^{\left\vert \alpha _{\eta }\right\vert \left(
		1-p\right) }}= 2^{1-p}\sum_{\eta =0}^{\infty }\frac{\Phi ^{1/2}\left(
	2^{\left\vert \alpha _{\eta }\right\vert +1}\right) }{2^{\left( \left\vert
		\alpha _{\eta }\right\vert +1\right) \left( 1-p\right) }}<c<\infty .
\end{eqnarray}

Let $f=\left(f_{n},\text{ }n\in\mathbb{N}_+\right) \in H_{p}(G)$ be a martingale from the Example \ref{example2.2.1}, where
\begin{equation*}
\lambda _{k}=\frac{\Phi ^{1/2p}\left( 2^{\left\vert \alpha _{k}\right\vert
		+1}\right) }{2^{\left( \left\vert \alpha _{k}\right\vert \right) \left(
		1/p-1\right) }}
\end{equation*}

By combining (\ref{3.3.2aa}) and (\ref{5.3.121}) we get that $f\in H_{p}(G).$
According to (\ref{3.3.10AA}) we have that
\begin{eqnarray} \label{5.3.6aa}
&&\widehat{f}(j)=\left\{
\begin{array}{l}
\Phi ^{1/2p}\left(2^{\left\vert \alpha _{k}\right\vert +1}\right) 
\text{ \ if \ }j\in \left\{ 2^{\left\vert \alpha
	_{k}\right\vert },...,2^{\left\vert \alpha _{k}\right\vert +1}-1\right\} ,%
\text{ }k\in\mathbb{N}_+, \\
0\text{ },\text{  \ if \ } j\notin
\bigcup\limits_{k=0}^{\infty }\left\{ 2^{\left\vert \alpha _{k}\right\vert
},...,2^{\left\vert \alpha _{k}\right\vert +1}-1\right\} .\text{ }
\end{array}
\right. 
\end{eqnarray}

Let $2^{\left\vert \alpha _{k}\right\vert }<n<2^{\left\vert \alpha
	_{k}\right\vert +1}.$ Then

\begin{equation} \label{5.3.7aa}
\sigma _{_{n}}f=\frac{1}{n}\sum_{j=1}^{2^{\left\vert \alpha _{k}\right\vert
}}S_{j}f+\frac{1}{n}\sum_{j=2^{\left\vert \alpha _{k}\right\vert
	}+1}^{n}S_{j}f=III+IV.
\end{equation}

It is evident that
\begin{equation} \label{5.3.7aaa}
S_{j}f=0,\,\ \text{ if \thinspace \thinspace }0\leq j\leq 2^{\left\vert \alpha
	_{1}\right\vert }
\end{equation}

Let $2^{\left\vert \alpha _{s}\right\vert }<j\leq 2^{\left\vert
	\alpha _{s}\right\vert +1},$ where $s=1,2,...,k.$ If we apply (\ref{3.3.11AA}) we get that
\begin{eqnarray} \label{5.3.8aa}
S_{j}f &=&\sum_{\eta =0}^{s-1}\Phi ^{1/2p}\left( 2^{\left\vert \alpha _{\eta}\right\vert +1}\right) \left( D_{2^{\left\vert \alpha _{\eta }\right\vert+1}}-D_{2^{\left\vert \alpha _{\eta }\right\vert }}\right)  \\
&&+\Phi ^{1/2p}\left( 2^{\left\vert \alpha _{s}\right\vert +1}\right) w_{2^{\left\vert \alpha _{s}\right\vert }}D_{j-2^{\left\vert \alpha _{s}\right\vert }}.  \notag
\end{eqnarray}

Let $2^{\left\vert \alpha _{s}\right\vert +1}\leq j\leq 2^{\left\vert \alpha_{s+1}\right\vert },$ $s=0,1,...k-1.$ Then if we use (\ref{3.3.12AA}) we can conclude that
\begin{equation} \label{5.3.10aaaa}
S_{j}f=\sum_{\eta =0}^{s}\Phi ^{1/2p}\left( 2^{\left\vert \alpha _{\eta
	}\right\vert +1}\right) \left( D_{2^{\left\vert \alpha _{\eta }\right\vert
		+1}}-D_{2^{\left\vert \alpha _{\eta }\right\vert }}\right).
\end{equation}

Let $x\in I_{2}\left( e_{0}+e_{1}\right).$ Since (see Lemmas \ref{lemma1} and \ref{lemma4})
\begin{equation} \label{5.3.40}
D_{2^{n}}\left( x\right) =K_{2^{n}}\left( x\right) =0,\text{ where }n\geq 2
\end{equation}
by combining (\ref{5.3.122}) and (\ref{5.3.7aaa}-\ref{5.3.40}) we get that
\begin{eqnarray}  \label{5.3.9aaa}
III&=&\frac{1}{n}\sum_{\eta =0}^{k-1}\Phi ^{1/2p}\left( 2^{\left\vert \alpha
	_{\eta }\right\vert +1}\right) \sum_{v=2^{\left\vert \alpha _{\eta
		}\right\vert }+1}^{2^{\left\vert \alpha _{\eta }\right\vert +1}}D_{v}\left(
x\right) \\
&=&\frac{1}{n}\sum_{\eta =0}^{k-1}\Phi ^{1/2p}\left( 2^{\left\vert \alpha
	_{\eta }\right\vert +1}\right) \left( 2^{\left\vert \alpha _{\eta
	}\right\vert +1}K_{2^{\left\vert \alpha _{\eta }\right\vert +1}}\left(
x\right) -2^{\left\vert \alpha _{\eta }\right\vert }K_{2^{\left\vert \alpha
		_{\eta }\right\vert }}\left( x\right) \right) =0. \notag
\end{eqnarray}

If we use (\ref{5.3.8aa}) when $s=k$  for $IV$ we can write that
\begin{eqnarray} \label{5.3.9aa}
IV &=&\frac{n-2^{\left\vert \alpha _{k}\right\vert }}{n}\sum_{\eta
	=0}^{k-1}\Phi ^{1/2p}\left( 2^{\left\vert \alpha _{\eta }\right\vert
	+1}\right) \left( D_{2^{\left\vert \alpha _{\eta }\right\vert
		+1}}-D_{2^{\left\vert \alpha _{\eta }\right\vert }}\right)   \\ \notag
&+&\frac{\Phi ^{1/2p}\left( 2^{\left\vert \alpha _{k}\right\vert +1}\right)
}{n}\sum_{j=2^{_{\left\vert \alpha _{k}\right\vert }}+1}^{n}w_{2^{\left\vert \alpha _{k}\right\vert }}
D_{j-2^{\left\vert \alpha _{k}\right\vert }}=IV_{1}+IV_{2}.
\notag
\end{eqnarray}

By combining (\ref{5.3.122}) and (\ref{5.3.40}) we can conclude that

\begin{equation}
IV_{1}=0,\text{ \ where \ }x\in I_{2}\left( e_{0}+e_{1}\right) .  \label{8aaaa}
\end{equation}

Let $\alpha _{k}\in \mathbb{A}_{0,2},$ $2^{\left\vert \alpha _{k}\right\vert
}<n<2^{\left\vert \alpha _{k}\right\vert +1}$ and $x\in I_{2}\left(
e_{0}+e_{1}\right) $. Since $n-2^{_{\left\vert \alpha _{k}\right\vert }}\in
\mathbb{A}_{0,2},$ Lemmas \ref{lemma0} and \ref{lemma3} and (\ref{5.3.40}) follows that
\begin{eqnarray}  \label{5.3.31}
\left\vert IV_{2}\right\vert
&=&\frac{\Phi ^{1/2p}\left( 2^{\left\vert \alpha _{k}\right\vert +1}\right)
}{n}\left\vert \sum_{j=1}^{n-2^{\left\vert _{\alpha _{k}}\right\vert
}}D_{_{j}}\left( x\right) \right\vert  \notag\\
&=&\frac{\Phi ^{1/2p}\left( 2^{\left\vert \alpha _{k}\right\vert +1}\right)
}{n}\left\vert \left( n-2^{_{\left\vert \alpha _{k}\right\vert }}\right)
K_{n-2^{_{\left\vert \alpha _{k}\right\vert }}}\left( x\right) \right\vert\geq\frac{\Phi ^{1/2p}\left( 2^{\left\vert \alpha _{k}\right\vert
		+1}\right) }{2^{\left\vert \alpha _{k}\right\vert +1}}. \notag
\end{eqnarray}

Let $0<p<1/2$ and $n\in \mathbb{A}_{0,2}.$ By combining (\ref{5.3.7aa}-\ref{5.3.31}) we get that

\begin{eqnarray} \label{10aaa}
&&\left\Vert \sigma _{n}f\right\Vert _{weak-L_{p}(G)}^{p}   \\
&\geq &\frac{c_{p}\Phi ^{1/2}\left( 2^{\left\vert \alpha _{k}\right\vert
		+1}\right) }{2^{p\left( \left\vert \alpha _{k}\right\vert +1\right) }}\mu
\left\{ x\in I_{2}\left( e_{0}+e_{1}\right) :\text{ }\left\vert \sigma
_{n}f\right\vert \geq \frac{c_{p}\Phi ^{1/2p}\left( 2^{\left\vert \alpha
		_{k}\right\vert +1}\right) }{2^{\left\vert \alpha _{k}\right\vert +1}}%
\right\}   \notag \\
&\geq &\frac{c_{p}\Phi ^{1/2}\left( 2^{\left\vert \alpha _{k}\right\vert
		+1}\right) }{2^{p\left( \left\vert \alpha _{k}\right\vert +1\right) }}\mu
\left\{ I_{2}\left( e_{0}+e_{1}\right) \right\}\geq  \frac{c_{p}\Phi
	^{1/2}\left( 2^{\left\vert \alpha _{k}\right\vert +1}\right) }{2^{p\left(
		\left\vert \alpha _{k}\right\vert +1\right) }}.  \notag
\end{eqnarray}

Hence,
\begin{eqnarray*}
	\underset{n=1}{\overset{\infty }{\sum }}\frac{\left\Vert \sigma
		_{n}f\right\Vert _{weak-L_{p}(G)}^{p}}{\Phi \left( n\right) }
	&\geq & \underset{\left\{ n\in \mathbb{A}_{0,2}:\text{ }2^{\left\vert \alpha _{k}\right\vert}<n<2^{\left\vert \alpha _{k}\right\vert +1}\right\} }{\sum }\frac{\left\Vert \sigma _{n}f\right\Vert _{weak-L_{p}(G)}^{p}}{\Phi \left(n\right) } \\
	&\geq &\frac{1}{\Phi ^{1/2}\left( 2^{\left\vert \alpha _{k}\right\vert
			+1}\right) }\underset{\left\{ n\in \mathbb{A}_{0,2}:\text{ }2^{\left\vert
			\alpha _{k}\right\vert }<n<2^{\left\vert \alpha _{k}\right\vert +1}\right\} }%
	{\sum }\frac{1}{2^{p\left( \left\vert \alpha _{k}\right\vert +1\right) }} \\
	&\geq &\frac{c_{p}2^{\left( 1-p\right) \left( \left\vert \alpha
			_{k}\right\vert +1\right) }}{\Phi ^{1/2}\left( 2^{\left\vert \alpha
			_{k}\right\vert +1}\right) }\rightarrow \infty ,\text{ \qquad as \ \ }%
	k\rightarrow \infty .
\end{eqnarray*}

The proof of Theorem \ref{theorem5.3.1} is complete.

\begin{theorem} \label{theorem5.3.2}
	Let $f\in H_{1/2}(G).$ Then
	\begin{equation*}
	\underset{n\in \mathbf{%
			\mathbb{N}
		}_{+}}{\sup }\underset{\left\Vert f\right\Vert _{H_{p}}\leq 1}{\sup }\frac{1%
	}{n}\underset{m=1}{\overset{n}{\sum }}\left\Vert \sigma _{m}f\right\Vert
	_{1/2}^{1/2}=\infty .
	\end{equation*}
\end{theorem}

{\bf Proof}:
Let $0<p\leq 1$ and
\begin{equation*}
f_{k}({x}):=2^{k}\left( D_{2^{k+1}}(x)-D_{2^k}(x)\right)
\end{equation*}%
Since
\begin{equation*}
\text{supp}(f_{k})=I_{k},\text{ \ \ }\int_{I_{k}}a_{k}d\mu =0
\end{equation*}
and
\begin{equation*}
\left\Vert f_{k}\right\Vert _{\infty }\leq  2^{2k}=(\text{supp}f_{k})^{-2},
\end{equation*}%
we conclude that $f_{k}$ is $ 1/2 $-atom, for every $k\in \mathbb{N}$

Moreover, if we use orthogonality of Walsh functions we get that 
\begin{equation*}
S_{2^{n}}\left(f_{k},x)\right)
\end{equation*}
\begin{equation*}
=\left\{
\begin{array}{ll}
0, & n=0,...,k, \\
\left( D_{2^{k+1}}(x)-D_{2^k}(x)\right), & n\geq k+1,
\end{array}%
\right.
\end{equation*}%
and
\begin{eqnarray*}
	&&\sup\limits_{n\in \mathbb{N}}\left\vert S_{2^{n}}\left(f_{k},x\right)\right\vert \\
	&=&\left\vert \left( D_{2^{k+1}}(x)-D_{2^k}(x)\right) \right\vert,
\end{eqnarray*}
where ${x}\in G$.

By combining first equality of Lemma \ref{lemma0} and Lemma \ref{lemma3.2.3} we obtain that
\begin{eqnarray*}
	\left\Vert a_{k}\right\Vert _{H_{p}(G)} 
	&=&{2^{k}}
	\left\Vert \sup\limits_{n\in \mathbb{N}}\left\vert S_{2^{n}}\left( D_{2^{k+1}}(x)-D_{2^k}(x)\right)\right\vert\right\Vert _{1/2} \\
	&=&{2^{k}}\left\Vert\left( D_{2^{k+1}}(x)-D_{2^k}(x)\right)\right\Vert _{1/2} ={2^{k}}\left\Vert D_{2^k}(x)\right\Vert _{1/2}\leq  {2^{k}}\cdot 2^{-k}\leq 1.
\end{eqnarray*}

It is easy to easy to show that
\begin{equation} \label{5.3.14a}
\widehat{f}_{m}\left( i\right) =\left\{
\begin{array}{l}
\text{ }2^{m},\text{ if }i=2^{m},...,2^{m+1}-1, \\
\text{ }0,\text{otherwise}%
\end{array}%
\right.
\end{equation}
and
\begin{equation}  \label{5.2.14}
S_{i}f_{m}=\left\{
\begin{array}{l}
2^{m}\left( D_{i}-D_{2^{m}}\right) ,\text{ \ if }i=2^{m}+1,...,2^{m+1}-1, \\
\text{ }f_{m},\text{ \ if }i\geq 2^{m+1}, \\
0,\text{ \ otherwise.}
\end{array}%
\right.
\end{equation}

Let $0<n<2^{m}.$ By using first equality of Lemma \ref{lemma0} we have that

\begin{eqnarray} \label{5.3.16b}
\left\vert \sigma _{n+2^{m}}f_{m}\right\vert &=&\frac{1}{n+2^{m}}\left\vert
\overset{n+2^{m}}{\underset{j=2^{m}+1}{\sum }}S_{j}f_{m}\right\vert  \notag\\
&=&\frac{1}{n+2^{m}}\left\vert 2^{m}\overset{n+2^{m}}{\underset{j=2^{m}+1}{%
		\sum }}\left( D_{j}-D_{2^{m}}\right) \right\vert  \notag \\
&=&\frac{1}{n+2^{m}}\left\vert 2^{m}\overset{n}{\underset{j=1}{\sum }}\left(
D_{j+2^{m}}-D_{2^{m}}\right) \right\vert  \notag \\
&=&\frac{1}{n+2^{m}}\left\vert 2^{m}\overset{n}{\underset{j=1}{\sum }}%
D_{j}\right\vert=\frac{2^{m}}{n+2^{m}}n\left\vert K_{n}\right\vert .  \notag
\end{eqnarray}

Let
$$n=\sum_{i=1}^{s}\sum_{k=l_{i}}^{m_{i}}2^{k},$$
where
$$0\leq l_{1}\leq m_{1}\leq l_{2}-2<l_{2}\leq m_{2}\leq ...\leq l_{s}-2<l_{s}\leq m_{s}.$$
By applying Lemma \ref{lemma7} and (\ref{5.3.16b}) we find that
\begin{equation*}
\left\vert\sigma_{n+2^{m}}f_{m}\left( x\right)\right\vert\geq c2^{2l_{i}}, \text{ \ \ \ where \ \ \ }x\in I_{l_{i}+1}\left(e_{l_{i}-1}+e_{l_{i}}\right).
\end{equation*}
Hence,
\begin{eqnarray*}
	&&\int_{G}\left\vert \sigma _{n+2^{m}}f_{m}(x)\right\vert ^{1/2}d\mu \left(
	x\right) \\
	&\geq &\text{ }\underset{i=0}{\overset{s}{\sum }}\int_{I_{l_{i}+1}\left(
		e_{l_{i}-1}+e_{l_{i}}\right) }\left\vert \sigma
	_{n+2^{m}}f_{m}(x)\right\vert ^{1/2}d\mu \left( x\right) \\
	&\geq &c\overset{s}{\underset{i=0}{\sum }}\frac{1}{2^{l_{i}}}2^{l_{i}}\geq cs\geq cV\left( n\right) .
\end{eqnarray*}

According to the second estimation of Lemma \ref{lemma2} we can conclude that
\begin{eqnarray*}
	\underset{n\in \mathbf{\mathbb{N}}_{+}}{\sup }\underset{\left\Vert f\right\Vert _{H_{p}}\leq 1}{\sup}\frac{1}{n}\underset{k=1}{\overset{n}{\sum }}\left\Vert \sigma _{k}f\right\Vert_{1/2}^{1/2}
	&\geq &\frac{1}{2^{m+1}}\overset{2^{m+1}-1}{\underset{k=2^{m}+1}{\sum }}
	\left\Vert \sigma _{k}f_{m}\right\Vert _{1/2}^{1/2} \\
	&\geq &\frac{c}{2^{m+1}}\underset{k=2^{m}+1}{\overset{2^{m+1}-1}{\sum }}
	V\left( k-2^{m}\right) \\
	&\geq& \frac{c}{2^{m+1}}\underset{k=1}{\overset{2^{m}-1}{\sum }}V\left( k\right)\geq c\log m\rightarrow \infty ,\text{ \ as \ }m\rightarrow \infty .
\end{eqnarray*}

The proof is complete.

\end{document}